\RequirePackage{ifpdf}
\ifpdf 
\documentclass[pdftex]{sigma}
\else
\documentclass{sigma}
\fi

\usepackage{ulem}

\def\one{\mbox{\rm 1}\hskip-2.8pt \mbox{\rm l}}

\newcommand{\nco}{\newcommand}
\nco{\munite}{\ensuremath{\,\,\mathrm{l}\!\!\!1}}
\nco{\ZZ}{\mathbb{Z}}
\nco{\CC}{\mathbb{C}}

\begin{document}

\allowdisplaybreaks

\renewcommand{\thefootnote}{$\star$}

\renewcommand{\PaperNumber}{044}

\FirstPageHeading

\ShortArticleName{Quantum Symmetries for Exceptional ${\rm SU}(4)$ Modular Invariants}

\ArticleName{Quantum Symmetries for Exceptional $\boldsymbol{{\rm SU}(4)}$ Modular Invariants Associated with Conformal Embeddings}

\Author{Robert COQUEREAUX and Gil SCHIEBER}

\AuthorNameForHeading{R.~Coquereaux and G.~Schieber}
\Address{Centre de Physique Th\'eorique (CPT), Luminy, Marseille, France\footnote{UMR 6207 du CNRS et des Universit\'es  Aix-Marseille I, Aix-Marseille II,
et du Sud Toulon-Var,  af\/f\/ili\'e \`a la FRUMAM (FR 2291)}}


\Email{\href{mailto:coque@cpt.univ-mrs.fr}{coque@cpt.univ-mrs.fr}, \href{mailto:schieber@cpt.univ-mrs.fr}{schieber@cpt.univ-mrs.fr}}

\ArticleDates{Received December 24, 2008, in f\/inal form March 31,
2009; Published online April 12, 2009}

\Abstract{Three exceptional modular invariants of ${\rm SU}(4)$ exist at levels $4$, $6$ and $8$.  They can be obtained from appropriate conformal embeddings and the corresponding graphs have self-fusion. From these embeddings, or from their associated modular invariants, we determine  the algebras of quantum symmetries, obtain their generators, and, as a by-product, recover the known graphs ${\mathcal E}_4 $, ${\mathcal E}_6 $ and ${\mathcal E}_8 $ describing  exceptional quantum subgroups of type~${\rm SU}(4)$.  We  also obtain characteristic numbers (quantum cardinalities, dimensions) for each of them and for their associated quantum groupo\"\i ds.}


\Classification{81R50; 16W30; 18D10}

\renewcommand{\thefootnote}{\arabic{footnote}}
\setcounter{footnote}{0}

\pdfbookmark[1]{Foreword}{Foreword}
\section*{Foreword}

\noindent{\bf General presentation.}
A classif\/ication of ${\rm SU}(4)$  graphs associated with WZW models,  or  ``quantum graphs'' for short, was presented by A. Ocneanu in \cite{Ocneanu:Bariloche} and claimed to be completed. These graphs generalize the ADE Dynkin diagrams that classify the ${\rm SU}(2)$ models \cite{CIZ},  and the Di Francesco--Zuber diagrams that classify the ${\rm SU}(3)$ models \cite{DiFrancescoZuber}. They describe modules over a ring of irreducible representations of  quantum ${\rm SU}(4)$ at roots of unity.
A particular  partition function associated with each of those quantum graphs is modular invariant.

According to  \cite{Ocneanu:Bariloche}, the ${\rm SU}(4)$ family includes  the ${\mathcal A}_k$ series (describing fusion algebras) and their conjugates for all $k$, two kinds of orbifolds, the ${\mathcal D}^{(2)}_{k} = \mathcal{A}_k/2$ series for all $k$ (with self-fusion when $k$ is even), and members of the ${\mathcal D}^{(4)}_{k} = \mathcal{A}_k/4$ series when $k$ is even  (with self-fusion when $k$ is divisible by 8), together with their conjugates. The orbifolds are constructed by using the $Z_4$ action on weigths generated by  ${\epsilon \{\lambda_1, \lambda_2, \lambda_3 \} = \{ k- \lambda_1 - \lambda_2 - \lambda_3, \lambda_1, \lambda_2\}}$ or the $Z_2$ action generated by $\epsilon^2$.
The ${\rm SU}(4)$ family also includes an exceptional case,  ${\mathcal{D}^{(4)t}_{8}}$,  without self-fusion (a generalization of $E_7$), and three exceptional quantum graphs with self-fusion, at levels $4$, $6$ and $8$, denoted ${\mathcal E}_4 $, ${\mathcal E}_6 $ and ${\mathcal E}_8 $, together with one exceptional module for each of the last two. The  modular invariant partition functions associated with ${\mathcal E}_4 $, ${\mathcal E}_6 $ and ${\mathcal E}_8 $ can be obtained from appropriate conformal embeddings, namely from ${\rm SU}(4)$ level $4$ in ${\rm Spin}(15)$,  from ${\rm SU}(4)$ level~$6$ in~${\rm SU}(10)$,  and from ${\rm SU}(4)$ level 8 in ${\rm Spin}(20)$.
There exists also a conformal embedding of~${\rm SU}(4)$, at level~2, in  ${\rm SU}(6)$, but this gives rise to $\mathcal{D}^{(2)}_2= \mathcal{A}_2/2$, the f\/irst member of the ${\mathcal D}^{(2)}_{k}$ series.   This exhausts the list of conformal embeddings of ${\rm SU}(4)$. The other ${\rm SU}(4)$ quantum graphs, besides the ${\mathcal A}_k$,  can either be obtained as modules over the exceptional ones,  or are associated (possibly using conjugations) with non-simple conformal embeddings followed by contraction, ${\rm SU}(4)$ appearing only as a direct summand of the embedded algebra.

Vertices $a,b,\ldots$ of a chosen quantum graph (denoted generically by ${\mathcal E}_k$)  describe boundary conditions for a WZW conformal f\/ield theory specif\/ied by ${\rm SU}(4)_k$. These irreducible objects span a vector space which is a module over the fusion algebra, itself spanned, as a linear space,  by the vertices $m,n,\ldots$ of the graph ${\mathcal A}_k({\rm SU}(4))$, or ${\mathcal A}_k$ for short since ${\rm SU}(4)$ is chosen once and for all, the truncated Weyl chamber at level~$k$ (a Weyl alcove). Vertices of ${\mathcal A}_k$ can be understood as integrable irreducible highest weight representations of the af\/f\/ine Lie algebra
$\widehat{su}(4)$ at level~$k$ or as irreducible representations with non-zero $q$-trace of the quantum group ${\rm SU}(4)_q$ at the root of unity $q=\exp(i \pi/(k+g))$, $g$ being the dual Coxeter number (for ${\rm SU}(4)$, $g=4$). Edges of  ${\mathcal E}_k$ describe action of the fundamental representations of ${\rm SU}(4)$, the generators of  ${\mathcal A}_k$.

 To every quantum graph one associates an algebra of quantum symmetries\footnote{Sometimes called ``fusion algebra of defect lines'' or ``full system''.}  ${\mathcal O}$, along the lines described in~\cite{Ocneanu:paths}. Its vertices $x, y, \ldots $ can be understood, in the interpretation of~\cite{PetkovaZuber:Oc}, as  describing  the same BCFT theory but with defects labelled by $x$.
 To every fundamental representation ($3$~of them for ${\rm SU}(4)$) one associates two generators of ${\mathcal O}$, respectively called ``left'' and ``right''. Multiplication by these generators is described by a graph\footnote{One should not confuse the quantum graph (or McKay graph) that refers to ${\mathcal E}_k$ with the graph of quantum symmetries (or Ocneanu graph) that refers to ${\mathcal O}.$}, called the Ocneanu graph. Its vertices span the algebra ${\mathcal O}$ as a linear space, and its edges describe multiplication by the left and right fundamental generators (we have $6 = 2 \times 3$ types of edges\footnote{Actually, since those associated with weights $\{100\}$ and $\{001\}$ are conjugated and $\{010\}$ is real (self-conjugated), we need only $2$ types of edges (the f\/irst is oriented, the other is not) for ${\mathcal E}_k$ or ${\mathcal A}_k$, and $4 = 2 \times 2$ types of edges for ${\mathcal O}$.} for ${\rm SU}(4)$).

The exceptional modular invariants at level $4$ and $6$ were found by \cite{SchellekensYankielowicz, AltschulerBauerItzykson}, and at level~$8$ by~\cite{AldazabalEtAl}. The corresponding quantum graphs\footnote{We often drop the reference to ${\rm SU}(4)$ since no confusion may arise: we are not discussing in this paper the usual $E_6=\mathcal{E}_{10}({\rm SU}(2))$  or  $E_8=\mathcal{E}_{28}({\rm SU}(2))$  Dynkin diagrams!} ${\mathcal E}_4$,  ${\mathcal E}_6$ and ${\mathcal E}_8$ were respectively obtained by \cite{PetkovaZuberNP1995, PetkovaZuberNP1997, Ocneanu:Bariloche}.
There are several techniques to determine quantum graphs. One of them, probably the most powerful but involving rather heavy calculations, is to obtain the quantum graph associated\footnote{We use the word ``associated'' here in a rather loose sense, since the relation between both concepts is not one-to-one.} with a modular invariant as a by-product of the determination of its algebra of quantum symmetries. This requires in particular the solution of the so-called modular splitting equation, which is
a huge collection of equations between matrices with non-negative integral entries, involving the known fusion algebra, the chosen modular invariant, and expressing the fact that ${\mathcal O}$ is a bi-module over ${\mathcal A}_k$.
Because of the heaviness of the calculation, a simplif\/ied method using only the f\/irst line of the modular invariant matrix was used  in \cite{Ocneanu:Bariloche} to achieve this goal, namely the determination of ${\rm SU}(4)$ quantum graphs (some of them, already mentioned, were already known) but the algebra of quantum symmetries was not obtained in all cases. To our knowledge, for exceptional modular invariants of ${\rm SU}(4)$ at levels $4$, $6$ and $8$,  the full modular splitting system had not been solved, the full torus structure had not been obtained, and the graph of quantum symmetries was not known.  This is what we did. We have recovered in particular the structure of the already known quantum graphs;  they now appear, together with their modules, as components of their respective Ocneanu graphs.

{\bf Categorical description.}
Category theory of\/fers a synthetic presentation of the whole subject and we present it here in a few lines,  for the benef\/ice of those readers who may f\/ind appealing such a description. However, it will not be used in the body of our article.
The starting point is the fusion category ${\mathcal A}_k$ associated with a Lie group $K$. This modular category, both monoidal and ribbon, can be def\/ined either in terms of representation theory of an af\/f\/ine Lie algebra (simple objects are highest weight integrable irreducible representations), or in terms of representation theory of a quantum group at roots of unity (simple objects are irreducible representations of non-vanishing quantum dimension).  In the case of ${\rm SU}(2)$, we refer to the description given in~\cite{Ostrik,EtingofOstrik}.
One should keep in mind the distinction between this category (with its objects and morphisms),  its Grothendieck ring (the fusion ring), and the graph describing multiplication by  its generators,  but they are denoted by the same symbol.
The next ingredient is an additive category ${\mathcal E}_k$, not modular usually, on which the previous one, ${\mathcal A}_k$,  acts.
In general this module-category  ${\mathcal E}_k$ has no-self-fusion (no compatible monoidal structure) but in the cases studied in the present paper, it does. Again, the category itself, its Grothendieck group, and the graph (here called McKay graph) describing the action of generators of  ${\mathcal A}_k$ are denoted by the same symbol. The last ingredient is the centralizer (or dual) category ${\mathcal O}={\mathcal O}({\mathcal E}_k)$ of ${\mathcal E}_k$ with respect to the action of ${\mathcal A}_k$.  It is monoidal and comes with its own ring (the algebra of quantum symmetries) and graph (the Ocneanu graph). One way to obtain a realization of this collection of data is to construct a f\/inite dimensional weak bialgebra ${\mathcal B}$, which should be such that  ${\mathcal A}_k$ can be realized as ${\rm Rep}({\mathcal B})$, and also such that ${\mathcal O}$ can be realized as ${\rm Rep}({\widehat{\mathcal B}})$,  where $\widehat{\mathcal B}$ is the dual of~${\mathcal B}$. These two algebras are f\/inite dimensional, actually semisimple in our case,  and one algebra structure (say $\widehat{\mathcal B}$) can be traded against a coalgebra structure on its dual. ${\mathcal B}$ is a~weak bialgebra, not a bialgebra, because $\Delta \one \neq \one \otimes \one$, the coproduct in  ${\mathcal B}$ being $\Delta$, and $\one$ its unit.  ${\mathcal B}$  is not only a weak bialgebra but a weak Hopf algebra: one can def\/ine an antipode, with the expected properties.

\begin{remark} \label{remark}
Given a graph def\/ining a module over a fusion ring ${\mathcal A}_k$ for some Lie group $K$,  the question is to know if it is a ``good graph'', i.e., if the corresponding module-category indeed exists.
Using A.~Ocneanu's terminology~\cite{Ocneanu:Bariloche}, this will be the case if and only if  one can associate, in a coherent manner,  a complex number to each triangle of the graph (when the rank of $K$ is $\geq 2$): this def\/ines, up to some kind of gauge choice,  {\it a self-connection on the set of triangular cells}. Here,  ``coherent manner'' means that there are two compatibility  equations, respectively nicknamed {\it the small and  large   pocket equations}, that this self-connection should obey.
These equations ref\/lect properties that hold for the intertwining operators of a fusion category,  they are sometimes called {``compatibility equations for Kuperberg spiders''} (see \cite{Kuperberg:spiders}).  The point is that exhibiting a module over a fusion ring does not necessarily entail existence of an underlying theory: when the graph (describing the module structure) does not admit any self-connection in the above sense, it should be rejected;
another way to express the same thing is to say that a~particular family of $6j$ symbols, expected to obey appropriate equations, fails to be found.
Such features are not going to be discussed further in the present paper.
\end{remark}

{\bf Historical remarks concerning $\boldsymbol{{\mathcal E}_8({\rm SU}(4))}$.} Not all
conformal embeddings $K\subset G$ correspond to isotropy-irreducible
pairs and not all isotropy-irreducible homogeneous spaces def\/ine
conformal embeddings. However, it is a known fact that most isotropy
irreducible spaces $G/K$ (given in \cite{Wolf}) indeed def\/ine conformal
embeddings. This is actually so in all examples stu\-died here, and in particular
for the ${\rm SU}(4) \subset {\rm Spin}(20)$ case
which can also be recognized as the smallest member ($n=1)$ of a
$D_{2n+1} \subset D_{(n+1)(4n+1)} $ family of conformal embeddings
appearing on table $4$ of the standard reference~\cite{BaisBouwknegt},
and on table II(a) of the standard reference~\cite{SchellekensWarner},  since ${\rm SU}(4) \simeq  {\rm Spin}(6)$.
This embedding, which is ``special''
(i.e., non  regular: unequal ranks and Dynkin index not equal to $1$),  does not seem to be quoted in other
standard references on conformal embeddings (for instance
\cite{KacWakimoto, LevsteinLiberati, Verstegen}),
although it is explicitly mentioned in~\cite{AldazabalEtAl}  and although its rank-level dual is
indirectly used in case~18 of~\cite{Schellekens:c24}, or in~\cite{Xu:mirror}. The corresponding
${\rm SU}(4)$ modular invariant was later recovered by
\cite{Ocneanu:Bariloche}, using arithmetical methods, and used to
determine the ${\mathcal E}_8({\rm SU}(4))$ quantum graph, but since the existence
of an associated conformal embedding had slipped into oblivion, it was
incorrectly stated that this particular example could not be obtained
from conformal embedding considerations.

{\bf Structure of the article.}
The technique relating modular invariants to conformal embeddings is standard~\cite{YellowBook} but the results concerning ${\rm SU}(4)$ are either scattered in the literature, or unpublished; for this reason, we devote the main part of the f\/irst section to it.
In the same section, we obtain characteristic numbers (quantum cardinality, quantum dimensions etc.) for the ${\mathcal E}_k$ graphs.  In the second section, after a  description of the structures at hand and a general presentation of our
method of resolution, we solve, in a f\/irst step, for the three exceptional cases ${\mathcal E}_4$, ${\mathcal E}_6$ and ${\mathcal E}_8$ of the ${\rm SU}(4)$ family,  the full modular splitting equation that determines the corresponding set of toric matrices (generalized partition functions)  and, in a second step,  the general intertwining equations that determine the structure of the generators of the algebra of quantum symmetries.
The size of calculations involved in this part is huge (quite intensive computer help was required) and, for reasons of size,  we can only present part of our results. On the other hand, each case being exceptional, there are no generic formulae.
For each case, we encode the structure of the algebra of quantum symmetries by displaying the Cayley graphs of multiplication by the fundamental generators, whose collection makes the Ocneanu graph. We also give a brief description of the structure of the corresponding quantum groupo\"\i ds.  In the appendices, after a short description of the Kac--Peterson formula,  we gather several explicit results, providing quantum dimensions for those irreducible representations of the various groups used in the text.

The interested reader may also consult the article \cite{RobertGil:E4SU4} which provides more information on the general theory and gives a more complete description of the $\mathcal{E}_4({\rm SU}(4))$ case.
Properties of quantum graphs of type ${\rm SU}(3)$ and their quantum symmetries are summarized in~\cite{RobertGilSL3Categories}, see also~\cite{GilDahmaneHassan} and references therein. Those of type ${\rm SU}(2)$ are certainly well known but many explicit results, like the explicit structure of toric matrices for exceptional diagrams,  can be found in~\cite{GilCoque:ADE}.

\section[Conformal embeddings of ${\rm SU}(4)$]{Conformal embeddings of $\boldsymbol{{\rm SU}(4)}$} \label{section1}

\subsection[Homogeneous spaces $G/K$]{Homogeneous spaces $\boldsymbol{G/K}$}

We describe the embeddings of $K={\rm SU}(4)$ in $G={\rm Spin}(15), \, {\rm SU}(10), \, {\rm Spin}(20)$.
The reduction of the adjoint representation of $G$ with respect to $K$ reads ${\rm Lie}(G) = {\rm Lie}(K) \oplus T(G/K)$.
The isotropy representation of $K$ on the tangent space at the origin of $G/K$ has dimension $\dim(G)-\dim(K)$.
In all three cases, the space $G/K$ is isotropy irreducible (but not symmetric): the isotropy representation is real irreducible. After extension to the f\/ield of complex numbers it may stay irreducible (strong irreducibility) or not.
The following are known results, already mentioned in~\cite{Wolf}.

\begin{itemize}\itemsep=0pt
\item[] ${\rm SU}(4) \subset {\rm SU}(6)$. This embedding leads to the lowest member of an orbifold series (the ${\mathcal D}_2^{(2)}={\mathcal A}_2/2$ graph) and, in this paper, we are not interested in it.

\item[] ${\rm SU}(4) \subset {\rm Spin}(15)$.
Reduction of the adjoint representation of $G$ with respect to $K$ reads $[105] \mapsto [15] + [90] $. After complexif\/ication,  $[90]$ is recognized as the reducible representation with highest weight  $\{0,1,2\} \oplus \{2,1,0\} = [45]\oplus[\overline{45}] $ so that $G/K$ is not strongly  irreducible.

\item[] ${\rm SU}(4) \subset {\rm SU}(10)$.
Reduction of the adjoint representation of $G$ with respect to $K$ reads $[99] \mapsto [15] + [84] $. After complexif\/ication,  $[84]$ is recognized as the irreducible representation with highest weight  $\{2,0,2\}$ so that $G/K$ is strongly  irreducible.

\item[] ${\rm SU}(4) \subset {\rm Spin}(20)$.
Reduction\footnote{The inclusion  ${\rm SU}(4)/Z_4 \subset {\rm SO}(20) \subset {\rm GL}(20, \CC)$ is associated with a representation of ${\rm SU}(4)$, of dimension~$20$, with highest weight $\{0,2,0\}$.} of the adjoint representation of $G$ with respect to $K$ reads $[190] \mapsto [15] + [175] $. After complexif\/ication,  $[175]$ is recognized as the irreducible representation with highest weight  $\{1,2,1\}$, so that $G/K$ is strongly  irreducible.
\end{itemize}

\subsection{The Dynkin index  of the embeddings} \label{cartan}

 The Dynkin index $k$ of an embedding $K \subset G$ def\/ined by a branching rule $\mu \mapsto \sum_j \alpha_j \nu_j$, where  $\mu$ refers to the adjoint representation of $G$ (one of the  $\nu_j$ on the right hand side is the adjoint representation of $K$),  $ \alpha_j $ being multiplicities,  is obtained in terms of the quadratic Dynkin indices~$I_\mu$,~$I_{\nu_j}$ of the representations:
\[k = \sum_j  \alpha_j     I_{\nu_j} / I_\mu \qquad  \hbox{with} \quad I_\lambda = \frac{\dim(\lambda)}{2 \dim(K)}    \langle \lambda, \lambda + 2 \rho \rangle.  \]
Here $\rho$ is the Weyl vector and  $ \langle \;  , \;  \rangle $ is def\/ined by the fundamental quadratic form.
For the three embeddings of ${\rm SU}(4)$ that we consider, into ${\rm Spin}(15)$, ${\rm SU}(10)$ and ${\rm Spin}(20)$,  one f\/inds respectively $k=4,6,8$.

\subsection{Those embeddings are conformal}

An embedding $K \subset G$, for which the Dynkin index is $k$,  is conformal  if the following equality is satisf\/ied\footnote{Warning: it is not dif\/f\/icult to f\/ind  embeddings $K \subset G$, and appropriate values of $k$  for which this equality is satisf\/ied,  but where $k$ is {\sl not} the Dynkin index! Such embeddings are, of course, non conformal.}:
\[
 \frac{\dim(K) \times  k }{ k + g_K} = \frac{\dim(G) \times 1}{ 1 + g_G},
\]
 where $g_K$ and $g_G$ are the dual Coxeter numbers of $K$ and $G$.  One denotes by  $c$ the common value of these two expressions.
 In the framework of af\/f\/ine Lie algebras,  $c$ is interpreted as a~central charge and the numbers $k$ and $1$ denote the respective levels for the af\/f\/ine algebras corresponding to $K$ and $G$.
The above def\/inition, however,  does not require the framework of af\/f\/ine Lie algebras (or of quantum groups at roots of unity) to make sense.

Using $\dim(G) = 105, \, 99, \, 190$, for $G= {\rm Spin}(15), \, {\rm SU}(10), \, {\rm Spin}(20)$, $\dim(K={\rm SU}(4))=15$ and the corresponding values for the dual Coxeter numbers $g_G = 13, \, 10,  \, 18$ and $g_K=4$, we see  immediately that the above equality is obeyed, for the levels $k=4,6,8$, with central charges $c=15/2$, $c=9$, and $c=10$.

The conformal embeddings of ${\rm SU}(4)$ into ${\rm SU}(6)$,  ${\rm SU}(10)$ and ${\rm Spin}(15)$
belong respectively to the series of  embeddings of ${\rm SU}(N)$ into ${\rm SU}(N(N-1)/2)$, ${\rm SU}(N(N+1)/2$ and ${\rm Spin}(N^2-1)$, at respective levels $N-2$, $N+2$ and $N$ (provided $N$ is big enough), whereas the last one, namely~${\rm SU}(4)$ into ${\rm Spin}(20)$, is recognized as the smallest member of the ${\rm Spin}(N) \subset {\rm Spin}((2N+2)(4N+1))$ series, since ${\rm SU}(4) \simeq {\rm Spin}(6)$.

\subsection{The modular invariants}
Here we reduce the diagonal modular invariants of $G= {\rm Spin}(15), \, {\rm SU}(10), \, {\rm Spin}(20)$,  at level $k=1$, to $K={\rm SU}(4)$, at levels $k=4,6,8$, and obtain exceptional modular invariants for ${\rm SU}(4)$ at those levels.
The previous section was somehow ``classical'' whereas this one is ``quantum''.
Since there is an equivalence of categories~\cite{Drinfeld,  KazhdanLusztig}  between the fusion category (integrable highest weight representations)  of an af\/f\/ine algebra at some level and  a category of representations with non-zero $q$-dimension for the corresponding quantum group at a root of unity determined by the level, we shall freely use both terminologies. From now on, simple objects will be called i-irreps,  for short.

\subsubsection{The method}\label{orderirrep}

\begin{itemize}\itemsep=0pt
\item One has f\/irst to determine what i-irreps  $\lambda$ appear at the chosen levels.
Given a level $k$, the  integrability condition reads
$
 \langle \lambda, \theta \rangle  \leq k
$, where $\theta$ is the highest root of the chosen Lie algebra. This is the simplest way of determining these representations.
One may notice that they will have non vanishing q-dimension when $q$ is specialized to the value $q = \exp(i\pi/\kappa)$,  with $\kappa = g_G + k$ (use the quantum Weyl formula together with the property $\langle \rho, \theta \rangle = g-1$, $\rho$ being the Weyl vector,  and the fact that $\kappa_q = 0$, see footnote~\ref{qnumdef}).
When $k=1$,  the i-irreps of $G={\rm SU}(n)$ are the fundamental representations, and the trivial.
For other Lie groups $G$, not all fundamental representations give rise to i-irreps at level $1$ (see Appendix).

\item To an i-irrep $\lambda$ of $G$ or of $K$, one associates a conformal weight def\/ined by
\begin{gather}
h_{  \lambda} = \frac{ \langle \lambda, \lambda + 2 \rho \rangle}{2(k+g)},
\label{confdim}
\end{gather}
where  $g$ is the dual Coxeter number of the chosen  Lie algebra,  $k$ is the level (for $G$, one chooses $k=1$), $\rho$ is the Weyl vector (of~$G$, or of~$K$).  The scalar product is given by the inverse of the Cartan matrix when the Lie algebra is simply laced ($A_3 \simeq {\rm SU}(4)$, $A_9 \simeq {\rm SU}(10)$ or $D_{10}\simeq {\rm Spin}(20)$), and is the inverse of the matrix obtained by multiplying the last line of the Cartan matrix by a coef\/f\/icient $2$ in the non simply laced case $B_7\simeq {\rm Spin}(15)$. Note that $h_{  \lambda}$ is related to the phase $m_{\lambda}$ of the modular $t$ matrix by $m_{\lambda} = h_{\lambda} - c/24$.
One builds the list of i-irreps $\lambda$  of $G$ at level $1$ and calculate their conformal weights $h_{\lambda}$; then, one builds the list of i-irreps $\mu$ of $K$ at level $k$ and calculate their conformal weights $h_{\mu}$.

\item  A necessary -- but not suf\/f\/icient -- condition for an (af\/f\/ine or quantum) branching from~$\lambda$ to~$\mu$ is that  $h_{\mu} = h_{\lambda} + m$ for some non-negative integer $m$.
So we can make a list of candidates for the branching rules ${  \lambda} \hookrightarrow \sum_n c_n   {\mu}_n $,  where $c_n$ are positive integers to be determined.

\item  There exist several techniques to determine the coef\/f\/icients $c_n$ (some of them can be $0$), for instance using information coming from the f\/inite branching rules. An ef\/f\/icient possibility\footnote{A drawback of this method is that it may lead to several solutions (an interesting fact, however).} is
to impose that the candidate for the modular invariant matrix should commute with the generators $s$ and $t$ of ${\rm SL}(2,\ZZ)$ (modularity constraint).

\item
We write the diagonal invariant of type $G$ as a sum $\sum_s  \lambda_{\overline s} \lambda_s $. Its associated quantum graph is denoted  ${\mathcal J}={\mathcal A}_1(G)$.
Using the above branching rules, we replace, in this expression,  each $\lambda_s$ by the corresponding sum of i-irreps for $K$. The modular invariant $\mathcal{M}$ of type $K$ that we are looking for is parametrized by
\[
{\mathcal Z} = \sum_{s\in {\mathcal J}}  \bigg(\sum_{n} c_n(\overline s)  {  \mu}_n (\overline s)\bigg) \bigg(\sum_{n} c_n(s)   {  \mu}_n (s)\bigg).
\]
\end{itemize}

In all three cases we shall need to compute conformal weights for ${\rm SU}(4)$ representations.
In the base of fundamental weights\footnote{We use sometimes the same notation $\lambda_i$ to denote a representation or to denote the Dynkin labels of a weight; this should be clear from the context.
We never write explicitly the af\/f\/ine component of a weight since it is equal to $k- \langle \lambda, \theta \rangle$.}, an arbitrary weight reads $\lambda=(\lambda_n)$, the Weyl vector is $\rho = \{1,1,1\}$, the scalar product of weights is  $\langle \lambda, \mu \rangle = (\lambda_m) Q_{mn} (\mu_n)$.  At level $k$, i-irreps $\lambda = \{\lambda_1, \lambda_2, \lambda_3\}$ are such that   $0 \leq \sum\limits_{n=1}^{n=3} \lambda_n  \leq k$; they build a set of cardinality $r_A=(k+1)(k+2)(k+3)/6$.
 We order the irreducible representations $\{i,j,k\}$ of ${\rm SU}(4)$ as follows: f\/irst of all, they are ordered by increasing level $i+j+k$, then, for a given level, we set $\{i,j,k\} < \{i',j',k'\} \Leftrightarrow
i+j+k < i'+j'+k'$ or ($ i+j+k = i'+j'+k'$ and $i> i'$)  or  ($i+j+k = i'+j'+k'$, $i=i'$  and  $j > j'$).
We now consider each case, in turn.

\subsubsection[${\rm SU}(4) \subset {\rm Spin}(15)$, $k=4$]{$\boldsymbol{{\rm SU}(4) \subset {\rm Spin}(15)}$, $\boldsymbol{k=4}$}

\begin{itemize}\itemsep=0pt

\item
At level 1, there are only three i-irreps for $B_7$, namely $\{0\}$ , $\{1,0,0,0,0,0,0\}$ and $\{0,0,0,0$, $0,0,1\}$, namely the trivial, the vectorial and the spinorial.
From equation (\ref{confdim}) we calculate their conformal weights: $\left\{0,\frac{1}{2},\frac{15}{16}\right\}$.

\item    At level 4, we calculate the $35$ conformal weights for ${\rm SU}(4)$ i-irreps and f\/ind (use ordering def\/ined previously):
\[\mbox{\tiny $\displaystyle 0,\frac{15}{64},\frac{5}{16},\frac{15}{64},\frac{9}{16},\frac{39}{
   64},\frac{1}{2},\frac{3}{4},\frac{39}{64},\frac{9}{16},\frac{63}{64},1
   ,\frac{55}{64},\frac{71}{64},\frac{15}{16},\frac{55}{64},\frac{21}{16}, \frac{71}{64},1,\frac{63}{64},\frac{3}{2},\frac{95}{64},\frac{21}{16}
   ,\frac{25}{16},
   \frac{87}{64},\frac{5}{4},\frac{111}{64},\frac{3}{2},\frac{87}{64},\frac{21}{16},2,\frac{111}{64},\frac{25}{16},\frac{95}{64}
   ,\frac{3}{2}.$}
\]
\item  The dif\/ference between conformal weights of $B_7$ and $A_3$ should be an integer. This selects
the three following possibilities:
\[
\mbox{\tiny $0000000 \stackrel{?}{\hookrightarrow} 000 + 210 + 012 + 040, \qquad
  1000000  \stackrel{?}{\hookrightarrow} 101 + 400 + 121 + 004, \qquad
  0000001 \stackrel{?}{\hookrightarrow}111.$}
\]
The above three possibilities give only necessary conditions for branching. Imposing the modularity constraint implies that the multiplicity of $(111)$ should be $4$, and that all the other coef\/f\/icients indeed appear, with multiplicity $1$.  This is actually a particular case of general branching rules already found in \cite{Kac-Peter-81, KacWakimoto, Fuchs-Schellekens-Schweigert}.
 \end{itemize}

The partition function obtained from the diagonal invariant $ \vert 0000000\vert^2\! +   \vert 1000000\vert^2\! +    \vert 0000001\vert^2\!$ of $B_7$ reads:
\[
\mathcal{Z}(\mathcal{E}_4) = \vert 000 + 210 + 012 + 040 \vert^2 + \vert 101 + 400 + 121 + 004 \vert^2 + 4 \vert 111 \vert^2.
\]
It introduces a partition on the set of exponents, def\/ined as  the i-irreps corresponding to the nine non-zero diagonal entries of $\mathcal{M}$:  $\{ 000, 210,  012, 040, 101, 400, 121, 004, 111\}$. To our knowledge, this invariant was f\/irst obtained in~\cite{SchellekensYankielowicz}.

\subsubsection[${\rm SU}(4) \subset {\rm SU}(10)$, $k=6$]{$\boldsymbol{{\rm SU}(4) \subset {\rm SU}(10)}$, $\boldsymbol{k=6}$}

\begin{itemize}\itemsep=0pt

\item
 At level 1, there are ten  i-irreps for $A_{9}$, namely $\{0,0,0,0,0,0,0,0,0\}$,  and $\{0,\ldots 0,1,0$, $\ldots,0\}$. From equation (\ref{confdim}) we calculate their conformal weights:
$\left\{0,\frac{9}{20},\frac{4}{5},\frac{21}{20},\frac{6}{5},\frac{5}{4},\frac{6}{5},\frac{21}{20}\right.$,
$\left.\frac{4}{5},\frac{9}{20}\right\}$.

\item    At  level $6$, we calculate the $84$ conformal weights for ${\rm SU}(4)$ i-irreps and f\/ind (use ordering def\/ined previously):
\[\mbox{\tiny
$\begin{array}{ccccccccccccccccc}
0&\frac{3}{16}&\frac{1}{4}&\frac{3}{16}&\frac{9}{20}&\frac{39}{80}&\frac{2}{5}&\frac{3}{5}&\frac{39}{80}&\frac{9}{20}&\frac{63}{80
   }&\frac{4}{5}&\frac{11}{16}&\frac{71}{80}&\frac{3}{4}&\frac{11}{16}&\frac{21}{20} \\
{} & \\
   \frac{71}{80}&\frac{4}{5}&\frac{63}{80}&\frac{6}{5}&
   \frac{19}{16}&\frac{21}{20}&\frac{5}{4}&\frac{87}{80}&1&\frac{111}{80}&\frac{6}{5}&\frac{87}{80}&\frac{21}{20}&\frac{8}{5}&\frac{111}{
   80}&\frac{5}{4}&\frac{19}{16} \\
{} & \\
   \frac{6}{5}&\frac{27}{16}&\frac{33}{20}&\frac{119}{80}&\frac{27}{16}&\frac{3}{2}&\frac{111}{80}&\frac{9}
   {5}&\frac{127}{80}&\frac{29}{20}&\frac{111}{80}&\frac{159}{80}&\frac{7}{4}&\frac{127}{80}&\frac{3}{2}&\frac{119}{80}&\frac{9}{4}\\
{} & \\
   \frac{159}{80}&\frac{9}{5}&\frac{27}{16}&\frac{33}{20}&\frac{27}{16}&\frac{9}{4}&\frac{35}{16}&2&\frac{11}{5}&\frac{159}{80}&\frac{37}{20}&
   \frac{183}{80}&\frac{41}{20}&\frac{151}{80}&\frac{9}{5}&\frac{49}{20}&\frac{35}{16}\\
{} & \\
   2&\frac{151}{80}&\frac{37}{20}&\frac{43}{16}&\frac
   {12}{5}&\frac{35}{16}&\frac{41}{20}&\frac{159}{80}&2&3&\frac{43}{16}&\frac{49}{20}&\frac{183}{80}&\frac{11}{5}&\frac{35}{16}&\frac{9}{4}&0
   \end{array}$}
\]

\item  The dif\/ference between conformal weights of $A_{9}$ and $A_3$ should be an integer. This constraint selects ten possibilities that give only necessary conditions for branching.
 Imposing the modularity constraint eliminates several entries (that we crossed-out in the next tab\-le). One f\/inds actually  two solutions but only one is a sum of squares (the other solution corresponds to the ``conjugated graph'' $\mathcal{E}_6^c$, see our discussion in Section~\ref{E6conj}):
\[\mbox{\tiny
$\begin{array}{lcl}
\begin{array}{rcl}
 000000000 &\stackrel{?}{\hookrightarrow}&  0 0 0+ 2 0 2+ \xout{5 0 1}+ 2 2 2+  \xout{1 0 5}+ 0 6 0 \\
100000000 &\stackrel{?}{\hookrightarrow}& 2 0 0+ \xout{0 0 2}+ 2 1 2+ \xout{2 4 0}+ 0 4 2 \\
010000000 &\stackrel{?}{\hookrightarrow}& \xout{2 1 0}+ 0 1 2+ 2 3 0+ \xout{0 3 2}+ 3 0 3 \\
001000000 &\stackrel{?}{\hookrightarrow}& 0 3 0+ \xout{3 0 1}+ 1 0 3+ 3 2 1+ \xout{1 2 3} \\
000100000 &\stackrel{?}{\hookrightarrow}& 4 0 0+ 1 2 1+ \xout{0 0 4}+ \xout{4 2 0}+ 0 2 4
\end{array}
&
\qquad
&
\begin{array}{rcl}
000010000 &\stackrel{?}{\hookrightarrow}& \xout{0 1 0}+ 2 2 0+ 0 2 2+\xout{ 0 5 0}+ 6 0 0+ 0 0 6 \\
000001000 &\stackrel{?}{\hookrightarrow}& \xout{4 0 0}+ 1 2 1+ 0 0 4+ 4 2 0+ \xout{0 2 4} \\
000000100 &\stackrel{?}{\hookrightarrow}& 0 3 0+ 3 0 1+ \xout{1 0 3}+ \xout{3 2 1}+ 1 2 3 \\
000000010 &\stackrel{?}{\hookrightarrow}& 2 1 0+\xout{ 0 1 2}+ \xout{2 3 0}+ 0 3 2+ 3 0 3 \\
000000001 &\stackrel{?}{\hookrightarrow}& \xout{2 0 0}+ 0 0 2+ 2 1 2+ 2 4 0+ \xout{0 4 2}
\end{array}
\end{array}$}
\]
\end{itemize}

The partition function obtained from the diagonal invariant of $A_{9}$  reads:
\begin{gather*}
\mathcal{Z}(\mathcal{E}_6)  =   \vert 000 + 060 + 202 + 222 \vert^2 + \vert 042 + 200 + 212 \vert^2 +  \vert 012 + 230 + 303 \vert^2   \\
\phantom{\mathcal{Z}(\mathcal{E}_6)  =}{}  + \vert 030 + 103 + 321\vert^2
 +     \vert 024 + 121 + 400\vert^2  +  \vert 006 + 022 + 220 + 600 \vert^2 \\
 \phantom{\mathcal{Z}(\mathcal{E}_6)  =}{}
 +  \vert 004 + 121 + 420  \vert^2  +   \vert 030 + 123 + 301 \vert^2
  +    \vert 032 + 210 + 303  \vert^2  \\
\phantom{\mathcal{Z}(\mathcal{E}_6)  =}{}   +  \vert 002 + 212 + 240 \vert^2.
\end{gather*}
It introduces a partition on the set of exponents, which are, by def\/inition,  the $32$  i-irreps corresponding to the non-zero diagonal entries of~$\mathcal{M}$.
To our knowledge, this invariant was f\/irst obtained in~\cite{AltschulerBauerItzykson}.

\subsubsection[${\rm SU}(4) \subset {\rm Spin}(20)$, $k=8$]{$\boldsymbol{{\rm SU}(4) \subset {\rm Spin}(20)}$, $\boldsymbol{k=8}$}

\begin{itemize}\itemsep=0pt

\item
At level 1, there are only four i-irreps for $D_{10}$, namely $\{0,0,0,0,0,0,0,0;0,0\}$,  $\{1,0,0,0,0$, $0,0,0;0,0\}$,  $\{0,0,0,0,0,0,0,0;1,0\}$, $\{0,0,0,0,0,0,0,0;0,1\}$; the last two entries refer to the fork of the $D$ graph. These i-irreps correspond to the trivial, the vectorial and the two half-spinorial representations.
From equation (\ref{confdim}) we calculate their conformal weights: $\left\{0,\frac{1}{2},\frac{5}{4},\frac{5}{4}\right\}$.

\item    At  level $8$, we calculate the $165$ conformal weights for ${\rm SU}(4)$ i-irreps and f\/ind (use ordering def\/ined previously):
\[\mbox{\tiny
$\begin{array}{ccccccccccccccccc}
 0 & \frac{5}{32} & \frac{5}{24} & \frac{5}{32} & \frac{3}{8} & \frac{13}{32} & \frac{1}{3} & \frac{1}{2} & \frac{13}{32} & \frac{3}{8} & \frac{21}{32} & \frac{2}{3} & \frac{55}{96} & \frac{71}{96} & \frac{5}{8} & \frac{55}{96} & \frac{7}{8} \\[0.1cm]
 \frac{71}{96} & \frac{2}{3} & \frac{21}{32} & 1 & \frac{95}{96} & \frac{7}{8} & \frac{25}{24} & \frac{29}{32} & \frac{5}{6} & \frac{37}{32} & 1 & \frac{29}{32} & \frac{7}{8} & \frac{4}{3} & \frac{37}{32} & \frac{25}{24} & \frac{95}{96} \\[0.1cm]
 1 & \frac{45}{32} & \frac{11}{8} & \frac{119}{96} & \frac{45}{32} & \frac{5}{4} & \frac{37}{32} & \frac{3}{2} & \frac{127}{96} & \frac{29}{24} & \frac{37}{32} & \frac{53}{32} & \frac{35}{24} & \frac{127}{96} & \frac{5}{4} & \frac{119}{96} & \frac{15}{8} \\[0.1cm]
 \frac{53}{32} & \frac{3}{2} & \frac{45}{32} & \frac{11}{8} & \frac{45}{32} & \frac{15}{8} & \frac{175}{96} & \frac{5}{3} & \frac{11}{6} & \frac{53}{32} & \frac{37}{24} & \frac{61}{32} & \frac{41}{24} & \frac{151}{96} & \frac{3}{2} & \frac{49}{24} & \frac{175}{96} \\[0.1cm]
 \frac{5}{3} & \frac{151}{96} & \frac{37}{24} & \frac{215}{96} & 2 & \frac{175}{96} & \frac{41}{24} & \frac{53}{32} & \frac{5}{3} & \frac{5}{2} & \frac{215}{96} & \frac{49}{24} & \frac{61}{32} & \frac{11}{6} & \frac{175}{96} & \frac{15}{8} & \frac{77}{32} \\[0.1cm]
 \frac{7}{3} & \frac{69}{32} & \frac{223}{96} & \frac{17}{8} & \frac{191}{96} & \frac{19}{8} & \frac{69}{32} & 2 & \frac{61}{32} & \frac{239}{96} & \frac{9}{4} & \frac{199}{96} & \frac{47}{24} & \frac{61}{32} & \frac{8}{3} & \frac{77}{32} & \frac{53}{24} \\[0.1cm]
 \frac{199}{96} & 2 & \frac{191}{96} & \frac{93}{32} & \frac{21}{8} & \frac{77}{32} & \frac{9}{4} & \frac{69}{32} & \frac{17}{8} & \frac{69}{32} & \frac{77}{24} & \frac{93}{32} & \frac{8}{3} & \frac{239}{96} & \frac{19}{8} & \frac{223}{96} & \frac{7}{3} \\[0.1cm]
 \frac{77}{32} & 3 & \frac{93}{32} & \frac{65}{24} & \frac{23}{8} & \frac{85}{32} & \frac{5}{2} & \frac{93}{32} & \frac{8}{3} & \frac{239}{96} & \frac{19}{8} & 3 & \frac{263}{96} & \frac{61}{24} & \frac{77}{32} & \frac{7}{3} & \frac{101}{32} \\[0.1cm]
 \frac{23}{8} & \frac{85}{32} & \frac{5}{2} & \frac{77}{32} & \frac{19}{8} & \frac{27}{8} & \frac{295}{96} & \frac{17}{6} & \frac{85}{32} & \frac{61}{24} & \frac{239}{96} & \frac{5}{2} & \frac{117}{32} & \frac{10}{3} & \frac{295}{96} & \frac{23}{8} & \frac{263}{96} \\[0.1cm]
 \frac{8}{3} & \frac{85}{32} & \frac{65}{24} & 4 & \frac{117}{32} & \frac{27}{8} & \frac{101}{32} & 3 & \frac{93}{32} & \frac{23}{8} & \frac{93}{32} & 3 &{}&{}&{}&{}&{}
\end{array}$}
\]

\item  The dif\/ference between conformal weights of $D_{10}$ and $A_3$ should be an integer. This selects four possibilities that give only necessary conditions for branching. Imposing the modularity constraint implies eliminating entries   400, 004, 440, 044  from the f\/irst line.
 \begin{gather*}
 \mbox{\tiny $0000000000  \ \stackrel{?}{\hookrightarrow} \  000+\xout{400}+121+\xout{ 004}+141+412+214+800+ \xout{440}+080+\xout{044}+008$},\\
 \mbox{\tiny $1000000000 \  \stackrel{?}{\hookrightarrow} \  020 + 230 + 032 + 303 + 060 + 602 + 323 + 206,$} \\
 \mbox{\tiny $0000000010 \ \stackrel{?}{\hookrightarrow} \ 311 + 113 + 331 + 133,$} \\
 \mbox{\tiny $0000000001 \ \stackrel{?}{\hookrightarrow} \ 311 + 113 + 331 + 133.$}
\end{gather*}
Notice that the contribution comes from $0$ (the trivial representation), from the f\/irst vertex of $D_{10}$ and from the two vertices of the fork (they have identical branching rules).
\end{itemize}

The partition function obtained from the diagonal invariant $ \vert 0000000000\vert^2 +   \vert 1000000000\vert^2 +    \vert 0000000010\vert^2 +  \vert 0000000001\vert^2$ of $D_{10}$ reads:
\begin{gather*}
\mathcal{Z}(\mathcal{E}_8)  =  \vert 000 + 121 + 141 + 412 + 214 + 800 + 080 + 008 \vert^2 + 2    \vert 311+113+331+133 \vert^2 \\
\phantom{\mathcal{Z}(\mathcal{E}_8)  =}{} + \vert 020 +230+032+060+303+602+323+206\vert^2.
\end{gather*}
It introduces a partition on the set of exponents, which are, by def\/inition,  the $20$  i-irreps corresponding to the non-zero diagonal entries of~$\mathcal{M}$.
To our knowledge, this invariant was f\/irst obtained in~\cite{AldazabalEtAl}.

\subsubsection{Quantum dimensions and cardinalities}

{\bf Quantum dimensions for $\boldsymbol{\mathcal{A}_k({\rm SU}(4))}$.}
  Multiplication by  its generators (associated with fundamental representations of ${\rm SU}(4)$) is encoded by a fusion matrix that may be considered as the adjacency matrix of a graph
 with three types of edges (self-conjugated fundamental representation corresponds to non-oriented edges).    Its vertices  build the Weyl  alcove of ${\rm SU}(4)$ at level $8$: a tetrahedron (in $3$-space) with $k$  f\/loors.
 It is convenient to think that $\mathcal{A}_k$ is a~quantum discrete group with $\vert \widehat {\mathcal{A}}_k \vert = r_A$ representations.   The quantum dimension $\dim(n)$ of a representation $n$ is calculated, for example, from the quantum Weyl formula. For the fundamental representations $f=\{1,0,0\}, \{0,1,0\}, \{0,0,1\}, $ (with classical dimensions  4, 6, 4)  one f\/inds: $\dim(f)=\{4_q ,\,  {3_q 4_q}/{2_q}, \, 4_q\} $.  In particular\footnote{\label{qnumdef}We set $n_q = (q^n - q^{-n}) /(q-q^{-1})$, with $q^{\kappa}=-1, \kappa = k+g$ and  $g=4$ for ${\rm SU}(4)$.}, $\beta = 4_q = 4 \cos \left(\frac{\pi }{\kappa }\right) \cos \left(\frac{2  \pi }{\kappa }\right)$ with $\kappa = k + 4$. The square of $\beta$ is the Jones index.  The quantum cardinality (also called quantum mass, quantum order, or ``global dimension''  like in \cite{EtingofOnVafa}) of this quantum discrete space, is obtained by summing the square of quantum dimensions for all $r_A$ simple objects: $\vert \mathcal{A}_k\vert = \sum_n \dim (n)^2$. Details are given in the Appendix.

 \begin{itemize}\footnotesize\itemsep=0pt
 \item If $k=4$, $r_A=35$,  $\dim (f)$:  $\left\{\beta=\sqrt{2 \left(2+\sqrt{2}\right)},2+\sqrt{2},\sqrt{2 \left(2+\sqrt{2}\right)}\right\}$,  $\vert \mathcal{A}_4\vert = 128 \left(3+2 \sqrt{2}\right)$.
  \item If $k=6$, $r_A=84$,   $\dim(f)$:  $\left \{ \beta=\sqrt{5+2 \sqrt{5}},2+\sqrt{5},\sqrt{5+2 \sqrt{5}}\right\}$, $\vert \mathcal{A}_6\vert = 800 \left(9+4 \sqrt{5}\right)$.
  \item  If $k=8$, $r_A=165$,   $\dim(f)$:  $ \left\{\beta=\sqrt{3 \left(2+\sqrt{3}\right)},3+\sqrt{3},\sqrt{3 \left(2+\sqrt{3}\right)}\right\}$,  $\vert \mathcal{A}_8\vert = 3456 \left(26+15 \sqrt{3}\right)$.
\end{itemize}

 {\bf Quantum dimensions for $\boldsymbol{\mathcal{E}_k({\rm SU}(4))}$, $\boldsymbol{\{k=4,6,8\}}$.\ First method.}
 Action of  $\mathcal{A}_k$ on~$\mathcal{E}_k$   is encoded by matrices generically called ``annular matrices'' (they are also called ``nimreps'' in the literature, but this last term is sometimes used to denote other types of matrices with non negative integer entries). In particular, action of the generators is described by annular matrices that we consider as adjacency matrices for the graph  $\mathcal{E}_k$ itself.  Once the later is obtained, one calculates quantum dimensions $\dim (a)$ for its $r_E$ vertices (the simple objects) by using for instance the Perron--Frobenius vector of the annular matrix associated with the genera\-tor~$F_{\{1,0,0\}}$.
Its quantum cardinality  is then def\/ined by $\vert \mathcal{E}_k \vert = \sum_{a}  \dim (a)^2$.  The problem is that we do not know, at this stage, the values $\dim (a)$ for all vertices $a$ of $ \mathcal{E}_k $, since this graph will only be determined later.

{\bf Quantum dimensions for $\boldsymbol{\mathcal{E}_k({\rm SU}(4))}$, $\boldsymbol{\{k=4,6,8\}}$.\ Second method.}
 It is convenient to think that $\mathcal{A}_k/\mathcal{E}_k$ is a homogenous space, both discrete and quantum.
Like in a classical situation, we have\footnote{These maps (actually functors) are described by the square annular matrices $F_n$ or by the rectangular essential matrices $E_a$ with $(E_a)_{nb}=(F_n)_{ab}$ that we shall introduce later.} restriction maps ${\mathcal A}_k \mapsto {\mathcal E}_k$ and induction maps $ {\mathcal E}_k  \mapsto {\mathcal A}_k$.  One may think that vertices of the quantum graph $\mathcal{E}_k$ do not only label irreducible objects $a$ of $\mathcal{E}$ but also space of sections of quantum vector bundles $\Gamma_a$ which can be decomposed, using induction, into irreducible objects of $\mathcal{A}_k$: we write $\Gamma_a = \bigoplus_{n  \uparrow \Gamma_a} n$. This implies, for quantum dimensions, the equality $\dim (\Gamma_a) = \bigoplus_{n  \uparrow \Gamma_a} \dim(n)$. The space of sections $\mathcal{F}=\Gamma_0$, associated with the identity representation, is special since it can be considered as the quantum algebra of functions over $\mathcal{A}_k/\mathcal{E}_k$.  Its dimension $\dim ( \Gamma_0) =  \vert \mathcal{A}_k/\mathcal{E}_k \vert$ is obtained by summing $q$-dimensions (not their squares!) of the ${n  \uparrow \Gamma_0}$ representations.  We are in a type-I situation (the modular invariant is a sum of blocks) and in this case, we make use of the following particular feature -- not true in general: the irreducible representations ${n  \uparrow \Gamma_0}$ that appear in the decomposition of $\Gamma_0$ are exactly those appearing in the f\/irst modular block of the partition function.  From the property  $\vert \mathcal{A}_k/\mathcal{E}_k \vert = \vert \mathcal{A}_k\vert /\vert \mathcal{E}_k \vert $, we f\/inally obtain  $\vert \mathcal{E}_k \vert$ by calculating  $\frac{\vert \mathcal{A}_k\vert}{\vert \mathcal{A}_k/\mathcal{E}_k \vert} $.

    \begin{itemize}\footnotesize\itemsep=0pt
 \item When $k=4$, we have  $\mathcal{F}=\Gamma_0  =  000 \oplus 210 \oplus 012 \oplus 040$ so that  $\dim  (\mathcal{F})   =  \dim  (\Gamma_0) = \vert \mathcal{A}/\mathcal{E} \vert = 8+4 \sqrt{2}$. Using the known value for $\vert \mathcal{A}\vert$ one obtains $\vert \mathcal{E} \vert  = 16 \left(2+\sqrt{2}\right)$.

  \item When $k=6$, we have  $\mathcal{F}=\Gamma_0  =  000 \oplus 060\oplus202 \oplus222$ so that  $\dim  (\mathcal{F})   =  \dim  (\Gamma_0) = \vert \mathcal{A}/\mathcal{E} \vert = 20+8 \sqrt{5}$. Using the known value for $\vert \mathcal{A}\vert$ one obtains $\vert \mathcal{E} \vert  = 40 \left(5+2 \sqrt{5}\right) $.

  \item  When $k=8$,  we have  $\mathcal{F}=\Gamma_0  =  000 \oplus  121 \oplus  141 \oplus  412 \oplus  214 \oplus  800 \oplus  080 \oplus  008$ so that  $\dim  (\mathcal{F})   =  \dim  (\Gamma_0) = \vert \mathcal{A}/\mathcal{E} \vert = 12(9 + 5 \sqrt 3)$. Using the known value for $\vert \mathcal{A}\vert$ one obtains $\vert \mathcal{E} \vert  =48 (9+5 \sqrt{3})$.
\end{itemize}

 {\bf Quantum dimensions for $\boldsymbol{\mathcal{E}_k({\rm SU}(4))}$, $\boldsymbol{\{k=4,6,8\}}$.\ Third method.}   \label{massE}
 The third method (which is probably the shortest, in the case of quantum graphs obtained from conformal embeddings) does not even use the expression of the f\/irst modular block but it uses some general results and concepts from the structure of the graph of quantum symmetries $\mathcal{O}(\mathcal{E}_k)$ that will be discussed in a coming section. In a nutshell, one uses the following known results:  1)~$\vert \mathcal{O}( \mathcal{E}_k)\vert =  \vert \mathcal{E} \vert \times \vert \mathcal{E} \vert / \vert {\mathcal J}_{\mathcal{O}} \vert$ where $ {\mathcal J}_O$ denotes the set of ambichiral vertices of the Ocneanu graph, 2)~$\vert \mathcal{A}_k\vert$ = $\vert \mathcal{O}(\mathcal{E}_k)\vert$,  and 3)~$ \vert {\mathcal J}_{{\mathcal O}} \vert = \vert {\mathcal J}_{{\mathcal E}} \vert$ where ${\mathcal J}_{{\mathcal E}}$ denote the sets of modular vertices of the graph~${\mathcal E}_k$. Finally, one notices that for a case coming from a conformal embedding $K={\rm SU}(4) \subset G$, one can identify vertices $c \in  \mathcal{J}_{{\mathcal E}} \subset  \mathcal{E}$ with vertices $c \in {\mathcal J}= {\mathcal A}_1(G)$.
The conclusion is that one can f\/irst calculate $\vert {\mathcal J} \vert = \sum_{s} \dim (s)$ as the mass of the small quantum group ${\mathcal A}_1(G)$, and f\/inally obtain  $ \vert \mathcal{E}_k \vert $ from the following relation\footnote{In this paper $G={\rm SU}(4)$ and $K={\rm Spin}(15), {\rm SU}(10), {\rm Spin}(20)$ for $k=4,6,8$,  but this relation is valid for any case stemming from a conformal embedding.}:
\[
 \vert \mathcal{E}_k (K) \vert = \sqrt{\vert  \mathcal{A}_k(K) \vert \times \vert  \mathcal{A}_1(G)  \vert }.
\]
Values for  quantum cardinality of the (very) small quantum groups $\vert  \mathcal{A}_1(G)  \vert $ at relevant\footnote{For a conformal embedding $K \subset G$,  the value of $q$ used to study $\mathcal{A}_k(K)$ is not the same as the value of $q$ used to study $\mathcal{A}_1(G)$, since $q$ is given by $\exp(i \pi /(k+g))$: for instance one uses $q^{12}=-1$ for  $\mathcal{A}_8({\rm SU}(4))$ but $q^{19}=-1$ for $\mathcal{A}_1({\rm Spin}(20))$.}
va\-lues of $q$ are obtained in the appendix. One f\/inds\footnote{One should not think that $\vert {\mathcal J} \vert$ is always an integer: compute for instance  $\vert {\mathcal A}_1(G_2) \vert$ which can be used to determine $E_8= \mathcal{E}_{28} ({\rm SU}(2))$. However, $\vert {\mathcal A}_1({\rm SU}(g)) \vert = g$.}
  $\vert {\mathcal A}_1({\rm Spin}(15)) \vert = 4$, $\vert {\mathcal A}_1({\rm SU}(10)) \vert = 10$, $\vert {\mathcal A}_1({\rm Spin}(20)) \vert = 4$ and we recover the already given results for $ \vert \mathcal{E}_k \vert$. Incidentally this provides another check that obtained branching rules are indeed correct.

\begin{remark}
 We stress the fact that the calculation of $\vert \mathcal{E}_k \vert$ can be done, using the second or the third method,  before having determined the quantum graph  $\mathcal{E}_k$ itself, in particular without using any knowledge of the quantum dimensions $\dim (a)$ of its vertices.  Once the graph is known, one can obtain these quantum dimensions from a Perron--Frobenius eigenvector, then check the consistency of calculations by using induction, from the relation $\dim(a) = \dim (\Gamma_a) /\dim (\Gamma_0)$,  and f\/inally recover the quantum cardinality of $\mathcal{E}$ by a direct calculation (f\/irst method).
\end{remark}

\section{Algebras of quantum symmetries}

\subsection{General terminology and notations}\label{commutativityofO}

We introduce some terminology and several notations used in the later sections.

Fusion ring ${\mathcal A}_k$: the commutative ring spanned by integrable irreducible representations $m,n,\ldots$ of the af\/f\/ine Lie algebra of ${\rm SU}(4)$ at level $k$, of dimension $r_A=(k+1)(k+2)(k+3)/3! $. Structure constants are encoded by fusion matrices $N_m$ of dimension $r_A \times r_A$:  $m \cdot n = \sum_p (N_m)_{np} \, p$. Indices refer to Young Tableaux or to weights. Existence of duals implies, for the fusion ring,  the rigidity property  $(N_{\overline m})_{np} = (N_m)_{pn}  $, where $\overline m$ refers to the conjugate of the irreducible representation $m$. In the case of ${\rm SU}(4)$, we have three generators (fundamental irreducible representations): one of them is real (self-conjugated) and the other two are conjugated from one another.

${\mathcal A}_k$ acts on the additive group spanned by vertices $a,b,\ldots$ of the quantum graph ${\mathcal E}_k$. This module action is encoded by annular matrices $F_m$: $m \cdot a= \sum_b (F_m)_{ab}\, b $. These are square matrices of dimension $r_E \times r_E$, where $r_E$ is the number of simple objects (i.e., vertices of the quantum graph) in  ${\mathcal E}_k$.  To the fundamental representations of ${\rm SU}(4)$ correspond particular annular matrices which are the adjacency matrices of the quantum graph. The rigidity property of  ${\mathcal A}_k$ implies\footnote{For the ${\rm SU}(2)$ theory, this property excludes non-ADE Dynkin diagrams.} $(F_{\overline n})_{ab} = (F_n)_{ba}  $.
 It is convenient to introduce a family of rectangular matrices called ``essential matrices'' \cite{Coque:Qtetra}, via the relation $(E_a)_{mb} = (F_m)_{ab}$.  When $a$ is the origin\footnote{A particular vertex of ${\mathcal E}$ is always distinguished.} $0$ of the quantum  graph, $E_0$ is usually called ``the intertwiner''.

In general there is no multiplication in ${\mathcal E}_k$, with non-negative integer structure constants,  compatible with the action of the fusion ring. When it exists, the quantum graph is said to possess self-fusion.  This is the case in the three examples under study. The multiplication is described by another family of matrices $G_a$ with non negative integer entries: we write $a \cdot  b =  \sum_c (G_a)_{bc}\, c $;   compatibility with the fusion algebra (ring)  reads $ m \cdot (a \cdot b) = (m \cdot a) \cdot b $, so that $(G_a \cdot F_m) = \sum_c (F_m)_{ac} \, G_c $.

The additive group ${\mathcal E}_k$ is not only a  $\ZZ_+$   module over the fusion ring ${\mathcal A}_k$, but also a $\ZZ_+$   module over the Ocneanu ring (or algebra) of quantum symmetries ${\mathcal O}$.  Linear generators of this ring are denoted $x,y,\ldots$ and its structure constants, def\/ined by $x\cdot y = \sum_z (O_x)_{yz} \, z $ are encoded by  the ``matrices of quantum symmetries''~$O_x$.  To each fundamental irreducible representation~$f$ of~${\rm SU}(4)$ one associates two fundamental generators of~${\mathcal O}$, called chiral left $f^L$ and chiral right~$f^R$.  So, ${\mathcal O}$ has $6 = 2 \times 3$ chiral generators.  Like in usual representation theory, all other linear generators of the algebra appear when we decompose  products of fundamental (chiral) generators. The Cayley graph of multiplication by the  chiral generators (several types of lines),  called the Ocneanu graph of ${\mathcal E}_k$, encodes the algebra structure of ${\mathcal O}$. Quantum symmetry matrices $O_x$ have dimension $r_O \times r_O$, where $r_O$ is the number of vertices of the Ocneanu graph.
Linear generators  that appear in the decomposition of products of left (right)  chiral  generators span a subalgebra called the left (right) chiral subalgebra. These two subalgebras are not necessarily commutative but the left and the right commute. Intersection of left and right chiral subalgebras is called the ambichiral subalgebra.
The module action of ${\mathcal O}$ on ${\mathcal E}_k$ is encoded by  ``dual annular matrices''~$S_x$, def\/ined by $x \cdot a = \sum_b (S_x)_{ab} \, b  $.

From general results obtained in operator algebra by \cite{Ocneanu:Unpublished} and  \cite{Evans-I, Evans-II, Evans-Kawahigashi}, translated to a categorical language by \cite{Ostrik}, one  shows that the ring of quantum symmetries ${\mathcal O}$ is a bimodule over the fusion ring  $ {\mathcal A}_k$. This action reads, in terms of generators,  $m \cdot x \cdot n = \sum_y   (V_{mn})_{xy} \, y $, where $m$, $n$ refer to irreducible objects of  ${\mathcal A}_k$ and $x,y$ to irreducible objects of ${\mathcal O}$.
Structure constants are encoded by the ``double-fusion matrices'' $V_{mn}$, with matrix elements $(V_{mn})_{xy}$, again non negative integers. To the fundamental representations $f$ of ${\rm SU}(4)$ correspond particular double fusion matrices encoding the multiplication by chiral generators in $\mathcal{O}$ (adjacency matrices of the Ocneanu graph): $V_{f0}=O_{f^L}$ and $V_{0f}=O_{f^R}$, where 0 is to the trivial representation of ${\rm SU}(4)$.

One also introduces the family of so-called toric matrices $W_{xy}$, with matrix elements $(W_{xy})_{mn}$ $=(V_{mn})_{xy}$.
When both $x$ and $y$ refer to the unit object  of $\mathcal{O}$ (that we label $0$), one recovers the modular invariant ${\mathcal M}=W_{00}$ encoded by the partition function $\mathcal{Z}$ of the corresponding conformal f\/ield theory.  As explained in~\cite{PetkovaZuber:Oc}, when one or two indices $x$ and $y$ are non trivial, toric matrices are interpreted as partition functions on a torus,  in a conformal theory of type $ {\mathcal A}_k$,  with boundary type conditions specif\/ied by~${\mathcal E}$,  but with defects specif\/ied by $x$ and $y$. Only ${\mathcal M}$ is modular invariant (it commutes with the generators $s$ and $t$ of~${\rm SL}(2,\ZZ)$ given by Kac--Peterson formulae). Toric matrices were f\/irst introduced and calculated by Ocneanu (unpublished) for theories of type~${\rm SU}(2)$. Various methods to compute or def\/ine them can be found in
\cite{FuchsRunkelSchweigert-I,Coque:Qtetra,PetkovaZuber:Oc}. Reference~\cite{GilCoque:ADE} gives explicit expressions for all $W_{x0}$, for all members of the ${\rm SU}(2)$ family ($ADE$ graphs).
Left and right associativity constraints $(m\cdot(n\cdot x \cdot p) \cdot q)=(m\cdot n)\cdot x \cdot (p\cdot q)  $ for the $ {\mathcal A} \times {\mathcal A}$ bimodule structure of ${\mathcal O}$ can be written in terms of fusion and toric matrices. A particular case of this equality reads\footnote{Equivalently, one can write $(N_\sigma \otimes N_\tau) {\mathcal M}_{\sigma \tau} = P_{1324} \sum_x W_x \otimes W_x$ where $W_x = W_{x0}$ and $P_{1324}$ is a permutation.}
$
\sum_x (W_{0x})_{m n} \, W_{x0} = N_{m}\, {\mathcal M}  \,N_{n}^{\rm tr}
$.
It was presented by A.~Ocneanu in \cite{Ocneanu:Bariloche} and called the  ``modular splitting equation''. Another particular case of the bimodule associativity constraints gives the following ``intertwining equations'':
$
\sum_y   (W_{xy})_{mn} \, W_{y0} = N_m \, W_{x0} \,  N_n^{\rm tr}
$.
A~practical method to solve this system (matrix elements should be non-negative integers) is discussed in~\cite{EstebanGil}, with several ${\rm SU}(3)$ examples. Given fusion matrices~$N_m$, known in general,  and a~modular invariant matrix~${\mathcal M}=W_{00}$,
solving the modular splitting equation, i.e., f\/inding the~$W_{x0}$, and subsequently solving the intertwining equations, allows one to construct the chiral generators of~${\mathcal O}$ and obtain the graph of quantum symmetries (and the graph ${\mathcal E}_k$ itself, as a by-product).  This is what we do in the next section, starting from the three exceptional partition functions of type ${\rm SU}(4)$ obtained previously.

\subsection{Method of resolution (summary)}\label{solvingintertwining}

The following program should be carried out for all examples:
\begin{enumerate}\itemsep=0pt
\item Solve the modular splitting equation (f\/ind toric matrices $W_{x0}$).
\item Solve the intertwining equations (i.e., f\/ind generators $O_x$ for the Ocneanu algebra and its graph ${\mathcal O}({\mathcal E}_k)$), obtain the permutations describing chiral transposition.
\item Determine the quantum graph ${\mathcal E}_k$ (f\/ind its adjacency matrices).
\item Determine the annular matrices $F_n$ describing ${\mathcal E}_k$ as a module over the fusion algebra ${\mathcal A}_k$.
\item Describe the self-fusion on ${\mathcal E}_k$ (f\/ind matrices $G_a$).
\item Reconstruct  ${\mathcal O}({\mathcal E}_k)$ in terms of ${\mathcal E}_k$.
\item Determine dual annular matrices $S_x$ describing ${\mathcal E}_k$ as a module over ${\mathcal O}({\mathcal E}_k)$.
\item Checks: reconstruct toric matrices from the previous realization of  ${\mathcal O}({\mathcal E}_k)$, verify the  relation\footnote{\label{alphainduction}
Since ${\mathcal O}({\mathcal E}_k)$  is an ${\mathcal A}_k$ bimodule, we obtain in particular two algebra homomorphisms from the later to the f\/irst (this should coincide with the notion of ``alpha induction'' introduced in \cite{Longo-Rehren, Xu:1998} and used in \cite{Evans-II}): $\alpha_L(m) = m \, 0_{Oc} \, 0_{\mathcal A}$ and
$\alpha_R(m) =  0_{\mathcal A} \, 0_{Oc} \, m$,  for $m \in {\mathcal A}_k$. They can be explicitly written in terms of toric matrices:  $\alpha_L(m) = \sum_y (W_{0y})_{m0} \, y$ and   $\alpha_R(m) = \sum_y (W_{0y})_{0m} \, y$.
If we compose these two maps with the homomorphism $S$ from
${\mathcal O}({\mathcal E}_k)$ to ${\mathcal E}_k$, described by dual annular matrices, we obtain a morphism from ${\mathcal A}_k$ to ${\mathcal E}_k$ that has to coincide with the one def\/ined by annular matrices, so that  we obtain  the identity $F_m =  \sum_y (W_{0y})_{m0} \,  S_y  = \sum_y  (W_{0y})_{0m} \,  S_y$ that implies, in particular $d_m =  \sum_y(W_{0y})_{m0} \,  d_y$.},  expressing $F_n$ in terms of $S_x$,  check  identities for quantum cardinalities, etc.
\item Describe the two multiplicative structures of the associated quantum groupo\"\i d ${\mathcal B}$.
\item Matrix units and block diagonalization of ${\mathcal E}_k$ and ${\mathcal O}({\mathcal E}_k)$.
\item Check consistency equations (self-connection on triangular cells).
\end{enumerate}

{\bf Determination of the toric matrices  $\boldsymbol{W_{z0}}$.}
\label{sec:eqsofmodsplit}
 These matrices, of size $r_A \times r_A$, are obtained by solving the modular splitting equation. This is done by using the following algorithm.
For each choice of the pair $({m , n})$ (i.e., $r_A^2$ possibilities), we  f\/irst def\/ine and calculate the matrices
${\mathcal K} _{mn} = N_{m}\,\mathcal{M}\,N_{n}^{\rm tr}$.
The modular splitting equation reads:
\begin{equation}
{\mathcal K} _{mn} = \sum_{x=0}^{r_O-1} (W_{0x})_{mn} \, W_{x0}  .
\label{Kmse}
\end{equation}
It can be viewed as the linear expansion of the matrix ${\mathcal K}
_{mn}$ over the set of toric matri\-ces~$W_{x0}$, where the coef\/f\/icients
of this expansion are the non-negative integers $(W_{0x})_{mn}$, and
$r_O={\rm Tr}(\mathcal{M}\mathcal{M}^{\dagger})$ is the dimension of the
quantum symmetry algebra.
This set of equations has to be solved for all possible values of $m$ and $n$. In other
words,  we have a single equation for a huge tensor ${\mathcal K}$ with
$r_A^2 \times r_A^2$ components viewed as a family of  $r_A^2$ vectors
${\mathcal K} _{mn}$, each vector being itself a $r_A\times r_A$ matrix.
 In general the family of toric matrices is not free: the $r_O$ toric
matrices~$W_{z0}$ are not linearly independent and span (like matrices
${\mathcal K}_{mn}$)  a vector space of dimension $r_W < r_O$; this
feature (related to the possible non-commutativity of ${\mathcal
O}({\mathcal E}_k)$) appears whenever the modular invariant
$\mathcal{M}$ has coef\/f\/icients bigger than $1$. Toric matrices $W_{z0}$
are obtained by using the following iterative algorithm already used in
\cite{EstebanGil,RobertGil:E4SU4}. The algebra of quantum symmetries comes with a
basis, made of the linear generators that we called $x$, which is
special because structure constants of the algebra, in this basis, are
non negative integers. We def\/ine a scalar product in the underlying
vector space for which the $x$ basis is orthonormal, and consider, for
each matrix ${\mathcal K} _{mn}$,  and because of equation~(\ref{Kmse}), the
vector  $\sum_x (W_{0x})_{mn}$, $x \in {\mathcal O}$,  whose norm,
abusively\footnote{Indeed, it may happen that two toric matrices $W_{x0}$,
$W_{y0}$  appearing on the r.h.s.\ of~(\ref{Kmse}) are equal, even though
$x\neq y$.} called norm of  ${\mathcal K} _{m n}$ and  denoted  $\vert
\vert {\mathcal K} _{m n} \vert \vert^2$,  is equal to $\sum_x
|(W_{0x})_{mn}|^2$. The relation
$V_{\overline{m}\overline{n}}=(V_{mn})^{\rm tr}$ (see later) and the modular
splitting equation imply that $\vert \vert {\mathcal K} _{m n} \vert
\vert^2 = ({\mathcal K} _{m n})_{\overline{m} \overline{n}}  =  \sum_{p,
q} (N_{m})_{\overline{m} p}\, ({\mathcal M})_{pq}  \,(N_{n}^{\rm tr})_{q
\overline{n}}$, i.e., for each $m$, $n$, this norm can be directly read from
the matrix ${\mathcal K} _{mn}$ itself.
The next task is therefore to calculate the norm of the matrices ${\mathcal K} _{m n}$.
The toric matrices $W_{x0}$
themselves are then obtained by considering matrices ${\mathcal K}_{mn}$
of increasing norms $1,2,3,\ldots$. For example those of norm $1$
immediately def\/ine toric matrices (since the sum appearing on the r.h.s.\ of the modular splitting equation involves only one term), in particular one recovers $W_{00} =
\mathcal M$.  Then we analyse those of norm $2$, and so on. The process
ultimately stops since the rank $r_W$ is f\/inite. A case dependent
complication, leading to ambiguities in the decomposition of  ${\mathcal
K}_{mn}$,  stems from the fact that the family of toric matrices is not
free (see our discussion of the specif\/ic cases).

In order to ease the discussion of the resolution of the equation of modular splitting, it is convenient to introduce the following notations and def\/initions.
We order the i-irrep as in  Section~\ref{orderirrep} and call $m^\#$ the position of $m$, so that $\{0,0,0\}^\#=1$, $\{1,0,0\}^\#=2$, $\ldots$.
For each possible square norm $u$,  we set ${\mathcal K} ^u = \{{\mathcal K} _{m n} / \vert\vert {\mathcal K} _{m n} \vert \vert^2= u \}$; notice that this is def\/ined as a set: it may be that, for a given $M \in {\mathcal K} ^u$,  there exist  distinct pairs $(m , n)$, $(m' , n')$  such that $M= {\mathcal K} _{m n} = {\mathcal K} _{m' n'}$, but this matrix appears only once in ${\mathcal K} ^u$.
The tensor ${\mathcal K} $ is a square array of square matrices of dimension $r_A^2 \times r_A^2$. Lines and columns are ordered by using the previously given ordering on the set of irreducible representations.  We scan ${\mathcal K} $ from left to right and from top to bottom. This allow us to order the sets ${\mathcal K} ^u$:
for a given  $M \in {\mathcal K} ^u$ we take note of its f\/irst
occurrence, i.e., the number ${\rm Inf}\{ (m^\# - 1)r_A + n^\# \}$ over all
$\{(m,n)\}$ such that  $M= {\mathcal K} _{m n}$, this def\/ines a strict
order on the set ${\mathcal K}^u$ (use the fact that $m^\# < r_A$ and
$n^\#<r_A$). We can therefore refer to elements $M$ of ${\mathcal K} ^u$
by their position $v$, and we shall write  $M = {\mathcal K} ^u[v]$.

{\bf Conjugations.}
Complex conjugation is def\/ined on the set of irreducible representations of~${\rm SU}(4)$ which, in terms of fusion matrices, reads
$N_{\overline m} = N_{m}^{\rm tr}$.  At the level of the tensor square, this star representation def\/ined by $(\overline{V_{mn}}) = V_{{\overline m}\, {\overline n}}$ reads $V_{{\overline m}\, {\overline n}} = (V_{mn})^{\rm tr}$ since all matrices have non negative integral entries (no need to take conjugate of complex numbers). In terms of toric matrices, this implies $(W_{xy})_{{\overline m}\, {\overline n}}=(W_{yx})_{mn}$.
We have also a conjugation (``bar operation'')  $x \mapsto \overline{x}$ on the algebra of quantum symmetries,
 that maps toric matrices to toric matrices, $W_x = W_{x0} \mapsto W_{\overline x} = W_{0x}$.  More generally we set
 $W_{ \overline{x}  \, \overline{y}} = W_{yx}$.  Real  generators are def\/ined by the property $x=\overline{x}$, in that case we have  $W_{x0}=W_{0x}$, i.e., $W_x = W_{\overline{x}}$.

Toric matrices are usually not symmetric:  transposition is not trivial, but it leaves invariant the set of toric matrices and induces an operation called ``chiral transposition''\footnote{Chiral transposition may be related to the ``conjugation of defect lines'', as in \cite{FuchsRunkelSchweigert-II}. This operation is not {\it a~priori\/} def\/ined for all CFT's, but its existence, along with the invariance property of the set of toric matrices under transposition is, for the cases studied here,  an observational fact, consequence of our explicit determination of toric matrices, and quantum symmetry generators.},
denoted  $x \mapsto x^c$ on the algebra of quantum symmetries. It reads $W_{x^c y^c}=(W_{x y})^{\rm tr}$.
In particular $(W_{x^c 0})_{{m}\, {n}}=(W_{x0})_{nm}$.
Using $\alpha_{L,R}$ morphisms introduced in footnote \ref{alphainduction}, it reads ${\alpha_L(m)}^c = \alpha_R(m)$.
Symmetric  generators are def\/ined by the property $x=x^c$, in that case the corresponding toric matrices are symmetric: $W_{x}= (W_{x})^{\rm tr}$.
It is this operation that maps chiral left to chiral right generators: $(f^L)^c = f^R$.

The operation  $x \mapsto x^{\dag}$ obtained by composing the above two operations is called ``chiral adjoint''.
It is such that ${(m \, x \, n)}^{\dag}  = {\overline n}  \, x^{\dag} \, {\overline m}$.
In terms of toric matrices, it reads $(W_{x^\dag y^\dag})_{{m}\, {n}}=(W_{yx})_{nm}$.  Self-adjoint (or hermitian)  generators are def\/ined by the property  $x=x^\dag$, in that case $W_{x0}= (W_{0x})^{\rm tr}$, i.e., $W_x = (W_{\overline{x}})^{\rm tr}$. Ambichiral generators (remember that they span a subalgeb\-ra~${\mathcal J}$ def\/ined as intersection of the left and right chiral subalgebras) can be recognized as those  self-adjoint  generators whose corresponding vertices belong to the f\/irst connected component\footnote{It is def\/ined as the connected component of the graph of quantum symmetries, using the generator $\{100\}$, that contains the identity of the algebra $\mathcal{O}$, whose corresponding toric matrix is the modular invariant $\mathcal{M}=W_{00}.$} in the graph of quantum symmetries.

We summarize the above discussion by the following collection of equalities:
\[
W_{x^\dag y^\dag} = W_{y^c x^c} = (W_{ \overline{x}  \, \overline{y}})^{\rm tr} = (W_{y x})^{\rm tr}.
\]
For each of the above three conjugations, one can introduce a permutation matrix acting on the set of generators of ${\mathcal O}$, intertwining between $x$ and $\overline{x}$, ${x}^c$ or $x^\dag$.

Determination of the conjugations is not straightforward when $r_W < r_O$.
What we do is to parametrize the solutions found after analysis of the set of toric matrices and we use them to solve the set of intertwining equations (see next paragraph). Imposing that the obtained generators of ${\mathcal O}$ obey the expected constraints (see later) restrict the possible choices for the conjugations, and ultimately f\/ixes all free parameters, up to possible graph automorphisms.

In many cases, and in particular in the three exceptional cases that we consider in this article, one can realize the algebra of quantum symmetries ${\mathcal O}$ as a quotient, over the ambichiral subalgebra, of an algebra def\/ined in terms of the tensor square of the algebra of the quantum graph ${\mathcal E}$ (in simple cases,  ${\mathcal O}$ can be identif\/ied with  ${\mathcal E} \otimes {\mathcal E}/{\mathcal J}$). Using this realization, i.e., writing $x = a \otimes b$, the above three operations read:  ${\overline x} = {\overline a} \otimes {\overline b}$,  $x^c= b \otimes a$ and $x^\dag = {\overline b} \otimes {\overline a}$. Actually the conjugation $a \mapsto {\overline a}$ in ${\mathcal E}$  (we could very well choose ${\mathcal E}={\mathcal A}$) can be deduced from the same operation in ${\mathcal O}$ via the identif\/ication $a \simeq a \otimes \one$.


{\bf Solving the intertwining equations.}
The family of toric matrices ``with one twist'', i.e., the $W_{x0}$ matrices, was determined in a previous step, but we should determine all the matri\-ces~$W_{xy}$. For each triplet $(m,n,x)$, we def\/ine the matrices $\mathcal{K}_{mn}^{x} = N_m \, W_{x0} \,  {N_n}^{\rm tr}$ and calculate them. The intertwining equations (one matrix equation for each triplet) read:
\[
\mathcal{K}_{mn}^{x}   = \sum_y \, (W_{xy})_{mn} \, W_{y0}.
\]
It can be viewed as the linear expansion of the matrix ${\mathcal K}_{mn}^x$ over the set of toric matrices $W_{x0}$, where the coef\/f\/icients of this expansion are the non-negative integers $(W_{xy})_{mn} = (V_{mn})_{xy}$, that we want to determine.
In order to f\/ind the algebra of quantum symmetries and its graph, it is enough to solve only those equations involving the six chiral generators\footnote{Remember that the notation
$O_{f^L}$  does not refer to the chiral adjoint of  $O_{f^R}$ but to  its chiral transpose.}, i.e., to determine the matrices $ V_{f0} = O_{f^L} $ and $V_{0f} = O_{f^R}$, where $f$ refer to the three fundamental representations of~${\rm SU}(4)$, $f=\{100\}$, $f=\{010\}$, $f=\{001\}$.
In other words,  we solve the  intertwining equations
\[
N_f \, W_{x0} \,  {N_0}^{\rm tr} = \sum_y  ( O_{f^L} )_{xy}\, W_{y0} \qquad  \hbox{and} \qquad    N_0 \, W_{x0} \,  {N_f}^{\rm tr} = \sum_y  ( O_{f^R} )_{xy}\, W_{y0}.
\]
When the toric matrices $W_{y0}$ are linearly independent, the linear expansions of $\mathcal{K}_{f0}^{x}$ and $\mathcal{K}_{0f}^{x}$ are unique and the determination of the chiral generators $O_{f^L}$ and $O_{f^R}$ is straightforward (for any chosen $f$ one sets the elements  of matrices $O_f$ to unknown parameters and solves a system of linear equations in ${r_O}^2$ unknowns).
But in general the family of toric matrices is not free,  and even after imposing that matrix coef\/f\/icients of $O_{f^L}$ and $O_{f^R}$ should be non-negative integers, we are still left with a solution with many free parameters.
Some of them are determined by imposing the bar-conjugation relations: $O_{{\overline f}^{L}} = (O_{f^L})^{\rm tr}$ (respectively $O_{{\overline f}^{R}} = (O_{f^R})^{\rm tr}$). In particular,  for the real generator  $f =  \{010\} = {\overline f}$, matrices of the corresponding  left and right chiral generators should be symmetric.

Other parameters are determined by imposing the chiral transposition relations $(f^L)^c = f^R$, so that  $O_{f^{R}} = P\, O_{f^L} \, P^{-1}$  where $P$ is the permutation matrix implementing chiral transposition (so $P^2=1$); it can be obtained from our knowledge of toric matrices since an equality $y=x^c$ among the generators of ${\mathcal O}$ implies $W_y = (W_x)^{\rm tr}$. The operation $c$ (or the matrix $P$) is usually not fully determined at that stage since it may happen that two distinct generators $x$ and $y$ are represented by identical toric matrices $W_x=W_y$.
One can then enforce the fact that left and right fundamental generators $O_f$ should commute and that they should commute with their complex conjugates (this does not imply that the algebra ${\mathcal O}$ is commutative (see remark in Section~\ref{commutativityofO}).
 However, some free parameters can still remain. In order to determine their value, we proceed as follows.
First of all, we notice that  fusion matrices $N_1$, $N_2$ and $N_3$ obey non-trivial polynomial relations  (see below)  ref\/lecting the fact that the fusion ring ${\mathcal A}_k({\rm SU}(4))$ is a quotient of the representation ring of ${\rm SU}(4)$. Since ${\mathcal O}$ and ${\mathcal E}_k$ are modules over the fusion ring, the same relations have to be satisf\/ied by the corresponding generators $O_{f^{L,R}}$ and $G_f$.
In general these equations allow us to determine the remaining parameters but it may be (see for instance our discussion of the ${\mathcal E}_8$ case in Section~\ref{sectionE8})  that the f\/inal solution,  after that last step, is not unique;  however, in our cases, this ref\/lects the existence of possible automorphisms\footnote{These are  permutations $\pi$ on vertices of ${\mathcal O}$ such that for all vertices,  $(\pi(x),\pi(y))$ is an edge if\/f $(x,y)$ is an edge.} of the graph of quantum symmetries.

 In the case of ${\rm SU}(4)$ at level $k$, one can use three non-trivial polynomial relations $y_{k+1}=0$, $y_{k+2}=0$, $y_{k+3}=0$ expressing the fact that irreducible representations associated with weights $\{k+1,0,0\}$, $\{k+2,0,0\}$ and $\{k+3,0,0\}$ do not exist in ${\mathcal A}_k({\rm SU}(4))$ (in terms of quantum groups at roots of unity, they correspond to representations with vanishing quantum dimension). Setting~$x_1$ for $\{100\}$,  $x_2$ for $\{010\}$,  $x_3$ for $\{001\}$, the polynomial $y_s$ can be expressed as the determinant of a square matrix $s\times s$, whose line number~$j$ is given by the vector  ${\ldots, 0, 1, x_1, x_2, x_3, 1, 0, \ldots }$, which should be truncated in such a way that $x_1$ belongs to the diagonal (for instance, line number~$1$ of~$y_6$ is $(x_1,x_2,x_3,1,0,0)$, line number $6$ of $y_7$ is  $(0, 0, 0, 0, 1, x_1, x_2)$, etc).
This property (Giambelli formula) is a consequence of the Littlewood--Richardson rule. Remark: One can always eliminate~$x_2$ between $y_{k+1}=0$, $y_{k+2}=0$, $y_{k+3}=0$  and express~$x_3$ as a (rational) polynomial in~$x_1$; one can instead eliminate $x_3$ and f\/ind a polynomial relation between~$x_1$ and~$x_2$ but one cannot express polynomially $x_2$ in terms of $x_1$ (it is known~\cite{YellowBook} that this is never possible for a fusion ring of ${\rm SU}(g)$ when $g$ and the chosen level $k$ are both even). We therefore use the vanishing of $y_5$, $y_6$, $y_7$ for ${\mathcal E}_4$,  of $y_7$, $y_8$, $y_9$ for ${\mathcal E}_6$ and of $y_{9}$, $y_{10}$, $y_{11}$ for ${\mathcal E}_8$ as a tool to determine the remaining parameters.
In the case of ${\mathcal E}_4$ for instance,  these polynomial relations read as follows:
\begin{gather*}
\mbox{\tiny $y_5 = {x_1}^5-4 {x_2} {x_1}^3+3 {x_3} {x_1}^2+3 {x_2}^2
   {x_1}-2 {x_1}-2 {x_2} {x_3} = 0,$} \\[-1ex]
\mbox{\tiny $y_6 =  {x_1}^6-5 {x_2} {x_1}^4+4 {x_3} {x_1}^3+6 {x_2}^2
   {x_1}^2-3 {x_1}^2-6 {x_2} {x_3}
   {x_1}-{x_2}^3+{x_3}^2+2 {x_2} = 0,$}\\[-1ex]
\mbox{\tiny $y_7  =   {x_1}^7-6 {x_2} {x_1}^5+5 {x_3} {x_1}^4+10
   {x_2}^2 {x_1}^3-4 {x_1}^3-12 {x_2} {x_3}
   {x_1}^2-4 {x_2}^3 {x_1}+3 {x_3}^2 {x_1}+6
   {x_2} {x_1}+3 {x_2}^2 {x_3}-2 {x_3} = 0.$}
\end{gather*}
Eliminating for example $x_3$, one f\/inds that  $ {x_1} {x_2}$ should be equal to
\begin{gather*}
\mbox{\tiny $\frac{1}{55611516017584128}({x_1}^3 (1572913848761 \,
   {x_1}^{28}-101219273794784\,  {x_1}^{24}+1519972607520288\,
   {x_1}^{20}   - 10071512027614400\,  {x_1}^{16}$}   \\
 \qquad\mbox{\tiny ${} -12849609824079344\,
   {x_1}^{12}   +189817789697417216\,  {x_1}^8-183010445962251264\,
   {x_1}^4-16377652617161728)).$}
\end{gather*}

\looseness=-1
Once the matrices describing the fundamental generators have been fully determined, up to possible graph automorphisms, we want to be able of giving explicitly  the permutations describing the complex conjugation, the chiral transposition and the chiral adjoint operations (the last one being the composition of the f\/irst two). We remember, however, that there is still some freedom  in the determination of these permutations. For example if $P$ is a matrix implementing the chiral transposition (so $O_{f^R}= P O_{f^L} P^{-1}$), and if  $U$ is a permutation matrix commuting both with  ${O}_{{f}^L}$ and  ${O}_{{f}^R}$ and such that $U \tilde U = 1$, we  f\/ind another acceptable ``chiral matrix''  $P^\prime$ by setting $P^\prime = P \, U$.
We shall restrict as follows the possible choices for $P$: whenever  $x \neq y$ are two distinct vertices of the graph ${\mathcal O}$ for which the associated toric matrices are both symmetric and equal, we decide that $x$ and $y$ should be f\/ixed (rather than interchanged) by the operation $c$.

From the knowledge of the six chiral generators, we can draw the two chiral subgraphs making the
Ocneanu graph of quantum symmetries. There are at least three ways to draw such a graph. The f\/irst one uses $r_O$ vertices and one type of line for each chiral generator; this is still readable in the ${\rm SU}(2)$ situation but not in our case, where we have six types of lines (actually four: two oriented ones, and two unoriented ones); another method (see the article \cite{RobertGil:E4SU4}  as an example) draws only the left graph that describes multiplication of an arbitrary vertex by a chiral left generator; chiral conjugated vertices are then related by a dashed line so that multiplication by chiral right generators is obtained by conjugating the left multiplication. In this paper, we shall use a third solution, that we f\/ind more readable: we only display the graphs describing the multiplication by left generators (see Figs.~\ref{grOc-E4},~\ref{grOc-E6} and~\ref{grOc-E8}), with some arbitrary labeling of the vertices, but we give, for each case, the permutation describing the chiral adjoint operation.
This allows the reader to obtain easily the multiplication by the bar-conjugated of the right generators, from the relations $O_{\overline{f}^R} = Q^{-1} \, O_{f^L} \, Q$
where $Q$ is the matrix implementing the chiral adjoint operation.
\label{solvingintertwining+}

{\bf About 4-ality.}
We have  $\mathbb{Z}_4$ grading $\tau$ (4-ality) def\/ined on the set of irreps, such that $\tau(\overline{\lambda}) = -\tau(\lambda) \mod 4$ given by
$\tau(\lambda_1,\lambda_2,\lambda_3) = \lambda_1+ 2 \lambda_2 + 3 \lambda_3 \mod 4$. It  is also obtained from the corresponding Young tableau by calculating the number of boxes modulo $4$. This 4-ality def\/ined on vertices of ${\mathcal A}_k$ induces a $\mathbb{Z}_4$ grading in the modules  ${\mathcal E}_k $, and in ${\mathcal O}$. It will be used to display their corresponding graphs.

{\bf Determination of the quantum graph $\boldsymbol{{\mathcal E}}$.} In all three cases, it is obtained  as one particular  component of the left (or right) graph of quantum symmetries ${\mathcal O}$,  where it coincides with the left (or right) chiral subgraph (this property is not generic but  holds for those quantum graphs obtained from direct\footnote{i.e., not followed by a contraction with respect to some simple component of the possibly non simple group~$K$ under study.} conformal embedding). Other components of ${\mathcal O}$ describe other quantum graphs (that do not have self-fusion in general) but are modules for ${\mathcal E}$, and of course for ${\mathcal A}_k$ as well.
The graph ${\mathcal E}$ is obtained as the union of three graphs $G_f$ (sharing the same vertices but with dif\/ferent types of edges) def\/ined by (three) adjacency matrices also denoted $G_f$ read from the adjacency matrices $O_{f^L}$ (or $O_{f^R}$) of ${\mathcal O}$. The graph $G_f$ is connected for $f=\{100\}$ (or $\{001\}$) but not, in general, for  $\{010\}$.
The fact that ${\mathcal E}$ has self-fusion, not necessarily commutative since it is isomorphic with the chiral subalgebras, follows from the multiplicative structure of ${\mathcal O}$.

Obtaining annular matrices $F_n$ is now straightforward since they obey the same recurrence relations as the fusion matrices $N_n$ of ${\rm SU}(4)$,  but with a dif\/ferent seed, namely $F_{000}= I_{r_E}$ and  $F_f = G_f$.
We shall not give explicitly these matrices for reasons of size (only $G_f$ will be given), but  the fact that their calculated matrix elements turn out to be  non-negative integers, as they should,  provides a compatibility check of the previous determination of the quantum graphs: indeed, any mistake in one of the adjacency matrices $G_f$ usually induces  the appearance of some negative integer coef\/f\/icients in one or several of the $F_n$'s.

{\bf What else is to be found, or not to be found, in the coming sections.} We have determined the toric structure (i.e., all toric matrices
$W_{x0}$)  for all three cases, using the modular splitting equation. This was a necessary step towards the determination of chiral generators for the graph of quantum symmetries. However, displaying for instance these $192$ matrices of size $165 \times 165$ (the case of ${\mathcal E}_8({\rm SU}(4)$) in a printed form is out of question. In order to keep the size of this paper reasonable, we shall only describe the structure of the chiral generators, by displaying the graphs of $O_{f^{L}}$ and the permutation $P$ that implements chiral transposition and allows one to reconstruct the graphs of $O_{f^{R}}$.  Matrices $O_{f^{L,R}}$ are adjacency matrices of those graphs.
We shall not describe the full multiplicative structure of ${\mathcal O}$ in terms of linear generators; this was done for~${\mathcal E}_4$ in~\cite{RobertGil:E4SU4}. For the same reason we shall not give the dual annular matrices $S_x$.
Once the quantum graph itself is known (adjacency matrices $G_f$), it is possible to ``reverse the machine'' and  realize explicitly the algebra ${\mathcal O}$ in terms of the algebra ${\mathcal E}$: it is a particular quotient of its tensor square.  Using then the annular matrices $F_{n}$ and the  realization of generators of ${\mathcal O}$ as tensor products,  there is a way to check that our determination of toric matrices $W_{x0}$ was indeed correct. This was  done explicitly~\cite{RobertGil:E4SU4} in the case of ${\mathcal E}_4$, and can be done for the other graphs along the same lines. This analysis will not be repeated here.

\looseness=-1
Along general lines discussed in \cite{Ocneanu:paths}, one can associate a quantum groupo\"\i d ${\mathcal B}$ to every quantum graph ${\mathcal E}$. More precisely,  ${\mathcal B}$ is  a f\/inite dimensional weak Hopf algebra which is simple and co-semisimple.
One can think of the algebra ${\mathcal B}$ as a direct sum of $r_A$ matrix simple components,  and of its dual, the algebra $\widehat {\mathcal B}$, as a sum of $r_O$ matrix simple components. The dimensions $d_n$ (and $d_x$) of these blocks, called horizontal or vertical dimensions, or dimensions of generalized spaces of essential paths, or spaces of admissible triples or generalized triangles, etc.,  can be obtained from the annular (or dual annular) matrices. We denote by $d_H = \sum_n d_n$ the total horizontal dimension.
We shall not provide more details about the structure of this quantum groupo\"\i d in the present paper but its total dimension $d_{\mathcal B} = \sum_n d_n^2 = \sum_x d_x^2 $ will be given in each case.

We  calculated the quantum dimensions of simple objects of ${\mathcal E}_k$  in two possible ways: using spectral properties of the adjacency matrix $F_{\{100\}}$, obtained as a by-product of the determination of the graph of quantum symmetries, and using induction/restriction from the fusion algebra~${\mathcal A}_4$. The quantum cardinality $\vert {\mathcal E}_k \vert $, already obtained at the end of the previous section, is then recovered by summing the squares of these quantum dimensions. This provides a non trivial check of the calculations.

Real-ambichiral partition functions:
As it was recalled already, all vertices of an Ocneanu graph are associated with partition functions. Among them, only one ($Z_1$, associated with the origin) is modular invariant: it commutes with $s$ and $t$.  The others are not, although they all commute with $s^{-1}   t   s$. It would be rather heavy to give tables for all of them, and the reader can certainly obtain these results by using the provided information (they can also be  obtained  from the authors, if needed). However we shall give explicit expressions for partition functions associated with those vertices that are both ambichiral (i.e., $x$ is such that it belongs to the f\/irst connected component of the graph ${\mathcal O}$ and such that $x=x^\dag$)  and real  (i.e., $x=\overline{x}$).  There are only three vertices of that type for ${\mathcal E}_4$, two for ${\mathcal E}_6$ (not ten\footnote{They coincide with ambichiral vertices for ${\mathcal E}_4$ and ${\mathcal E}_8$, but not for ${\mathcal E}_6$.}) and four for~${\mathcal E}_8$.

In some cases, exceptional modules can be found among  the connected components of the graph of quantum symmetries of a quantum graph with self-fusion.
They provide new quantum graphs, in general without self-fusion, and they can be themselves associated with modular invariant partition functions (they may be new or not).
At this point we also discuss possible conjugate invariants by using  the permutation matrix $\mathfrak{C}$
of size $r_A \times r_A$ that intertwines
representations $n$ and $\overline n$ of ${\mathcal A}_k({\rm SU}(4))$.
This matrix can also be considered as the modular invariant matrix associated with the conjugated quantum graph ${{\mathcal A}_k}^{\mathfrak c}$.
See our discussion in the dif\/ferent cases.

\subsection[Tables for ${\mathcal E}_k$]{Tables for $\boldsymbol{{\mathcal E}_k}$} \label{tableofdimensions}

From the modular invariant, we read immediately the following:
\begin{gather*}
  r_A = (k+1)(k+2)(k+3)/3!,  \qquad  r_E = {\rm Tr}({\mathcal M} ), \\
  r_O = {\rm Tr}({\mathcal M}^\dag {\mathcal M}),\qquad   r_W =  \#\{(i,j) / {\mathcal M}_{ij} \neq 0\}.
\end{gather*}

We gather in the following table the values of $r_A$ (number of i-irreps i.e., number of vertices of the graph  ${\mathcal A}_k$ ), $r_E$ (number of vertices of the graph ${\mathcal E}_k$), $r_O$ (number of vertices of the graph of quantum symmetries), $r_W$ (rank of the family of toric matrices, in general $r_W < r_O$), $\nu({\mathcal K} _{m n})$  (number of distinct norms for the matrices ${\mathcal K} _{m n}$ relative to the equations of modular splitting), $d_H$ (total horizontal dimension),  $d _{{\mathcal B}}$ (dimension of the associated quantum groupoid), and quantum cardinalities $\vert {\mathcal E}_k \vert$.
\[\begin{array}{c|cccccccc}
 {}  & r_A & r_E & r_O & r_W & \nu({\mathcal K} _{m n}) & d_H & d _{{\mathcal B}} & \vert {\mathcal E}_k \vert \\
 \hline
 {\mathcal E}_4   & 35 & 12 & 48 & 33 & 17 &  2^5 7^2  &  2^5 2713^1 & 16 (2 + \sqrt{2})  \\
  {\mathcal E}_6  &  84 & 32 & 112 & 100 & 46 & 2^6 5^1 73^1 &  2^9 5^2 13^1 53^1 & 40 (5 + 2 \sqrt{5}) \\
 {\mathcal E}_8  & 165 & 24 & 192 & 144 & 142 &  2^6 7^1 179^1 & 2^{13} 6997^1 & 48 (9 + 5 \sqrt{3})
\end{array}\]
Specif\/ic details concerning the dif\/ferent cases are given in the following subsections.

\subsection[${\rm Oc}({\mathcal E}_4({\rm SU}(4))$]{$\boldsymbol{{\rm Oc}({\mathcal E}_4({\rm SU}(4))}$}

   We consider the matrices ${\mathcal K} _{m n}$, keeping only those that are distinct,
   corresponding to each one of the possible norms\footnote{$u=1, 2, 3,
   4, 5, 6, 8, 10, 12, 13, 15, 16, 18, 32, 40, 48, 128$.} $u$. For
   instance there are $8$ (distinct) matrices in norm 1 (i.e.,
   $\#({\mathcal K} ^1)=8$), $11$ in norm $2$  (i.e., $\#({\mathcal K}
   ^2)=11$), then $8, 5, 6, 12, 3, 2, 4, 2, 4 ,6 ,4,2,2,4, 1$ of them
   for the next possible norms. A f\/irst analysis gives immediately $8$
   toric matrices in norm $1$, therefore all elements of ${\mathcal
   K}^1$, in particular ${\mathcal K} ^1[1] = W_{00} =  {\mathcal M}$ ,
   then we f\/ind $11$ new ones in norm $2$ (with multiplicity $2$), $4$
   others in norm $3$ (with multiplicity $2$). Elements of matrices
   ${\mathcal K}$ of norm $4$ are multiple of 4, so these matrices are
   either a sum of 4 toric matrices (the same toric matrix but with
   multiplicity 4), or 2 times a toric matrix with elements multiple of
   2. As the total number of toric matrices is limited (equal to~48), we
   select the second possibility, and therefore we f\/ind 5 toric matrices
   in norm 4 (with elements multiple of~2 and multiplicity~1), then with
   the same arguments we f\/ind 4 others in norm~6 (with elements multiple
   of 2 and multiplicity 1) and f\/inally the last toric matrix in norm~8
   (again with elements multiple of 2 and multiplicity~1).
All other equations, for the $17$ possible norms,  are then satisf\/ied,  and we check that the equation of modular splitting, itself, holds.
Altogether we have therefore $18=8+5+4+1$ toric matrices with multiplicity $1$ and $15=11+4$ toric matrices with multiplicity $2$.  The total number of toric matrices is $18+15+15=48$, as it should, but the rank is only $18+15=33$, as expected.

 The next step is to solve the intertwining equations that determine the $6$ matrices expressing the generators of ${\mathcal O}$,
 using the methods described in the previous section. As the family of toric matrices is not free ($r_w<r_O$), there are some free parameters left in these matrices. Elementary considerations bring their number down to $4$ for $O_{100}^L$,  and $4$ for  $O_{010}^L$;  each of them (say $\alpha$) could a priori have values equal to  0, 1  or~2 because matrix elements such as $\alpha, 1-\alpha, 2 - \alpha$ do appear in matrices $O_f$, but requiring that polynomials $y_5$, $y_6$ and $y_7$ should vanish imposes that all of these coef\/f\/icients $\alpha$ are equal to $1$. Generators $O_f^{L,R}$ are then fully  determined.
 
 \begin{figure}[t]
\centerline{\includegraphics[width=13.5cm]{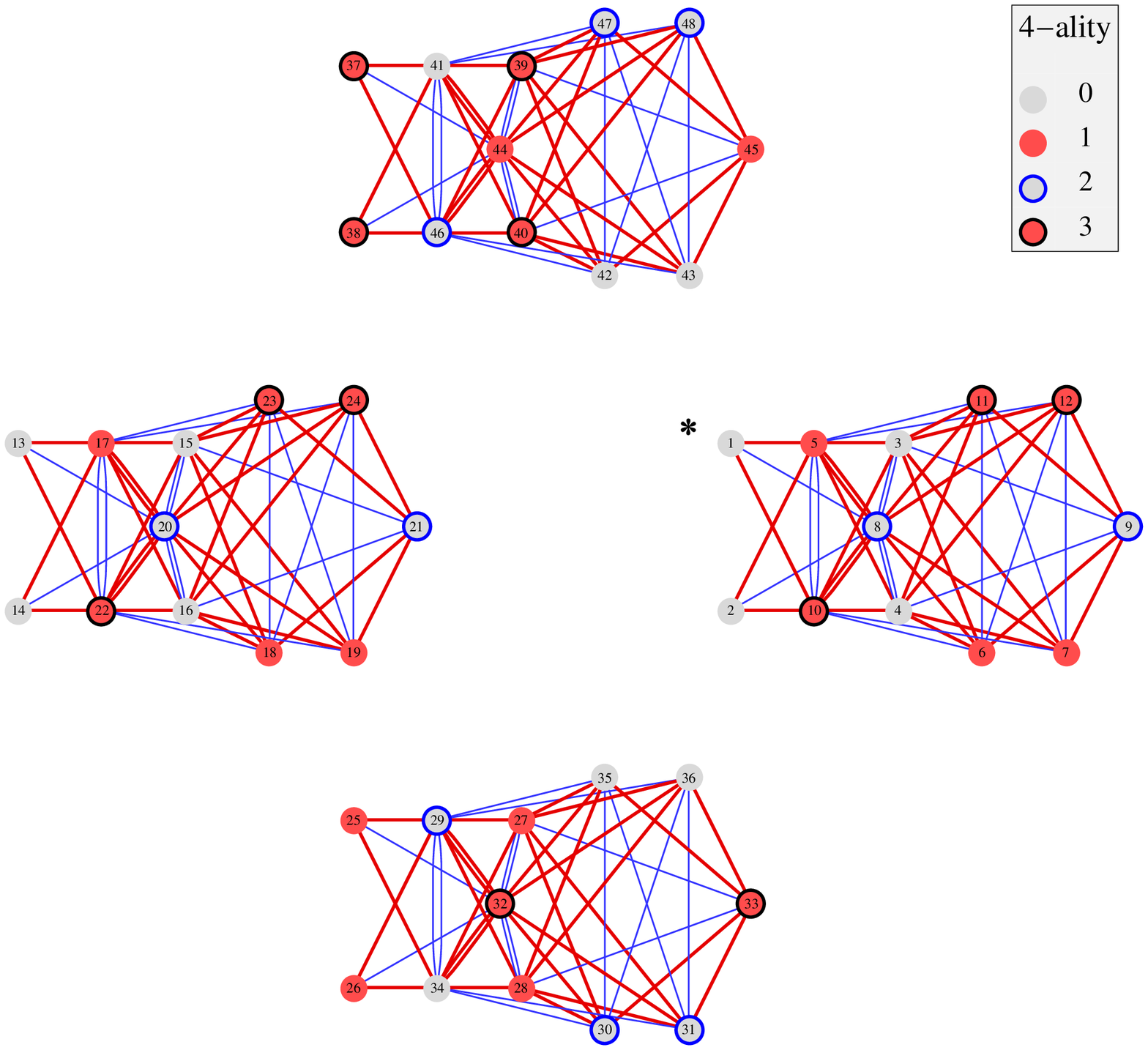}}

\caption{The left chiral graph of quantum symmetries ${\rm Oc}(\mathcal{E}_4)$.  Multiplication by the left chiral generator $100$ labeled 5 (resp.~$001$ labeled~10) is encoded by oriented red edges (thick lines), in the direction of increasing (resp. decreasing) 4-ality. Multiplication by $010$ labeled 8 is encoded by unoriented blue egdes (thin lines).}\label{grOc-E4}\vspace{-3mm}
\end{figure}

We display in Fig.~\ref{grOc-E4} the graph (with $48$ vertices)  describing the multiplication by the chiral left generators.
Multiplication by $100$ (resp.~$001$) is encoded by oriented red edges (thick lines), in the direction of increasing (resp. decreasing) 4-ality, and multiplication by $010$ is encoded by unoriented blue edges (thin lines). The identity in $\mathcal{O}$ is marked with a star on the graph. The vertex representing the fundamental left generator of type $f$ is the neighbour of the identity along the corresponding edge (of type $f$) in the left chiral graph of quantum symmetry. The chiral adjoint operation (that interchanges matrices $O_{f^{L}}$ and ${O_{{\overline f}^{R}}}$) is given by the following table, where we list only the non trivial pairs that are interchanged by this operation:{\samepage
\[\mbox{\small $\begin{array}{c|cccccccccccccccccc}
\hline
 & & \\
x              & 3  & 4  & 5  & 6  & 7  & 8   & 10 & 11 & 12 & 17 & 18 & 19  & 22 & 23 & 24 & 30 & 31  & 34  \\
x^{\dag} & 13 & 14 & 33 & 37 & 38 & 21  & 45 & 25 & 26 & 32 & 39 & 40  & 44 & 27 & 28   & 47 & 48   & 41  \\
 & & \\
\hline
\end{array}$}\]
There are 12 self-adjoint generators ($x=x^{\dag}$), ambichiral ones are 1, 2 and 9.}

{\bf Adjacency matrices of $\boldsymbol{{\mathcal E}_4}$.}
We order vertices $x$ of the left quantum graph in such a way that connected components are separated in blocks. With this choice, the left chiral generator matrices take the following block diagonal form  (see also \cite{RobertGil:E4SU4}):  $O_{100}^{L} = {\rm diag}(G_{100},G_{100},G_{100},G_{100})$, $O_{010}^{L} ={\rm diag}(G_{010},G_{010},G_{010},G_{010})$.
$G_{100}$ is obtained from the decomposition of $O_{100}^{L}$ (or of~$O_{100}^{R}$) in connected components, and
 $G_{010}$  from the decomposition of  $O_{010}^{L}$ (or of~$O_{010}^{R}$).  For the $G_f$ matrices, we order the basis elements by increasing 4-ality
\[
\mbox{\tiny $G_{100}=\left(
\begin{array}{@{\,\,}c@{\,\,}c@{\,\,}c@{\,\,}c@{\,\,}c@{\,\,}c@{\,\,}c@{\,\,}c@{\,\,}c@{\,\,}c@{\,\,}c@{\,\,}c@{\,\,}}
.& .& .& .& 1& .& .& .& .& .& .& . \\
.& .& .& .& 1& .& .& .& .& .& .& . \\
.& .& .& .& 1& 1& 1& .& .& .& .& . \\
.& .& .& .& 1& 1& 1& .& .& .& .& . \\
.& .& .& .& .& .& .& 2& .& .& .& . \\
.& .& .& .& .& .& .& 1& 1& .& .& . \\
.& .& .& .& .& .& .& 1& 1& .& .& . \\
.& .& .& .& .& .& .& .& .& 2& 1& 1 \\
.& .& .& .& .& .& .& .& .& .& 1& 1 \\
1& 1& 1& 1& .& .& .& .& .& .& .& . \\
.& .& 1& 1& .& .& .& .& .& .& .& . \\
.& .& 1& 1& .& .& .& .& .& .& .& .
\end{array}
\right),
\qquad
G_{010}=
\left(
\begin{array}{@{\,\,}c@{\,\,}c@{\,\,}c@{\,\,}c@{\,\,}c@{\,\,}c@{\,\,}c@{\,\,}c@{\,\,}c@{\,\,}c@{\,\,}c@{\,\,}c@{\,\,}}
.&.&.&.&.&.&.&1&.&.&.&.\\
.&.&.&.&.&.&.&1&.&.&.&.\\
.&.&.&.&.&.&.&2&1&.&.&.\\
.&.&.&.&.&.&.&2&1&.&.&.\\
.&.&.&.&.&.&.&.&.&2&1&1\\
.&.&.&.&.&.&.&.&.&1&1&1\\
.&.&.&.&.&.&.&.&.&1&1&1\\
1&1&2&2&.&.&.&.&.&.&.&.\\
.&.&1&1&.&.&.&.&.&.&.&.\\
.&.&.&.&2&1&1&.&.&.&.&.\\
.&.&.&.&1&1&1&.&.&.&.&.\\
.&.&.&.&1&1&1&.&.&.&.&.
\end{array}
\right).$}
\]

  It is now easy to determine the annular matrices $F_n$,  the
  horizontal dimensions  $d_n$, the total
  horizontal dimension $d_H$ and the
  dimension $d_{{\mathcal B}}$ of the quantum groupo\"\i d ${\mathcal B}({\mathcal E}_4)$. One can also
  calculate the quantum dimensions of simple objects of ${\mathcal E}_4$  and check the already obtained value for the  quantum cardinality $\vert {\mathcal E}_4 \vert$. Results are summarized in the table of Section~\ref{tableofdimensions}.

{\bf Real-ambichiral  partition functions of $\boldsymbol{{\mathcal E}_4}$}  (i.e., $x=\overline{x}$,   $x=x^{\dag}$  and $x \in {\mathcal J}$). Here,  such $x$'s coincide with ambichiral ones.
 Setting $U = 004+101+121+400$, $V = 000+012+040+210$, $W=U+V$ and $X = 111$, we get
\begin{gather*}
Z_1  =   4   \overline{X} \;  X + \overline{U} U +  \overline{V}  V   \qquad \mbox{(the modular invariant partition function $\mathcal{Z}$)},\\
Z_2  =  4    \overline{X}\; X +U \overline{V} +  \overline{U}  V , \qquad
Z_{9}  =   2    \overline{X} W+ 2 \; \overline{W}  X.
\end{gather*}

{\bf No exceptional module for ${\mathcal E}_4$.}
Since the permutation matrix $\mathfrak{C}$ commutes with $s$ and $t$, we may think of considering the new modular invariant matrix ${\mathcal M}^{\mathfrak c} = \mathfrak{C} \, {\mathcal M}$.
However, in this particular case, ${\mathcal M}^{\mathfrak c} = {\mathcal M}$, and we do not discover any new invariant in this way.
Moreover  the graph ${\mathcal O}$ of quantum symmetries contains only  copies of the quantum graph ${\mathcal E}_4$.
So we do not f\/ind any exceptional module in this case.

\subsection[${\rm Oc}({\mathcal E}_6({\rm SU}(4))$]{$\boldsymbol{{\rm Oc}({\mathcal E}_6({\rm SU}(4))}$}

 We consider the matrices ${\mathcal K} _{m n}$, keeping only those that are distinct,  corresponding to each one of the possible
 norms\footnote{$u=1,2,3,4,5,6,8,10,12, \ldots, 326$.} $u$. For instance
 there are $30$ (distinct) matrices in norm 1 (i.e., $\#({\mathcal K}
 ^1)=30$), then $104$ in norm $2$  (i.e., $\#({\mathcal K} ^2)=104$), then
 $32,130,26,50,70,64,44, \ldots, 2$ of them for the other possible norms.
 A f\/irst analysis gives immediately $30$ toric matrices in norm~$1$,
 therefore all elements of ${\mathcal K}^1$, in particular ${\mathcal K}
 ^1[1] = W_{00} =  {\mathcal M}$, then we f\/ind $36$ new ones in norm
 $2$ (many cases remaining unsettled at that stage), $4$ other in norm~$3$ and $2$ in norm~$12$ (but with multiplicity~$5$), therefore a total
 of $72$ independent ones.  A ref\/ined analysis of the norm $2$ case
 gives us $28$ more toric matrices, so that we have now reached the
 expected rank $r_W=100$. We are still missing~$12=r_O-r_W$ others that
 should not be independent of those already found.  $8$ $(= 2 \times 5
 -2)$ among them are immediately obtained from the fact that those
 coming from the norm $12$ analysis had multiplicity $5$. The last four
 are harder to f\/ind and are obtained after a deeper analysis of the norm
 $2$ case: they can be expressed as dif\/ferences between a matrix
 belonging to the set  ${\mathcal K}^2$ and a toric matrix previously
 determined in our analysis of the norm $3$ case (their matrix elements
 are nevertheless non-negative integers, of course). We then verify that
 all equations for the $46$ possible norms can be satisf\/ied with the
 obtained $112$ toric matrices, and that the equation of modular
 splitting, itself, holds.

 The next step is to solve the intertwining equations that determine the $6$ matrices expres\-sing the generators of ${\mathcal O}$ using the methods described previously.  Here again there are  several free parameters still remaining after resolution of these equations, but they are subsequently determined by using non negativity and integrality of coef\/f\/icients, commutation properties of generators, and imposing that the ${\rm SU}(4)$ polynomials $y_7$, $y_8$, $y_9$ should vanish.
 
\begin{figure}[t]
\centerline{\includegraphics[width=13.8cm]{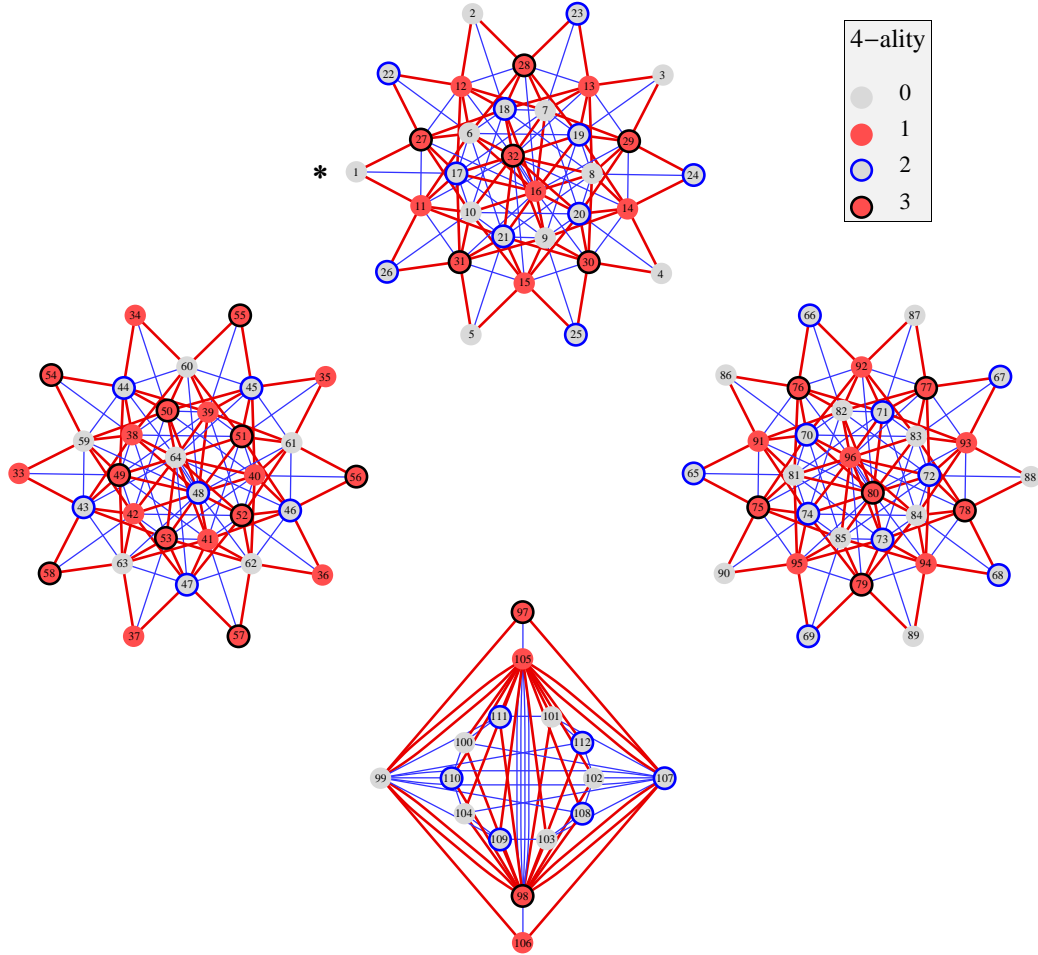}}

\caption{The left chiral graph of quantum symmetries ${\rm Oc}(\mathcal{E}_6)$.  Multiplication by the left chiral generator $100$ labeled 11 (resp.~$001$ labeled 27) is encoded by oriented red edges (thick lines), in the direction of increasing (resp.\ decreasing) 4-ality. Multiplication by $010$ labeled 17 is encoded by unoriented blue edges (thin lines).}\label{grOc-E6}\vspace{-3mm}
\end{figure}

We display in Fig.~\ref{grOc-E6} the graph  (with $112$ vertices) describing the multiplication by the chiral left generators, where we adopt the same conventions as for the $\mathcal{E}_4$ case. The graph of quantum symmetries  contains 4 connected components, the quantum graph $\mathcal{E}_6$ appears 3 times, and a~module, that we call $\mathcal{E}_6^{\mathfrak c}$, appears once.
 The chiral adjoint operation  (that interchanges matri\-ces~$O_{f^{L}}$ and ${O_{{\overline f}^{R}}}$) is given by the following table, where we list only the non trivial pairs that are interchanged by this operation,
\[\mbox{\tiny $\begin{array}{c|cccccccccccccccccccccccc}
\hline
 & & \\
x               & 6  & 7  & 8  & 9 & 10 & 11 & 12 & 13 & 14 & 15 & 16 & 17 & 18 & 19 & 20  &  21  & 27 & 28 & 29 & 30  \\
x^{\dag} & 86 & 87 & 88 & 89 & 90 & 58 & 54 & 55 & 56 & 57 & 97 & 65 & 66 & 67 & 68 &  69  & 33 & 34 & 35 & 36 \\
 & & \\
x             & 31 & 32  & 38 & 39 & 40 &  41 & 42  & 48 & 49 & 50 & 51 & 52 & 53 & 64 & 80 & 96  & & & &  \\
x^{\dag} & 37 & 106 & 76 & 77 & 78 & 79 & 75  & 107 & 91 & 92 & 93 & 94 & 95 & 99 & 105  & 98 & & & &  \\
 & & \\
\hline
\end{array}$}\]
There are 40 self-adjoint vertices ($x=x^{\dag}$), ambichiral ones are (1, 2, 3, 4, 5, 22, 23, 24, 25, 26).

{\bf Adjacency matrices of $\boldsymbol{{\mathcal E}_6}$ and $\boldsymbol{{\mathcal E}_6^{\mathfrak c}}$.}
We order vertices $x$ of the left quantum graph in such a way that connected components are separated in blocks. The left chiral gene\-ra\-tor matrices  take the following block diagonal form:  $O_{100}^{L} = {\rm diag}(G_{100},G_{100},G_{100},{G_{100}}^{\mathfrak c})$, $O_{010}^{L} ={\rm diag}(G_{010},G_{010},G_{010},{G_{010}}^{\mathfrak c})$.
Matrices $G_{100}$ and ${G_{100}}^{\mathfrak c}$ are obtained from the decomposition of $O_{100}^{L}$ (or of $O_{100}^{R}$) in connected components,
 $G_{010}$ and  ${G_{010}}^{\mathfrak c}$ from the decomposition of~$O_{010}^{L}$ (or of~$O_{010}^{R}$).  For the $G_f$ matrices, we order the basis elements by increasing 4-ality.

For the ${\mathcal E}_6$ graph we obtain:
\begin{gather*}
\mbox{\tiny $G_{100} = \left(
\begin{array}{@{\,\,}c@{\,\,}c@{\,\,}c@{\,\,}c@{\,\,}c@{\,\,}c@{\,\,}c@{\,\,}c@{\,\,}c@{\,\,}c@{\,\,}c@{\,\,}c@{\,\,}c@{\,\,}c@{\,\,}c@{\,\,}c@{\,\,}c@{\,\,}c@{\,\,}c@{\,\,}c@{\,\,}c@{\,\,}c@{\,\,}c@{\,\,}c@{\,\,}c@{\,\,}c@{\,\,}c@{\,\,}c@{\,\,}c@{\,\,}c@{\,\,}c@{\,\,}c@{\,\,}}
. & . & . & . & . & . & . & . & . & . & 1 & . & . & . & . & . & . & . & . & . & . & . & . & . & . & . & . & . & . & . & . & .\\ . & . & . & . & . & . & . & . & . & . & . & 1 & . & . & . & . & . & . & . & . & . & . & . & . & . & . & . & . & . & . & . & .\\ . & . & . & . & . & . & . & . & . & . & . & . & 1 & . & . & . & . & . & . & . & . & . & . & . & . & . & . & . & . & . & . & .\\ . & . & . & . & . & . & . & . & . & . & . & . & . & 1 & . & . & . & . & . & . & . & . & . & . & . & . & . & . & . & . & . & .\\ . & . & . & . & . & . & . & . & . & . & . & . & . & . & 1 & . & . & . & . & . & . & . & . & . & . & . & . & . & . & . & . & .\\ . & . & . & . & . & . & . & . & . & . & 1 & 1 & . & . & . & 1 & . & . & . & . & . & . & . & . & . & . & . & . & . & . & . & .\\ . & . & . & . & . & . & . & . & . & . & . & 1 & 1 & . & . & 1 & . & . & . & . & . & . & . & . & . & . & . & . & . & . & . & .\\ . & . & . & . & . & . & . & . & . & . & . & . & 1 & 1 & . & 1 & . & . & . & . & . & . & . & . & . & . & . & . & . & . & . & .\\ . & . & . & . & . & . & . & . & . & . & . & . & . & 1 & 1 & 1 & . & . & . & . & . & . & . & . & . & . & . & . & . & . & . & .\\ . & . & . & . & . & . & . & . & . & . & 1 & . & . & . & 1 & 1 & . & . & . & . & . & . & . & . & . & . & . & . & . & . & . & .\\ . & . & . & . & . & . & . & . & . & . & . & . & . & . & . & . & 1 & . & . & . & 1 & . & . & . & . & 1 & . & . & . & . & . & .\\ . & . & . & . & . & . & . & . & . & . & . & . & . & . & . & . & 1 & 1 & . & . & . & 1 & . & . & . & . & . & . & . & . & . & .\\ . & . & . & . & . & . & . & . & . & . & . & . & . & . & . & . & . & 1 & 1 & . & . & . & 1 & . & . & . & . & . & . & . & . & .\\ . & . & . & . & . & . & . & . & . & . & . & . & . & . & . & . & . & . & 1 & 1 & . & . & . & 1 & . & . & . & . & . & . & . & .\\ . & . & . & . & . & . & . & . & . & . & . & . & . & . & . & . & . & . & . & 1 & 1 & . & . & . & 1 & . & . & . & . & . & . & .\\ . & . & . & . & . & . & . & . & . & . & . & . & . & . & . & . & 1 & 1 & 1 & 1 & 1 & . & . & . & . & . & . & . & . & . & . & .\\ . & . & . & . & . & . & . & . & . & . & . & . & . & . & . & . & . & . & . & . & . & . & . & . & . & . & 1 & . & . & . & 1 & 1\\ . & . & . & . & . & . & . & . & . & . & . & . & . & . & . & . & . & . & . & . & . & . & . & . & . & . & 1 & 1 & . & . & . & 1\\ . & . & . & . & . & . & . & . & . & . & . & . & . & . & . & . & . & . & . & . & . & . & . & . & . & . & . & 1 & 1 & . & . & 1\\ . & . & . & . & . & . & . & . & . & . & . & . & . & . & . & . & . & . & . & . & . & . & . & . & . & . & . & . & 1 & 1 & . & 1\\ . & . & . & . & . & . & . & . & . & . & . & . & . & . & . & . & . & . & . & . & . & . & . & . & . & . & . & . & . & 1 & 1 & 1\\ . & . & . & . & . & . & . & . & . & . & . & . & . & . & . & . & . & . & . & . & . & . & . & . & . & . & 1 & . & . & . & . & .\\ . & . & . & . & . & . & . & . & . & . & . & . & . & . & . & . & . & . & . & . & . & . & . & . & . & . & . & 1 & . & . & . & .\\ . & . & . & . & . & . & . & . & . & . & . & . & . & . & . & . & . & . & . & . & . & . & . & . & . & . & . & . & 1 & . & . & .\\ . & . & . & . & . & . & . & . & . & . & . & . & . & . & . & . & . & . & . & . & . & . & . & . & . & . & . & . & . & 1 & . & .\\ . & . & . & . & . & . & . & . & . & . & . & . & . & . & . & . & . & . & . & . & . & . & . & . & . & . & . & . & . & . & 1 & .\\ 1 & . & . & . & . & 1 & . & . & . & 1 & . & . & . & . & . & . & . & . & . & . & . & . & . & . & . & . & . & . & . & . & . & .\\ . & 1 & . & . & . & 1 & 1 & . & . & . & . & . & . & . & . & . & . & . & . & . & . & . & . & . & . & . & . & . & . & . & . & .\\ . & . & 1 & . & . & . & 1 & 1 & . & . & . & . & . & . & . & . & . & . & . & . & . & . & . & . & . & . & . & . & . & . & . & .\\ . & . & . & 1 & . & . & . & 1 & 1 & . & . & . & . & . & . & . & . & . & . & . & . & . & . & . & . & . & . & . & . & . & . & .\\ . & . & . & . & 1 & . & . & . & 1 & 1 & . & . & . & . & . & . & . & . & . & . & . & . & . & . & . & . & . & . & . & . & . & .\\ . & . & . & . & . & 1 & 1 & 1 & 1 & 1 & . & . & . & . & . & . & . & . & . & . & . & . & . & . & . & . & . & . & . & . & . & .
\end{array}
\right)$,}
\\
\mbox{\tiny $G_{010} = \left(
\begin{array}{@{\,\,}c@{\,\,}c@{\,\,}c@{\,\,}c@{\,\,}c@{\,\,}c@{\,\,}c@{\,\,}c@{\,\,}c@{\,\,}c@{\,\,}c@{\,\,}c@{\,\,}c@{\,\,}c@{\,\,}c@{\,\,}c@{\,\,}c@{\,\,}c@{\,\,}c@{\,\,}c@{\,\,}c@{\,\,}c@{\,\,}c@{\,\,}c@{\,\,}c@{\,\,}c@{\,\,}c@{\,\,}c@{\,\,}c@{\,\,}c@{\,\,}c@{\,\,}c@{\,\,}}
. & . & . & . & . & . & . & . & . & . & . & . & . & . & . & . & 1 & . & . & . & . & . & . & . & . & . & . & . & . & . & . & .\\. & . & . & . & . & . & . & . & . & . & . & . & . & . & . & . & . & 1 & . & . & . & . & . & . & . & . & . & . & . & . & . & .\\. & . & . & . & . & . & . & . & . & . & . & . & . & . & . & . & . & . & 1 & . & . & . & . & . & . & . & . & . & . & . & . & .\\. & . & . & . & . & . & . & . & . & . & . & . & . & . & . & . & . & . & . & 1 & . & . & . & . & . & . & . & . & . & . & . & .\\. & . & . & . & . & . & . & . & . & . & . & . & . & . & . & . & . & . & . & . & 1 & . & . & . & . & . & . & . & . & . & . & .\\. & . & . & . & . & . & . & . & . & . & . & . & . & . & . & . & 1 & 1 & 1 & . & 1 & 1 & . & . & . & . & . & . & . & . & . & .\\. & . & . & . & . & . & . & . & . & . & . & . & . & . & . & . & 1 & 1 & 1 & 1 & . & . & 1 & . & . & . & . & . & . & . & . & .\\. & . & . & . & . & . & . & . & . & . & . & . & . & . & . & . & . & 1 & 1 & 1 & 1 & . & . & 1 & . & . & . & . & . & . & . & .\\. & . & . & . & . & . & . & . & . & . & . & . & . & . & . & . & 1 & . & 1 & 1 & 1 & . & . & . & 1 & . & . & . & . & . & . & .\\. & . & . & . & . & . & . & . & . & . & . & . & . & . & . & . & 1 & 1 & . & 1 & 1 & . & . & . & . & 1 & . & . & . & . & . & .\\. & . & . & . & . & . & . & . & . & . & . & . & . & . & . & . & . & . & . & . & . & . & . & . & . & . & 1 & . & . & . & 1 & 1\\. & . & . & . & . & . & . & . & . & . & . & . & . & . & . & . & . & . & . & . & . & . & . & . & . & . & 1 & 1 & . & . & . & 1\\. & . & . & . & . & . & . & . & . & . & . & . & . & . & . & . & . & . & . & . & . & . & . & . & . & . & . & 1 & 1 & . & . & 1\\. & . & . & . & . & . & . & . & . & . & . & . & . & . & . & . & . & . & . & . & . & . & . & . & . & . & . & . & 1 & 1 & . & 1\\. & . & . & . & . & . & . & . & . & . & . & . & . & . & . & . & . & . & . & . & . & . & . & . & . & . & . & . & . & 1 & 1 & 1\\. & . & . & . & . & . & . & . & . & . & . & . & . & . & . & . & . & . & . & . & . & . & . & . & . & . & 1 & 1 & 1 & 1 & 1 & 2\\1 & . & . & . & . & 1 & 1 & . & 1 & 1 & . & . & . & . & . & . & . & . & . & . & . & . & . & . & . & . & . & . & . & . & . & .\\. & 1 & . & . & . & 1 & 1 & 1 & . & 1 & . & . & . & . & . & . & . & . & . & . & . & . & . & . & . & . & . & . & . & . & . & .\\. & . & 1 & . & . & 1 & 1 & 1 & 1 & . & . & . & . & . & . & . & . & . & . & . & . & . & . & . & . & . & . & . & . & . & . & .\\. & . & . & 1 & . & . & 1 & 1 & 1 & 1 & . & . & . & . & . & . & . & . & . & . & . & . & . & . & . & . & . & . & . & . & . & .\\. & . & . & . & 1 & 1 & . & 1 & 1 & 1 & . & . & . & . & . & . & . & . & . & . & . & . & . & . & . & . & . & . & . & . & . & .\\. & . & . & . & . & 1 & . & . & . & . & . & . & . & . & . & . & . & . & . & . & . & . & . & . & . & . & . & . & . & . & . & .\\. & . & . & . & . & . & 1 & . & . & . & . & . & . & . & . & . & . & . & . & . & . & . & . & . & . & . & . & . & . & . & . & .\\. & . & . & . & . & . & . & 1 & . & . & . & . & . & . & . & . & . & . & . & . & . & . & . & . & . & . & . & . & . & . & . & .\\. & . & . & . & . & . & . & . & 1 & . & . & . & . & . & . & . & . & . & . & . & . & . & . & . & . & . & . & . & . & . & . & .\\. & . & . & . & . & . & . & . & . & 1 & . & . & . & . & . & . & . & . & . & . & . & . & . & . & . & . & . & . & . & . & . & .\\. & . & . & . & . & . & . & . & . & . & 1 & 1 & . & . & . & 1 & . & . & . & . & . & . & . & . & . & . & . & . & . & . & . & .\\. & . & . & . & . & . & . & . & . & . & . & 1 & 1 & . & . & 1 & . & . & . & . & . & . & . & . & . & . & . & . & . & . & . & .\\. & . & . & . & . & . & . & . & . & . & . & . & 1 & 1 & . & 1 & . & . & . & . & . & . & . & . & . & . & . & . & . & . & . & .\\. & . & . & . & . & . & . & . & . & . & . & . & . & 1 & 1 & 1 & . & . & . & . & . & . & . & . & . & . & . & . & . & . & . & .\\. & . & . & . & . & . & . & . & . & . & 1 & . & . & . & 1 & 1 & . & . & . & . & . & . & . & . & . & . & . & . & . & . & . & .\\. & . & . & . & . & . & . & . & . & . & 1 & 1 & 1 & 1 & 1 & 2 & . & . & . & . & . & . & . & . & . & . & . & . & . & . & . & .
\end{array}
\right)$}
\end{gather*}
and for the ${\mathcal E}_6^{\mathfrak c}$ graph:
\begin{gather*}
\mbox{\tiny $G_{100}^{\mathfrak c} = \left(
\begin{array}{@{\,\,}c@{\,\,}c@{\,\,}c@{\,\,}c@{\,\,}c@{\,\,}c@{\,\,}c@{\,\,}c@{\,\,}c@{\,\,}c@{\,\,}c@{\,\,}c@{\,\,}c@{\,\,}c@{\,\,}c@{\,\,}c@{\,\,}}
. & . & 1 & . & . & . & . & . & . & . & . & . & . & . & . & .\\. & . & 2 & 1 & 1 & 1 & 1 & 1 & . & . & . & . & . & . & . & .\\. & . & . & . & . & . & . & . & 2 & 1 & . & . & . & . & . & .\\. & . & . & . & . & . & . & . & 1 & . & . & . & . & . & . & .\\. & . & . & . & . & . & . & . & 1 & . & . & . & . & . & . & .\\. & . & . & . & . & . & . & . & 1 & . & . & . & . & . & . & .\\. & . & . & . & . & . & . & . & 1 & . & . & . & . & . & . & .\\. & . & . & . & . & . & . & . & 1 & . & . & . & . & . & . & .\\. & . & . & . & . & . & . & . & . & . & 2 & 1 & 1 & 1 & 1 & 1\\. & . & . & . & . & . & . & . & . & . & 1 & . & . & . & . & .\\1 & 2 & . & . & . & . & . & . & . & . & . & . & . & . & . & .\\. & 1 & . & . & . & . & . & . & . & . & . & . & . & . & . & .\\. & 1 & . & . & . & . & . & . & . & . & . & . & . & . & . & .\\. & 1 & . & . & . & . & . & . & . & . & . & . & . & . & . & .\\. & 1 & . & . & . & . & . & . & . & . & . & . & . & . & . & .\\. & 1 & . & . & . & . & . & . & . & . & . & . & . & . & . & .
\end{array}
\right),
\qquad
G_{010}^{\mathfrak c} = \left(
\begin{array}{@{\,\,}c@{\,\,}c@{\,\,}c@{\,\,}c@{\,\,}c@{\,\,}c@{\,\,}c@{\,\,}c@{\,\,}c@{\,\,}c@{\,\,}c@{\,\,}c@{\,\,}c@{\,\,}c@{\,\,}c@{\,\,}c@{\,\,}}
. & . & . & . & . & . & . & . & 1 & . & . & . & . & . & . & .\\. & . & . & . & . & . & . & . & 4 & 1 & . & . & . & . & . & .\\. & . & . & . & . & . & . & . & . & . & 2 & 1 & 1 & 1 & 1 & 1\\. & . & . & . & . & . & . & . & . & . & 1 & . & . & 1 & 1 & .\\. & . & . & . & . & . & . & . & . & . & 1 & . & . & . & 1 & 1\\. & . & . & . & . & . & . & . & . & . & 1 & 1 & . & . & . & 1\\. & . & . & . & . & . & . & . & . & . & 1 & 1 & 1 & . & . & .\\. & . & . & . & . & . & . & . & . & . & 1 & . & 1 & 1 & . & .\\1 & 4 & . & . & . & . & . & . & . & . & . & . & . & . & . & .\\. & 1 & . & . & . & . & . & . & . & . & . & . & . & . & . & .\\. & . & 2 & 1 & 1 & 1 & 1 & 1 & . & . & . & . & . & . & . & .\\. & . & 1 & . & . & 1 & 1 & . & . & . & . & . & . & . & . & .\\. & . & 1 & . & . & . & 1 & 1 & . & . & . & . & . & . & . & .\\. & . & 1 & 1 & . & . & . & 1 & . & . & . & . & . & . & . & .\\. & . & 1 & 1 & 1 & . & . & . & . & . & . & . & . & . & . & .\\. & . & 1 & . & 1 & 1 & . & . & . & . & . & . & . & . & . & .
\end{array}
\right)$.}
\end{gather*}

  It is now easy to determine the annular matrices $F_n$,  the
  horizontal dimensions  $d_n$, the total
  horizontal dimension $d_H$ and the
  dimension $d_{{\mathcal B}}$ of the quantum groupo\"\i d ${\mathcal B}({\mathcal E}_6)$. One can also
  calculate the quantum dimensions of simple objects of ${\mathcal E}_6$  and check the already obtained value for the  quantum cardinality $\vert {\mathcal E}_6 \vert$. Results are summarized in the table of Section~\ref{tableofdimensions}.

{\bf Real-ambichiral   partition functions of $\boldsymbol{{\mathcal E}_6}$}   (i.e., $x=\overline{x}$,  $x=x^{\dag}$ and $x \in {\mathcal J}$). Here, there are only two such $x$'s, but $10$ ambichiral vertices.
   Setting  \begin{gather*}
  U_1  =   000 + 060 + 202 + 222   ,  \qquad U_2 =  006 + 022 + 220 + 600  ,    \\
 V_1  =  042 + 200 + 212   ,   \qquad V_2=  012 + 230 + 303  ,   \qquad V_3=  030 + 103 + 321  ,  \\
 V_4= 004 + 121 + 420  ,   \qquad
 V_5  =   030 + 123 + 301 ,   \qquad V_6= 024 + 121 + 400  ,   \\ V_7 = 032 + 210 + 303   ,\qquad V_{8}=  002 + 212 + 240
    \end{gather*}
    we obtain
      \begin{gather*}
Z_1  =    \overline{U_1}U_1 + \overline{U_2} U_2 +  \overline{V_1}  V_1 + \overline{V_2}  V_2 + \overline{V_3}  V_3 + \overline{V_4}  V_4 + \overline{V_5}  V_5 +  \overline{V_6}  V_6  +  \overline{V_7}  V_7  +  \overline{V_8}  V_8,  \\
Z_{24}  =     \overline{U_1}U_2 + \overline{U_2} U_1 + \overline{V_1}  V_4 + \overline{V_4}  V_1 +
\overline{V_2}  V_5 + \overline{V_5}  V_2  + \overline{V_3}  V_7 +  \overline{V_7}  V_3  +    \overline{V_6}  V_8 +    \overline{V_8}  V_6
    \end{gather*}
and in particular we recognize the modular invariant  $\mathcal{Z}=Z_1$.

{\bf An exceptional module for $\boldsymbol{{\mathcal E}_6}$.}
\label{E6conj}
Since the permutation matrix $\mathfrak{C}$ commutes with $s$ and~$t$, we may think of considering the new modular invariant matrix ${\mathcal M}^{\mathfrak c} = \mathfrak{C} \, {\mathcal M}$.
Here, ${\mathcal M}^{\mathfrak c} \neq {\mathcal M}$ and we f\/ind a new invariant. Using the above notations,
the corresponding partition function  $\mathcal{Z}^{\mathfrak c} = \sum_{\lambda}\chi_{\lambda}\,{\mathcal{M}}^{\mathfrak c}_{\lambda\mu}\,\bar{\chi}_{\mu}$, which is of type-II,
 (the modular invariant is not a sum of blocks) reads:
\[
{Z_1}^{\mathfrak c}  =      \overline{U_1}U_1 + \overline{U_2} U_2 +  \overline{V_1}  V_8 + \overline{V_2}  V_7 + \overline{V_3}  V_5 + \overline{V_4}  V_6 + \overline{V_5}  V_3 +  \overline{V_6}  V_4  +  \overline{V_7}  V_2  +  \overline{V_8}  V_1.
\]
Its own quantum graph, denoted ${{\mathcal E}_6}^{\mathfrak c}$ appears as a module in the graph of quantum symmetries of  ${\mathcal O}({\mathcal E}_6)$. It has $16 = {\rm Tr}({\mathcal M}^{\mathfrak c}) = r_E/2$ vertices.
One can then study it directly, i.e., determine its own annular matrices, its own algebra of quantum symmetries etc.
Since the quantum graph~${{\mathcal E}_6}^{\mathfrak c}$ is known from the very beginning (its adjacency matrices ${G_f}^{\mathfrak c}$, given previously,  are obtained from  the connected components of  ${\mathcal O}({\mathcal E}_6)$ which are not of type ${\mathcal E}_6$), the analysis is much easier than for ${\mathcal E}_6$. The ${\mathcal A}$ module structure of ${{\mathcal E}_6}$ and ${{\mathcal E}_6}^{\mathfrak c}$ dif\/fer (not the same annular matrices, of course).  One f\/inds $d_{H} = \sum_n \, d_n =11456=2^6 179^1$ and $d_{\mathcal B} = \sum_n \, {d_n}^2  = 2152960 = 2^9 5^1 29^2$. Notice that ${{\mathcal E}_6}^{\mathfrak c}$ has no self-fusion.  As expected, ${\mathcal O}({\mathcal E}_6)$ and ${\mathcal O}({{\mathcal E}_6}^{\mathfrak c})$ are isomorphic algebras, but their realizations, in terms of the graph algebra of ${\mathcal E}_6$, are dif\/ferent.

\subsection[${\rm Oc}({\mathcal E}_8({\rm SU}(4))$]{$\boldsymbol{{\rm Oc}({\mathcal E}_8({\rm SU}(4))}$}
\label{sectionE8}

\looseness=-1
 We consider the matrices ${\mathcal K} _{m n}$, keeping only those that are distinct,  corresponding to each one of the possible norms\footnote{  $u=1,2,3,4,5,6,8,9,10,11,12,14,15,16, \ldots, 4096$.} $u$. For instance there are $63$ (distinct) matrices in norm 1 (i.e., $\#({\mathcal K} ^1)=63$), then $48$ in norm $2$  (i.e., $\#({\mathcal K} ^2)=48$), then $38,71,25,36,26,60,16,18,32,38,30,36,9, \ldots $ of them for the other possible norms. A f\/irst analysis gives immediately $63$ toric matrices in norm~$1$, therefore all elements of ${\mathcal K}^1$, in particular ${\mathcal K} ^1[1] = W_{00} =  {\mathcal M}$ , then we f\/ind $48$ new ones in norm~$2$ (among them, $23$ appear with multiplicity~$2$), $8$ others in norm $3$ (all of them with multiplicity~$2$), $16$ in norm~$4$ (among them, $12$ with multiplicity~$2$,  and $8$ entries remain unsettled cases at that stage), $3$~in norm~$5$ (among them, $3$ with multiplicity~$2$), $4$ in norm $6$ (no multiplicities), $1$~in norm~$7$ (with multiplicity~$2$), nothing in norms~8, 9, 10, 11, 12, 14, but $1$ in norm $15$ (with multiplicity $2$) so that we reach a total of $63+48+8+16+3+4+1+1=144=r_W$, the expected rank. After having checked that the $8$ unsettled cases  remaining in length $4$ can be reexpressed in terms of the others, we take into account the already determined multiplicities and see that $23+8+12+3+1+1=48$, so that we can easily complete the family since $144 + 48 = 192 = r_O$.
 We then verify that all equations for the $142$ possible norms can be satisf\/ied with the obtained $192$ toric matrices, and that the equation of modular splitting, itself, holds.

This case ${\mathcal E}_8$ involves huge computations, compared with the previous cases (the tensor ${\mathcal K}$ contains $40869849$ non-zero entries), but,  fortunately, the resolution of the equation of modular splitting is somehow easier than in the  ${\mathcal E}_6$ case (in particular, the determination of the $r_O -r_W$ matrices needed to complete the family).

\begin{figure}[t]
\centerline{\includegraphics[width=13.7cm]{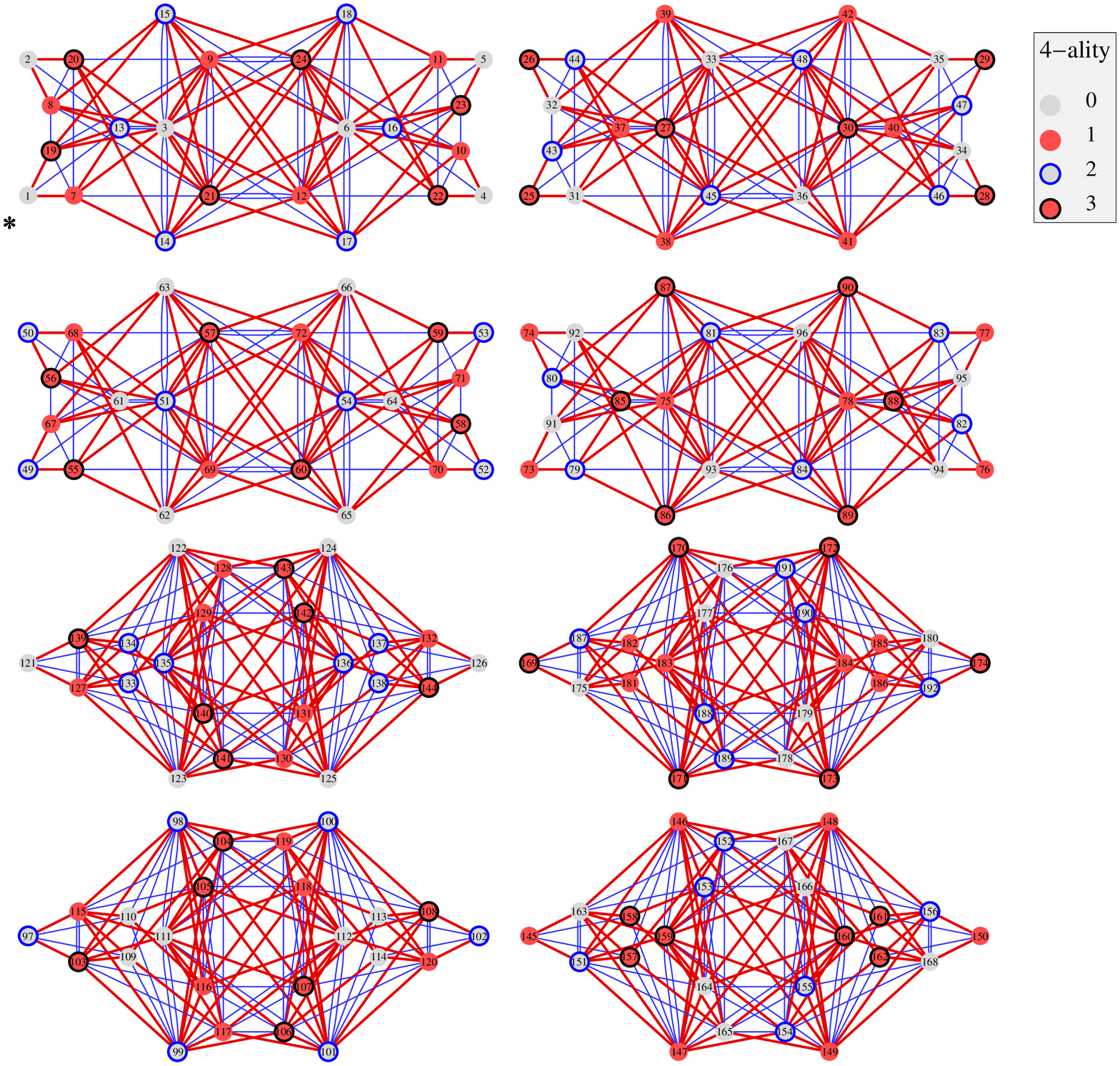}}

\caption{The left chiral graph of quantum symmetries ${\rm Oc}(\mathcal{E}_8)$.  Multiplication by the left chiral genera\-tor~$100$ labeled~7 (resp.~$001$ labeled 19) is encoded by oriented red edges (thick lines), in the direction of increasing (resp.\ decreasing) 4-ality. Multiplication by $010$ labeled 13 is encoded by unoriented blue edges (thin lines).}\label{grOc-E8}\vspace{-3mm}
\end{figure}

\looseness=-1
The next step is to solve the intertwining equations that determine the $6$ matrices expressing the fundamental generators of ${\mathcal O}$.
  Eliminating arbitrary coef\/f\/icients in those matrices is actually a~rather hard and tedious task involving simultaneously all the constraints and methods described previously (integrality, positivity, intertwining equations, commutation with complex conjugates and with chiral partners, polynomial identities in degrees 9, 10, 11).
  At the very end we obtain, up to reordering of the $192$ vertices, a single solution $O_{{010}^{L,R}}$ for the matrices of the left and right symmetric generators, a single solution for the left generator~$O_{{100}^L}$ (that we block diagonalize in the form given below) but several solutions for the corresponding right generator $O_{{100}^R}$. However, the non unicity of the solutions for the pair of fundamental generators
$(O_{{100}^L}, O_{{100}^R})$ ref\/lects the existence of graph automorphisms (see our discussion in Section~\ref{solvingintertwining}); for instance one of them exchanges vertices $17$ with $18$, or $122$ with $123$ (see Fig.~\ref{grOc-E8}). For each of these choices (for def\/initeness we select one of these equivalent solutions and just call it ${O}_{{100}^R}$) one can f\/ind several distinct permutation matrices  $Q$, with $Q =   \tilde Q = Q^{-1}$ such that ${O}_{{100}^R} = Q \, {O}_{{001}^L} \, Q^{-1}$. We then restrict  the possible choices for this matrix in the manner discussed at
page~\pageref{solvingintertwining+}.
The permutation associated with the chiral adjoint matrix $Q$ that interchanges matrices~$O_{f^{L}}$ and~${O_{{\overline f}^{R}}}$, up to graph isomorphisms, is given below (as before, we list only the non trivial pairs that are interchanged by this operation).
 \[\mbox{\tiny $\begin{array}{c|cccccccccccccccccccc}
\hline
 & & \\
x              & 3    & 6  & 7  & 8  & 9 & 10 & 11 & 12 & 13 & 14 & 15 & 16 & 17 & 18 & 19 & 20 & 21 & 22 & 23 & 24 \\
x^{\dag} & 121 & 126 & 25 & 26 & 174 & 28 & 29 & 169 & 97 & 49 & 50 & 102 & 52 & 53 & 73 & 74 & 150 & 76 & 77 & 145
\\[1mm]
x                & 27  & 30 & 33  & 36 & 37 & 38 & 39 & 40 & 41 & 42 & 43 & 44 & 45 & 46 & 47 & 48  & 51  & 54  & 57  & 60    \\
x^{\dag}  & 127 & 132 & 180 & 175 & 103 & 55 & 56 & 108 & 58 & 59 & 79 & 80 & 156 & 82 & 83 & 151& 135   & 136  & 184  & 183   \\[1mm]
x             & 61  & 64 & 67 & 68 & 69 & 70 & 71 & 72  & 75  & 78  &  81 & 84 & 85  & 88  & 93& 96 & 98 & 99 & 100 & 101  \\
x^{\dag}   & 111  & 112   & 86 & 87 & 160 & 89 & 90 & 159  & 139  & 144  & 192 & 187 & 115  & 120  & 168  & 163 & 133 & 134 & 138 & 137   \\[1mm]
x              & 104 & 105 & 106 & 107 & 116 & 117 & 118 & 119  & 128 & 129 & 130 & 131 & 140 & 141 & 142 & 143 & 152 & 153 & 154 & 155  \\
x^{\dag}  & 185 & 186 & 181 & 182  & 162 & 161 & 158 & 157 & 173  & 172 & 170 & 171 & 148 & 149 & 147 & 146   & 189 & 188 & 191 & 190    \\[1mm]
\hline
\end{array}$}\] 
There are $32$ self-adjoint vertices ($x=x^{\dag}$), ambichiral ones are  1, 2, 4, 5  (the toric matrices associated with vertices $4$ and $5$ are actually equal).

{\bf Adjacency matrices of $\boldsymbol{{\mathcal E}_8}$.}
We order vertices $x$ of the left quantum graph in such a way that connected components are separated in blocks. The left chiral generator matrices  take the following block diagonal form:  $O_{100}^{L} = {\rm diag}(G_{100},G_{100},G_{100},G_{100},{G_{100}}^{\mathfrak c},{G_{100}}^{\mathfrak c},{G_{100}}^{\mathfrak c},{G_{100}}^{\mathfrak c})$,
$O_{010}^{L} ={\rm diag}(G_{010},G_{010},G_{010},G_{010},{G_{010}}^{\mathfrak c},{G_{010}}^{\mathfrak c},{G_{010}}^{\mathfrak c},{G_{010}}^{\mathfrak c})$.
Matrices $G_{100}$ and ${G_{100}}^{\mathfrak c}$ are obtained from the decomposition of $O_{100}^{L}$ (or of $O_{100}^{R}$) in connected components, $G_{010}$ and  ${G_{010}}^{\mathfrak c}$ from the decomposition of  $O_{010}^{L}$ (or of $O_{010}^{R}$).  For the $G_f$ matrices, we order the basis elements by increasing 4-ality,
\begin{gather*}
\mbox{\tiny $G_{100} = \left(
\begin{array}{@{\,\,}c@{\,\,}c@{\,\,}c@{\,\,}c@{\,\,}c@{\,\,}c@{\,\,}c@{\,\,}c@{\,\,}c@{\,\,}c@{\,\,}c@{\,\,}c@{\,\,}c@{\,\,}c@{\,\,}c@{\,\,}c@{\,\,}c@{\,\,}c@{\,\,}c@{\,\,}c@{\,\,}c@{\,\,}c@{\,\,}c@{\,\,}c@{\,\,}}
  . &  . &  . &  . &  . &  . & 1 &  . &  . &  . &  . &  . &  . &  . &  . &  . &  . &  . &  . &  . &  . &  . &  . &  .\\ . &  . &  . &  . &  . &  . &  . & 1 &  . &  . &  . &  . &  . &  . &  . &  . &  . &  . &  . &  . &  . &  . &  . &  .\\ . &  . &  . &  . &  . &  . & 1 & 1 & 2 &  . &  . & 1 &  . &  . &  . &  . &  . &  . &  . &  . &  . &  . &  . &  .\\ . &  . &  . &  . &  . &  . &  . &  . &  . & 1 &  . &  . &  . &  . &  . &  . &  . &  . &  . &  . &  . &  . &  . &  .\\ . &  . &  . &  . &  . &  . &  . &  . &  . &  . & 1 &  . &  . &  . &  . &  . &  . &  . &  . &  . &  . &  . &  . &  .\\ . &  . &  . &  . &  . &  . &  . &  . & 1 & 1 & 1 & 2 &  . &  . &  . &  . &  . &  . &  . &  . &  . &  . &  . &  .\\ . &  . &  . &  . &  . &  . &  . &  . &  . &  . &  . &  . & 1 & 1 &  . &  . &  . &  . &  . &  . &  . &  . &  . &  .\\ . &  . &  . &  . &  . &  . &  . &  . &  . &  . &  . &  . & 1 &  . & 1 &  . &  . &  . &  . &  . &  . &  . &  . &  .\\ . &  . &  . &  . &  . &  . &  . &  . &  . &  . &  . &  . & 1 & 1 & 1 &  . & 1 & 1 &  . &  . &  . &  . &  . &  .\\ . &  . &  . &  . &  . &  . &  . &  . &  . &  . &  . &  . &  . &  . &  . & 1 & 1 &  . &  . &  . &  . &  . &  . &  .\\ . &  . &  . &  . &  . &  . &  . &  . &  . &  . &  . &  . &  . &  . &  . & 1 &  . & 1 &  . &  . &  . &  . &  . &  .\\ . &  . &  . &  . &  . &  . &  . &  . &  . &  . &  . &  . &  . & 1 & 1 & 1 & 1 & 1 &  . &  . &  . &  . &  . &  .\\ . &  . &  . &  . &  . &  . &  . &  . &  . &  . &  . &  . &  . &  . &  . &  . &  . &  . & 1 & 1 & 1 &  . &  . &  .\\ . &  . &  . &  . &  . &  . &  . &  . &  . &  . &  . &  . &  . &  . &  . &  . &  . &  . &  . & 1 & 1 &  . &  . & 1\\ . &  . &  . &  . &  . &  . &  . &  . &  . &  . &  . &  . &  . &  . &  . &  . &  . &  . & 1 &  . & 1 &  . &  . & 1\\ . &  . &  . &  . &  . &  . &  . &  . &  . &  . &  . &  . &  . &  . &  . &  . &  . &  . &  . &  . &  . & 1 & 1 & 1\\ . &  . &  . &  . &  . &  . &  . &  . &  . &  . &  . &  . &  . &  . &  . &  . &  . &  . &  . &  . & 1 &  . & 1 & 1\\ . &  . &  . &  . &  . &  . &  . &  . &  . &  . &  . &  . &  . &  . &  . &  . &  . &  . &  . &  . & 1 & 1 &  . & 1\\1 &  . & 1 &  . &  . &  . &  . &  . &  . &  . &  . &  . &  . &  . &  . &  . &  . &  . &  . &  . &  . &  . &  . &  .\\ . & 1 & 1 &  . &  . &  . &  . &  . &  . &  . &  . &  . &  . &  . &  . &  . &  . &  . &  . &  . &  . &  . &  . &  .\\ . &  . & 2 &  . &  . & 1 &  . &  . &  . &  . &  . &  . &  . &  . &  . &  . &  . &  . &  . &  . &  . &  . &  . &  .\\ . &  . &  . & 1 &  . & 1 &  . &  . &  . &  . &  . &  . &  . &  . &  . &  . &  . &  . &  . &  . &  . &  . &  . &  .\\ . &  . &  . &  . & 1 & 1 &  . &  . &  . &  . &  . &  . &  . &  . &  . &  . &  . &  . &  . &  . &  . &  . &  . &  .\\ . &  . & 1 &  . &  . & 2 &  . &  . &  . &  . &  . &  . &  . &  . &  . &  . &  . &  . &  . &  . &  . &  . &  . &  .
\end{array}
\right),$}
\qquad
\mbox{\tiny $G_{010} = \left(
\begin{array}{@{\,\,}c@{\,\,}c@{\,\,}c@{\,\,}c@{\,\,}c@{\,\,}c@{\,\,}c@{\,\,}c@{\,\,}c@{\,\,}c@{\,\,}c@{\,\,}c@{\,\,}c@{\,\,}c@{\,\,}c@{\,\,}c@{\,\,}c@{\,\,}c@{\,\,}c@{\,\,}c@{\,\,}c@{\,\,}c@{\,\,}c@{\,\,}c@{\,\,}}
 . & . & . & . & . & . & . & . & . & . & . & . & 1 & . & . & . & . & . & . & . & . & . & . & .\\. & . & . & . & . & . & . & . & . & . & . & . & 1 & . & . & . & . & . & . & . & . & . & . & .\\. & . & . & . & . & . & . & . & . & . & . & . & 2 & 2 & 2 & . & 1 & 1 & . & . & . & . & . & .\\. & . & . & . & . & . & . & . & . & . & . & . & . & . & . & 1 & . & . & . & . & . & . & . & .\\. & . & . & . & . & . & . & . & . & . & . & . & . & . & . & 1 & . & . & . & . & . & . & . & .\\. & . & . & . & . & . & . & . & . & . & . & . & . & 1 & 1 & 2 & 2 & 2 & . & . & . & . & . & .\\. & . & . & . & . & . & . & . & . & . & . & . & . & . & . & . & . & . & 1 & 1 & 1 & . & . & .\\. & . & . & . & . & . & . & . & . & . & . & . & . & . & . & . & . & . & 1 & 1 & 1 & . & . & .\\. & . & . & . & . & . & . & . & . & . & . & . & . & . & . & . & . & . & 1 & 1 & 2 & . & . & 2\\. & . & . & . & . & . & . & . & . & . & . & . & . & . & . & . & . & . & . & . & . & 1 & 1 & 1\\. & . & . & . & . & . & . & . & . & . & . & . & . & . & . & . & . & . & . & . & . & 1 & 1 & 1\\. & . & . & . & . & . & . & . & . & . & . & . & . & . & . & . & . & . & . & . & 2 & 1 & 1 & 2\\1 & 1 & 2 & . & . & . & . & . & . & . & . & . & . & . & . & . & . & . & . & . & . & . & . & .\\. & . & 2 & . & . & 1 & . & . & . & . & . & . & . & . & . & . & . & . & . & . & . & . & . & .\\. & . & 2 & . & . & 1 & . & . & . & . & . & . & . & . & . & . & . & . & . & . & . & . & . & .\\. & . & . & 1 & 1 & 2 & . & . & . & . & . & . & . & . & . & . & . & . & . & . & . & . & . & .\\. & . & 1 & . & . & 2 & . & . & . & . & . & . & . & . & . & . & . & . & . & . & . & . & . & .\\. & . & 1 & . & . & 2 & . & . & . & . & . & . & . & . & . & . & . & . & . & . & . & . & . & .\\. & . & . & . & . & . & 1 & 1 & 1 & . & . & . & . & . & . & . & . & . & . & . & . & . & . & .\\. & . & . & . & . & . & 1 & 1 & 1 & . & . & . & . & . & . & . & . & . & . & . & . & . & . & .\\. & . & . & . & . & . & 1 & 1 & 2 & . & . & 2 & . & . & . & . & . & . & . & . & . & . & . & .\\. & . & . & . & . & . & . & . & . & 1 & 1 & 1 & . & . & . & . & . & . & . & . & . & . & . & .\\. & . & . & . & . & . & . & . & . & 1 & 1 & 1 & . & . & . & . & . & . & . & . & . & . & . & .\\. & . & . & . & . & . & . & . & 2 & 1 & 1 & 2 & . & . & . & . & . & . & . & . & . & . & . & .
\end{array}
\right)$}
\end{gather*}
and for the ${\mathcal E}_8^{\mathfrak c}$:
\begin{gather*}
\mbox{\tiny $G_{100}^{\mathfrak c} = \left(
\begin{array}{@{\,\,}c@{\,\,}c@{\,\,}c@{\,\,}c@{\,\,}c@{\,\,}c@{\,\,}c@{\,\,}c@{\,\,}c@{\,\,}c@{\,\,}c@{\,\,}c@{\,\,}c@{\,\,}c@{\,\,}c@{\,\,}c@{\,\,}c@{\,\,}c@{\,\,}c@{\,\,}c@{\,\,}c@{\,\,}c@{\,\,}c@{\,\,}c@{\,\,}}
  . &  . &  . &  . &  . &  . & 1 &  . &  . &  . &  . &  . &  . &  . &  . &  . &  . &  . &  . &  . &  . &  . &  . &  .\\ . &  . &  . &  . &  . &  . & 1 & 1 & 1 &  . & 1 &  . &  . &  . &  . &  . &  . &  . &  . &  . &  . &  . &  . &  .\\ . &  . &  . &  . &  . &  . & 1 & 1 & 1 & 1 &  . &  . &  . &  . &  . &  . &  . &  . &  . &  . &  . &  . &  . &  .\\ . &  . &  . &  . &  . &  . &  . & 1 &  . & 1 & 1 & 1 &  . &  . &  . &  . &  . &  . &  . &  . &  . &  . &  . &  .\\ . &  . &  . &  . &  . &  . &  . &  . & 1 & 1 & 1 & 1 &  . &  . &  . &  . &  . &  . &  . &  . &  . &  . &  . &  .\\ . &  . &  . &  . &  . &  . &  . &  . &  . &  . &  . & 1 &  . &  . &  . &  . &  . &  . &  . &  . &  . &  . &  . &  .\\ . &  . &  . &  . &  . &  . &  . &  . &  . &  . &  . &  . & 1 & 1 & 1 &  . &  . &  . &  . &  . &  . &  . &  . &  .\\ . &  . &  . &  . &  . &  . &  . &  . &  . &  . &  . &  . & 1 &  . & 1 & 1 &  . &  . &  . &  . &  . &  . &  . &  .\\ . &  . &  . &  . &  . &  . &  . &  . &  . &  . &  . &  . &  . & 1 & 1 & 1 &  . &  . &  . &  . &  . &  . &  . &  .\\ . &  . &  . &  . &  . &  . &  . &  . &  . &  . &  . &  . &  . &  . & 1 & 1 &  . & 1 &  . &  . &  . &  . &  . &  .\\ . &  . &  . &  . &  . &  . &  . &  . &  . &  . &  . &  . &  . &  . & 1 & 1 & 1 &  . &  . &  . &  . &  . &  . &  .\\ . &  . &  . &  . &  . &  . &  . &  . &  . &  . &  . &  . &  . &  . &  . & 1 & 1 & 1 &  . &  . &  . &  . &  . &  .\\ . &  . &  . &  . &  . &  . &  . &  . &  . &  . &  . &  . &  . &  . &  . &  . &  . &  . & 1 &  . & 1 &  . &  . &  .\\ . &  . &  . &  . &  . &  . &  . &  . &  . &  . &  . &  . &  . &  . &  . &  . &  . &  . & 1 & 1 &  . &  . &  . &  .\\ . &  . &  . &  . &  . &  . &  . &  . &  . &  . &  . &  . &  . &  . &  . &  . &  . &  . & 1 & 1 & 1 & 1 & 1 &  .\\ . &  . &  . &  . &  . &  . &  . &  . &  . &  . &  . &  . &  . &  . &  . &  . &  . &  . &  . & 1 & 1 & 1 & 1 & 1\\ . &  . &  . &  . &  . &  . &  . &  . &  . &  . &  . &  . &  . &  . &  . &  . &  . &  . &  . &  . &  . & 1 &  . & 1\\ . &  . &  . &  . &  . &  . &  . &  . &  . &  . &  . &  . &  . &  . &  . &  . &  . &  . &  . &  . &  . &  . & 1 & 1\\1 & 1 & 1 &  . &  . &  . &  . &  . &  . &  . &  . &  . &  . &  . &  . &  . &  . &  . &  . &  . &  . &  . &  . &  .\\ . & 1 & 1 & 1 &  . &  . &  . &  . &  . &  . &  . &  . &  . &  . &  . &  . &  . &  . &  . &  . &  . &  . &  . &  .\\ . & 1 & 1 &  . & 1 &  . &  . &  . &  . &  . &  . &  . &  . &  . &  . &  . &  . &  . &  . &  . &  . &  . &  . &  .\\ . &  . & 1 & 1 & 1 &  . &  . &  . &  . &  . &  . &  . &  . &  . &  . &  . &  . &  . &  . &  . &  . &  . &  . &  .\\ . & 1 &  . & 1 & 1 &  . &  . &  . &  . &  . &  . &  . &  . &  . &  . &  . &  . &  . &  . &  . &  . &  . &  . &  .\\ . &  . &  . & 1 & 1 & 1 &  . &  . &  . &  . &  . &  . &  . &  . &  . &  . &  . &  . &  . &  . &  . &  . &  . &  .
\end{array}
\right),$}\qquad
\mbox{\tiny $G_{010}^{\mathfrak c} = \left(
\begin{array}{@{\,\,}c@{\,\,}c@{\,\,}c@{\,\,}c@{\,\,}c@{\,\,}c@{\,\,}c@{\,\,}c@{\,\,}c@{\,\,}c@{\,\,}c@{\,\,}c@{\,\,}c@{\,\,}c@{\,\,}c@{\,\,}c@{\,\,}c@{\,\,}c@{\,\,}c@{\,\,}c@{\,\,}c@{\,\,}c@{\,\,}c@{\,\,}c@{\,\,}}
 .  & .  & .  & .  & .  & .  & .  & .  & .  & .  & .  & .  & 1  & 1  & .  & .  & .  & .  & .  & .  & .  & .  & .  & . \\.  & .  & .  & .  & .  & .  & .  & .  & .  & .  & .  & .  & 1  & 1  & 2  & 1  & .  & .  & .  & .  & .  & .  & .  & . \\.  & .  & .  & .  & .  & .  & .  & .  & .  & .  & .  & .  & 1  & 1  & 2  & 1  & .  & .  & .  & .  & .  & .  & .  & . \\.  & .  & .  & .  & .  & .  & .  & .  & .  & .  & .  & .  & .  & .  & 1  & 2  & 1  & 1  & .  & .  & .  & .  & .  & . \\.  & .  & .  & .  & .  & .  & .  & .  & .  & .  & .  & .  & .  & .  & 1  & 2  & 1  & 1  & .  & .  & .  & .  & .  & . \\.  & .  & .  & .  & .  & .  & .  & .  & .  & .  & .  & .  & .  & .  & .  & .  & 1  & 1  & .  & .  & .  & .  & .  & . \\.  & .  & .  & .  & .  & .  & .  & .  & .  & .  & .  & .  & .  & .  & .  & .  & .  & .  & 2  & 1  & 1  & .  & .  & . \\.  & .  & .  & .  & .  & .  & .  & .  & .  & .  & .  & .  & .  & .  & .  & .  & .  & .  & 1  & 1  & 1  & 1  & 1  & . \\.  & .  & .  & .  & .  & .  & .  & .  & .  & .  & .  & .  & .  & .  & .  & .  & .  & .  & 1  & 1  & 1  & 1  & 1  & . \\.  & .  & .  & .  & .  & .  & .  & .  & .  & .  & .  & .  & .  & .  & .  & .  & .  & .  & .  & 1  & 1  & 1  & 1  & 1 \\.  & .  & .  & .  & .  & .  & .  & .  & .  & .  & .  & .  & .  & .  & .  & .  & .  & .  & .  & 1  & 1  & 1  & 1  & 1 \\.  & .  & .  & .  & .  & .  & .  & .  & .  & .  & .  & .  & .  & .  & .  & .  & .  & .  & .  & .  & .  & 1  & 1  & 2 \\1  & 1  & 1  & .  & .  & .  & .  & .  & .  & .  & .  & .  & .  & .  & .  & .  & .  & .  & .  & .  & .  & .  & .  & . \\1  & 1  & 1  & .  & .  & .  & .  & .  & .  & .  & .  & .  & .  & .  & .  & .  & .  & .  & .  & .  & .  & .  & .  & . \\.  & 2  & 2  & 1  & 1  & .  & .  & .  & .  & .  & .  & .  & .  & .  & .  & .  & .  & .  & .  & .  & .  & .  & .  & . \\.  & 1  & 1  & 2  & 2  & .  & .  & .  & .  & .  & .  & .  & .  & .  & .  & .  & .  & .  & .  & .  & .  & .  & .  & . \\.  & .  & .  & 1  & 1  & 1  & .  & .  & .  & .  & .  & .  & .  & .  & .  & .  & .  & .  & .  & .  & .  & .  & .  & . \\.  & .  & .  & 1  & 1  & 1  & .  & .  & .  & .  & .  & .  & .  & .  & .  & .  & .  & .  & .  & .  & .  & .  & .  & . \\.  & .  & .  & .  & .  & .  & 2  & 1  & 1  & .  & .  & .  & .  & .  & .  & .  & .  & .  & .  & .  & .  & .  & .  & . \\.  & .  & .  & .  & .  & .  & 1  & 1  & 1  & 1  & 1  & .  & .  & .  & .  & .  & .  & .  & .  & .  & .  & .  & .  & . \\.  & .  & .  & .  & .  & .  & 1  & 1  & 1  & 1  & 1  & .  & .  & .  & .  & .  & .  & .  & .  & .  & .  & .  & .  & . \\.  & .  & .  & .  & .  & .  & .  & 1  & 1  & 1  & 1  & 1  & .  & .  & .  & .  & .  & .  & .  & .  & .  & .  & .  & . \\.  & .  & .  & .  & .  & .  & .  & 1  & 1  & 1  & 1  & 1  & .  & .  & .  & .  & .  & .  & .  & .  & .  & .  & .  & . \\.  & .  & .  & .  & .  & .  & .  & .  & .  & 1  & 1  & 2  & .  & .  & .  & .  & .  & .  & .  & .  & .  & .  & .  & .
\end{array}
\right).$}
\end{gather*}

  It is now easy to determine the annular matrices $F_n$,  the
  horizontal dimensions  $d_n$, the total
  horizontal dimension $d_H$ and the
  dimension $d_{{\mathcal B}}$ of the quantum groupo\"\i d ${\mathcal B}({\mathcal E}_8)$. One can also
  calculate the quantum dimensions of simple objects of ${\mathcal E}_8$  and check the already obtained value for the  quantum cardinality $\vert {\mathcal E}_8 \vert$. Results are summarized in the table of Section~\ref{tableofdimensions}.

{\bf Real-ambichiral partition functions of $\boldsymbol{{\mathcal E}_8}$}   (i.e., $x=\overline{x}$,  $x=x^{\dag}$  and $x \in {\mathcal J}$).
Here such $x$'s coincide with ambichiral ones.
With $U_1 =  113 + 133 + 311 + 331$, $U_2 = 020 + 032 + 060 + 206 + 230 + 303 + 323 + 602$,   $U_3 = 000 + 008 + 080 + 121 + 141 + 214 + 412 + 800$, one f\/inds:
\begin{gather*}
    Z_1  =   2   \overline{U_1}U_1 + \overline{U_2} U_2 + \overline{U_3} U_3 \qquad   \mbox{(the modular invariant  $\mathcal{Z}$)},\\
    Z_{2}  =   2   \overline{U_1}U_1 + \overline{U_2} U_3 + \overline{U_3} U_2,  \qquad
    Z_{4} = Z_{5}  =   \overline{U_1}U_2 + \overline{U_2} U_1 + \overline{U_1} U_3 + \overline{U_3} U_1.
\end{gather*}

{\bf An exceptional module for $\boldsymbol{{\mathcal E}_8}$.}
Since the permutation matrix $\mathfrak{C}$ commutes with $s$ and $t$, we may think of considering the new modular invariant matrix ${\mathcal M}^{\mathfrak c} = \mathfrak{C} \, {\mathcal M}$.
However, in this case, and like for ${\mathcal E}_4$, we f\/ind  that ${\mathcal M}^{\mathfrak c} = \mathfrak{C} \, {\mathcal M}$ is equal to  ${\mathcal M}$.
Therefore, we do not f\/ind a new modular invariant in this way.
However, the graph ${\mathcal O}$ of quantum symmetries of ${\mathcal E}_8$ contains
not only four copies of ${\mathcal E}_8$ but also four copies of a module
that we call ${{\mathcal E}_8}^{\mathfrak c}$, with $24$ vertices as well. One can then
study it directly, like we did in the previous case (its adjacency
matrices ${G_f}^{\mathfrak c}$, given previously,  are obtained from the connected
components of  ${\mathcal O}({\mathcal E}_8)$ which are not of type ${\mathcal
E}_8$). In particular, we can determine its own annular
matrices, since, although they are associated with the same invariant,
the ${\mathcal A}$ module structure of ${\mathcal E}_8$ and ${{\mathcal E}_8}^{\mathfrak c}$
dif\/fer. One can deduce, from this study, the dimensions of the dif\/ferent
blocks $d_n$ of the associated bialgebra ${\mathcal B}({{\mathcal E}_8}^{\mathfrak c})$,
for its f\/irst multiplication, and f\/ind  $d_{H} = \sum_n \, d_n = 95040 =
 2^6 3^3 5^1 11^1$ and $d_{\mathcal B} = \sum_n \, {d_n}^2  = 80547840 = 2^12  3^2 5^1 19^1 23^1$.
In the present case ${\mathcal O}({\mathcal E}_8)$ and ${\mathcal O}({{\mathcal E}_8}^{\mathfrak c})$ are identical.

\pdfbookmark[1]{Afterword}{Afterword}
\section*{Afterword}

Although this was already discussed in our introduction, we conclude this article by comparing what was already
known and what, to our knowledge, is new in the present paper.

 What was already known:
\begin{itemize}\itemsep=0pt
\item The modular invariant partition functions corresponding to the three exceptional cases discussed in this paper \cite{SchellekensYankielowicz, AltschulerBauerItzykson, AldazabalEtAl}.
\item The corresponding quantum graphs \cite{PetkovaZuberNP1995, PetkovaZuberNP1997, Ocneanu:Bariloche}.
\item The fact that every  such graph passes the self-connection test \cite{Ocneanu:Bariloche}, see also our remark on page~\pageref{remark}. Concerning this last point, we believe that checking this condition should not be necessary for those quantum subgroups obtained, as here, from direct conformal embeddings, at least in those cases where the solution, obtained after resolution of the equations of modular splitting and intertwining,  is unique (the solution should exist, and if the one obtained is unique \dots\ it is it!).
\end{itemize}

  What was not known\footnote{Notice however that a detailed presentation of the $k=4$ case is given in our paper \cite{RobertGil:E4SU4}.}:

\begin{itemize}\itemsep=0pt
\item The full resolution of the equation of modular splitting for these three cases (only the chiral part of it\footnote{The chiral equations of modular splitting is a simplif\/ied form of the full system of equations described in Section~\ref{sec:eqsofmodsplit} ref\/lecting the fact that the following equality, for $m,n \in {\mathcal A}$, $a \in {\mathcal E}$, holds:
$
(m  (n a) ) = (m \cdot n) a
$.
For cases where ${(F_p)}_{00}= {\mathcal M}_{p0}$, a condition that holds in the cases studied in this paper, this associativity constraint implies immediately
$
\sum_p ({N_n})_{mp} \,  {\mathcal M}_{p0} =  \sum_b \, {(F_n)}_{0b} \, {(F_n)}_{b0}
$.
The left hand side (that involves only the f\/irst line of the modular invariant matrix) is known and the right hand side (the annular matrices) can be determined thanks to methods analog to those used in this paper, but this is technically simpler since the previous identity  describes only $r_A^2$ equations instead of $r_A^4$. Adjacency matrices $G_f$ can then be obtained, but not the quantum symmetry matrices $O_{f^{L,R}}$.} was used to obtain the ${\rm SU}(4)$ graphs in \cite{Ocneanu:Bariloche}).
\item The structure of the algebra of quantum symmetries and the graphs of its fundamental generators for these three cases.
\end{itemize}
One can also f\/ind in the previous sections many details concerning quantum dimensions and cardinalities  for quantum graphs obtained from conformal embeddings, as well as a description of general techniques that, to our knowledge, are not discussed elsewhere.

\appendix

\pdfbookmark[1]{Appendix}{Appendix}

\section*{Appendix}

We f\/irst remind the reader what are the expressions for representatives of the generators $s$ and $t$ of the modular group in the particular case of ${\rm SU}(g)$ groups and we give recurrence formulae for the fusion matrices of ${\rm SU}(4)$. Then, we give the quantum dimensions of the spaces of sections associated with the modular blocks of the partition function, for the three exceptional cases studied in this article. The quantum dimensions of ${\rm SU}(4)$ irreducible representations needed in this calculation are obtained from the quantum version of the Weyl formula. Using the Kac--Peterson formula for modular generators, we also obtain a general expression for the quantum cardinality of ${\mathcal A}_k$ when $q^{k+4}=-1$.
Quantum dimensions for irreducible representations of $B_7\simeq {\rm Spin}(15)$, $A_9 \simeq {\rm SU}(10)$ and $D_{10} \simeq {\rm Spin}(20)$ at level $k$ can be calculated in a standard way from the quantum version of the Weyl formula; those at level $1$ have been used in the text (Section~\ref{massE}). In the case of ${\rm SU}(10)$, or of ${\rm SU}(g)$ in general, the calculation at level $1$ is very simple, and one f\/inds $\dim (n)=1$ for all $n\in {\mathcal A}_1({\rm SU}(g))$. In the case of ${\rm Spin}(15)$ we refer to a~discussion in~\cite{RobertGil:E4SU4}. Our last appendix gives some details concerning the case of ${\rm Spin}(20)$.

\section[Generators $s$ and $t$ for ${\rm SL}(2,\mathbb{Z})$]{Generators $\boldsymbol{s}$ and $\boldsymbol{t}$ for $\boldsymbol{{\rm SL}(2,\mathbb{Z})}$}

Expressions for representatives of the generators $s$ and $t$ of ${\rm SL}(2,\mathbb{Z})$ are given, for any simple Lie algebra, and for a given level $k$,  by the Kac--Peterson formulae \cite{Kac-Peter}. These general expressions, that we recall below in the case of ${\rm SU}(g)$,  involves a summation over the elements of the Weyl group, which, for ${\rm SU}(4)$,  is the symmetric group on four objects, and scalar products between fundamental weights shifted by the Weyl vector. These explicit expressions for $s$ and $t$ are matrices of size $r_A \times r_A$. Indices $m, n, \ldots $ range over all ${\rm SU}(g)$ Young tableaux (including the trivial)  corresponding to i-irreps with levels up to $k$.
As usual,  $s$ and $t$ are unitary and such that $(s\,t)^3 = s^2 = \mathcal{C}$, the ``charge matrix'' satisfying $\mathcal{C}^2=\munite$.  The diagonal $t$ matrix obeys $t^{2\,  g \,\kappa} = \munite$, where $\kappa$ is the altitude, def\/ined as $\kappa=k + g$, (in our ${\rm SU}(4)$ case, the dual Coxeter number is  $g=4$). Using these expressions, we checked that $\mathcal{M}$ indeed commutes with the modular generators $s$ and $t$.
\[
s_{mn}= \frac{i^{rg/2} \sqrt{\Delta (r)} }{ (g+k)^{r/2} }\left(\sum _{s=1}^{g!}   \epsilon_{w(s)}    e^{-\frac{2 i \pi  \langle w(s)(m+\rho ),n+\rho \rangle }{g+k}}\right)
, \qquad
t_{mn}=e^{2 i \pi \left[\frac{\langle m+\rho ,m+\rho \rangle }{2 (g+k)}-\frac{\langle \rho ,\rho \rangle }{2 g}\right]}   \delta_{mn},
\]
where $w(s)$ runs over the $g!$ permutations of the Weyl group of ${\rm SU}(g)$, $\epsilon_{w(s)}$ is its signature, $r=g-1$ is the rank, $\rho$ is the Weyl vector, and $\Delta(r)$ is the determinant of the fundamental quadratic form. For ${\rm SU}(g)$, each entry $s_{mn}$ of the $s$ matrix can be written in terms of the determinant of a $g \times g$ matrix $A(m,n)$, with matrix elements $A(m,n)_{a,b = 1,\ldots, g}$, as follows \cite{KacWakimoto}:
\begin{gather*}
 s_{mn} = \frac{i^{rg/2}}{\sqrt{g} (g+k)^{r/2}}  \det A(m,n), \qquad \hbox{with} \quad A(m,n)_{ab} = \exp[-\frac{2 i \pi}{k+g}   \phi_a(m) \phi_b(n)]
\\
\hbox{and} \quad \phi_a(m)=\sum_{s=a}^{g-1}   (m_s + 1) - \frac{1}{g}   \sum_{s=1}^{g-1}   s   (m_s+1),
\end{gather*}
where $m_s$ are the components of $m$ along the basis of fundamental weights, and with the understanding that the f\/irst sum gives $0$ if $a=g$.

\section[Recurrence formulae for fusion and annular matrices of ${\rm SU}(4)$]{Recurrence formulae for fusion and annular matrices of $\boldsymbol{{\rm SU}(4)}$}

Fusion matrices for the fundamental irreducible representations can be obtained in several ways. One possibility is to use the previous expressions for $s$ and $t$ together with the Verlinde formulae, but in the case of ${\rm SU}(g)$ groups, and in particular for ${\rm SU}(4)$,  it is much simpler to use standard Young tableaux techniques (at a f\/ixed level, the horizontal size of the tableaux is bounded) or, equivalently, the Pieri rules describing the reduction of tensor products of arbitrary irreducible representations by the fundamental ones (see for instance \cite[p.~695]{YellowBook}). The same rules can be used to obtain the recursion formulae given below.
Once the fusion matrices of the fundamental irreducible representations are known, the others can be determined from the {\it truncated} recursion formulae of ${\rm SU}(4)$ irreps, applied for increasing level $\ell$, up to $k$ ($2 \leq \ell \leq k$):
\begin{gather*}
N_{(\ell -p ,p-q,q)}  =  N_{(1,0,0)}   N_{(\ell-p-1,p-q,q)} - N_{(\ell-p-2,p-q+1,q)} - N_{(\ell-p-1,p-q-1,q+1)} \\
 \phantom{N_{(\ell -p ,p-q,q)}  =}{} -  N_{(\ell-p-1,p-q,q-1)}    \qquad   \textrm{for} \quad 0 \leq q \leq p \leq \ell - 1,   \\
N_{(0,\ell-q,q)}  =   (N_{(q,\ell-q,0)})^{\rm tr} \qquad  \textrm{for} \quad 1 \leq q \leq \ell,   \\
N_{(0,\ell,0)}  =  N_{(0,1,0)}   N_{(0,\ell-1,0)} - N_{(1,\ell-2,1)} - N_{(0,\ell-2,0)}.
\end{gather*}

Annular matrices $F_n = (F_n)_{ab}$ relative to a quantum graph ${\mathcal E}_k$ are calculated with exactly the same recurrence relations. Only the seed used to start the recurrence is dif\/ferent: $F_{(1,0,0)}$,  $F_{(0,1,0)}$ and $F_{(0,0,1)} = F_{(1,0,0)}^{\rm tr}$ are the adjacency matrices of the chosen graph.

\section{Quantum dimensions}

{\bf Quantum dimensions of the modular blocks.}
As discussed at the end of Section~\ref{section1}, every vertex $a$ of ${\mathcal E}_k$ has a dimension $\dim (a) = \dim (\Gamma_a) /\dim (\Gamma_0)$, where the dimension of  the quantum space of sections  $\Gamma_a$  is calculated by using induction from ${\mathcal A}_k({\rm SU}(4))$:  $\dim (\Gamma_a) = \bigoplus_{n  \uparrow \Gamma_a} \dim (n)$.
In particular, the dimensions of spaces of sections $\Gamma_s$ associated with the modular vertices $s$ of~${\mathcal E}_k$ (or with the various modular blocks of the partition function)  are obtained by summing the quantum dimensions of the irreducible representations that appear in the modular blocks for the dif\/ferent cases:
When $k=4$, i.e.,  $q = \exp(i \pi / 8)$,  the f\/irst two modular blocks have dimension $4 \left(2+\sqrt{2}\right)$ and the third has dimension  $4 \left(1+\sqrt{2}\right)$.
 When $k=6$, i.e.,  $q = \exp(i \pi / 10)$, the ten modular blocks have dimension $20+8 \sqrt{5}$.
When $k=8$, i.e.,  $q = \exp(i \pi / 12)$,  the four modular blocks (remember that the last one appears twice in the modular invariant) have dimension $12 \left(9+5 \sqrt{3}\right)$.

{\bf Complements and checks.}
Quantum dimensions of irreducible representations of ${\rm SU}(4)$, when $k=4,6,8$, that are needed in the previous calculation, have been calculated  in two ways: from the Perron--Frobenius eigenvector of the adjacency matrix of the graph ${\mathcal A}_k$ associated with the def\/ining representation, and  from the quantum version of the Weyl formula,
\begin{gather*}
\dim (m)= \prod_{\alpha > 0} \frac{\langle m + \rho, \alpha\rangle_q}{\langle\rho, \alpha\rangle_q}
\end{gather*} taking $q^{k+4}=-1$ at the end.
Scalar products between a chosen weight $m= \{ m_1, m_2, m_3 \}$ and  the six positive roots $\alpha$ of $A_3$ are displayed below in the half-ribbon diagram (seen as a~generalized root set  \cite{Ocneanu:MSRI}) stemming from the $A_3 = {\mathcal A}_2 ({\rm SU}(2))$ graph. We also give the Weyl vector $\rho$
\[
\mbox{\tiny $\begin{array}{ccc}
 m_1 &  & 0 \\
  & m_1 &  \\
 0 &  & m_1 \\
  & 0 &
\end{array}
\; + \;
\begin{array}{ccc}
 0 &  & 0 \\
  & m_2 &  \\
 m_2 &  & m_2 \\
  & m_2 &
\end{array}
\; + \;
\begin{array}{ccc}
 0 &  & m_3 \\
  & m_3 &  \\
 m_3 &  & 0 \\
  & 0 &
\end{array}
,
\qquad
 \rho =
\begin{array}{ccc}
 \underline{1} &  &  \underline{1} \\
  &  \underline{3} &  \\
 2 &  & 2 \\
  & 1 &
\end{array}.$}
\]

The quantum Weyl denominator is $1_q^3 2_q^2 3_q$ (equal to $2$ when $q^{8}=-1$, to $ 5+2 \sqrt{5}$ when $q^{10}=-1$ and to $ 5+3 \sqrt{3}$ when $q^{12}=-1$).
The q-dimensions of the fundamental irreducible representations have been already given in the text.

{\bf The quantum cardinality (quantum mass) $\boldsymbol{\vert {\mathcal A}_k \vert}$ of $\boldsymbol{{\rm SU}(4)}$ at level $\boldsymbol{k}$.} This quantity  is obtained by summing the square of the quantum dimensions of the $r_A$ irreducible representations. One can also use the property  $\vert {\mathcal A}_k \vert = 1/s_{00}^2$ where $s$ is the f\/irst generator of the modular group, together with the Kac--Peterson formula \cite{Kac-Peter}, therefore expressing  $s_{00}$ in terms of weighted Weyl group averages of the norm of the Weyl vector. In this way, one f\/inds:
\begin{gather*}
\vert {\mathcal A}_k({\rm SU}(4))\vert = \frac{4 (k+4)^3}{2^{16}\cos ^4\left(\frac{\pi }{k+4}\right) \left(2 \cos \left(\frac{2 \pi
   }{k+4}\right)+1\right)^2 \sin ^{12}\left(\frac{\pi }{k+4}\right)}
\end{gather*}
and one recovers in particular the given expressions directly calculated for specif\/ic values of $k$.

{\bf Quantum dimensions for $\boldsymbol{D_{10} \simeq {\rm Spin}(20)}$ at level 1.}
We already know what the integrable representations of $D_{10}$ at level $1$ are, namely the trivial, the vectorial, and the two half-spinorial, but it is nice to recover this from a direct calculation of quantum dimensions.
We use the quantum version of the Weyl formula and, like we did with $A_3$, we display the scalar products between an arbitrary  weight and the $90$ positive roots of $D_{10}$ in the half-ribbon diagram (seen as a generalized root set  \cite{Ocneanu:MSRI}) stemming from the $D_{10} = {\mathcal D}_{16} ({\rm SU}(2))$ graph. We only display the Weyl vector $\rho$ and, to its right,  the q-dimensions for the fundamental irreducible representations.

The quantum Weyl denominator is
$1_q^{10} 2_q^9 3_q^9 4_q^8 5_q^8 6_q^7 7_q^7 8_q^6 9_q^6 10_q^4 11_q^4 12_q^3 13_q^3 14_q^2 15_q^2 16_q 17_q$,
\[
\mbox{\tiny $\rho =
\begin{array}{cccccccccc}
\underline {1} &  &\underline { 1} &  & \underline {1} &  & \underline {1} &  & \underline {1} & \underline {1} \\
  & \underline {3} &  & \underline {3} &  & \underline {3} &  & \underline {4} &  & \\
 2 &  & 5 &  & 5 &  & 6 &  & 3 & 3 \\
  & 4 &  & 7 &  & 8 &  & 8 &  &  \\
 2 &  & 6 &  & 1 &  & 1 &  & 5 & 5 \\
  & 4 &  & 9 &  & 12 &  & 12 &  &  \\
 2 &  & 7 &  & 11 &  & 14 &  & 7 & 7 \\
  & 5 &  & 9 &  & 13 &  & 16 &  &  \\
 3 &  & 7 &  & 11 &  & 15 &  & 9 & 9 \\
  & 5 &  & 9 &  & 13 &  & 17 &  &  \\
 2 &  & 7 &  & 11 &  & 15 &  & 8 & 8 \\
  & 4 &  & 9 &  & 13 &  & 14 &  &  \\
 2 &  & 6 &  & 11 &  & 12 &  & 6 & 6 \\
  & 4 &  & 8 &  & 1 &  & 1 &  &  \\
 2 &  & 6 &  & 7 &  & 8 &  & 4 & 4 \\
  & 4 &  & 5 &  & 5 &  & 6 &  &  \\
 2 &  & 3 &  & 3 &  & 3 &  & 2 & 2 \\
  & 1 &  & 1 &  & 1 &  & 1 &  &
\end{array}
 , \qquad \quad
\begin{array}{c}
 \frac{10_q 18_q}{1_q 9_q} \\ \\
 \frac{10_q 16_q 19_q}{1_q 2_q 8_q} \\ \\
 \frac{10_q 14_q 18_q 19_q}{1_q 2_q 3_q 7_q} \\ \\
 \frac{10_q 12_q 17_q 18_q 19_q}{1_q 2_q 3_q 4_q 6_q} \\ \\
 \frac{10_q^2 16_q 17_q 18_q 19_q}{1_q 2_q 3_q 4_q 5_q^2} \\ \\
 \frac{8_q 10_q 15_q 16_q 17_q 18_q 19_q}{1_q 2_q 3_q 4_q^2 5_q 6_q} \\ \\
 \frac{10_q 14_q 15_q 16_q 17_q 18_q 19_q}{1_q 2_q 3_q^2 4_q 5_q 7_q} \\ \\
 \frac{10_q 13_q 14_q 15_q 16_q 17_q 18_q 19_q}{1_q 2_q^2 3_q 5_q 6_q 7_q 8_q} \\ \\
 \frac{10_q 12_q 14_q 16_q 18_q}{1_q 3_q 5_q 7_q 9_q} \\ \\
 \frac{10_q 12_q 14_q 16_q 18_q}{1_q 3_q 5_q 7_q 9_q}
\end{array}$}
\]

With an altitude of $\kappa = 18 +1 = 19$, so that $q = \exp(i \pi / 19)$, one f\/inds that the $q$-dimensions of the fundamental irreducible representations are  1, 0, 0, 0, 0, 0, 0, 0, 1, 1. Taking into account the trivial representation, one f\/inds $\vert {\mathcal A}_1({\rm Spin}(20)) \vert = 1 + 3 = 4,$ as expected.

\subsection*{Acknowledgements}
This research was supported in part by the ANR program ``Geometry and Integrability in Mathematical Physics'', GIMP, ANR-05-BLAN-0029-0. One of us (R.C.) thanks the Centre de Recerca Matem\`atica (CRM), Bellaterra, Universitat Aut\`onoma de Barcelona, where part of this work was done, for its support.

\pdfbookmark[1]{References}{ref}
\LastPageEnding


\begin{thebibliography}{99}

\footnotesize\itemsep=0pt

\bibitem{AldazabalEtAl} Aldazabal G., Allekote I.,  Font A., Nu{\~ n}ez C.,
$N=2$ coset compactif\/ications with non-diagonal invariants, {\it Internat. J. Modern Phys. A} {\bf 7} (1992), 6273--6297, \href{http://arxiv.org/abs/hep-th/9111018}{hep-th/9111018}.

\bibitem{AltschulerBauerItzykson}   Altschuler D.,  Bauer M.,   Itzykson C.,
The branching rules of conformal embeddings,   {\it Comm.  Math.  Phys.}  {\bf 132} (1990), 349--364.

\bibitem{BaisBouwknegt}  Bais F.,  Bouwknegt P.,
A classif\/ication of subgroup truncations of the bosonic string,   {\it Nuclear  Phys.~B}  {\bf 279} (1987), 561--570.

\bibitem{Evans-I}   B\"ockenhauer J.,   Evans D.,
Modular invariants from subfactors: Type I coupling matrices and intermediate subfactors,  {\it Comm.  Math.  Phys.}   {\bf 213} (2000), 267--289,  \href{http://arxiv.org/abs/math.OA/9911239}{math.OA/9911239}.

\bibitem{Evans-II}  B\"ockenhauer J.,   Evans D.,
Modular invariants,  graphs and $\alpha$-induction for nets of subfactors.~II, {\it Comm.  Math.  Phys.} {\bf  200} (1999), 57--103, \href{http://arxiv.org/abs/hep-th/9805023}{hep-th/9805023}.

\bibitem{Evans-Kawahigashi}   B\"ockenhauer  J., Evans D.,  Kawahigashi Y.,
Chiral structure of modular invariants for subfactors, {\it Comm.  Math.  Phys.} {\bf  210} (2000), 733--784,
\href{http://arxiv.org/abs/math.QA/9907149}{math.QA/9907149}.

\bibitem{CIZ}    Cappelli A.,  Itzykson  C.,  Zuber J.-B.,
The ADE classif\/ication of minimal and $A_{1}^{(1)}$ conformal invariant theories, {\it Comm.  Math.  Phys.}  {\bf 13} (1987),  1--26.

\bibitem{Coque:Qtetra}    Coquereaux R.,
Notes on the quantum tetrahedron,  {\it Mosc. Math. J.} {\bf  2} (2002), 41--80,  \href{http://arxiv.org/abs/hep-th/0011006}{hep-th/0011006}.

\bibitem{GilCoque:ADE}    Coquereaux R., Schieber  G.,
Twisted partition functions for $ADE$ boundary conformal f\/ield theories and Ocneanu algebras of
quantum symmetries,  {\it  J. Geom.  Phys.}  {\bf 42} (2002), 216--258,  \href{http://arxiv.org/abs/hep-th/0107001}{hep-th/0107001}.

\bibitem{RobertGilSL3Categories}    Coquereaux R., Schieber G.,
Orders and dimensions for sl(2) or sl(3) module categories and boundary
conformal f\/ield theories on a torus,  {\it J.  Math.  Phys.}  {\bf 48} (2007), 043511, 17~pages, \href{http://arxiv.org/abs/math-ph/0610073}{math-ph/0610073}.

\bibitem{RobertGil:E4SU4}   Coquereaux R., Schieber G.,
From conformal embeddings to quantum symmetries: an exceptional ${\rm SU}(4)$ example,
 {\it J. Phys. Conf. Ser.} {\bf 103} (2008), 012006, 24~pages,  \href{http://arxiv.org/abs/0710.1397}{arXiv:0710.1397}.

\bibitem{YellowBook} Di Francesco  P., Matthieu P.,  S\'en\'echal   D., Conformal f\/ield theory,  {\it Graduate Texts in Contemporary Physics}, Springer-Verlag, New York, 1997.

\bibitem{DiFrancescoZuber}  Di Francesco P.,  Zuber J.-B.,
${\rm SU}(N)$ lattice integrable models associated with graphs,   {\it Nuclear  Phys. B}  {\bf 338} (1990),  602--646.

\bibitem{Drinfeld} Drinfeld V.G.,
On quasitriangular quasi-Hopf algebras and on a group that is closely connected with $\mathrm{Gal}(\overline{\mathbb Q}/\mathbb Q)$,  {\it Algebra i Analiz} {\bf 2} (1990), 149--181 (English transl.: {\it Leningrad Math. J.} {\bf 2} (1991), 829--860).

\bibitem{EtingofOnVafa} Etingof P.,
On Vafa's theorem or tensor categories, {\it Math. Res. Lett.} {\bf 9} (2002), 651--657,  \href{http://arxiv.org/abs/math.QA/0207007}{math.QA/0207007}.

\bibitem{EtingofOstrik}   Etingof  P., Ostrik  V.,
Finite tensor categories,  {\it Mosc. Math. J.}   {\bf  4} (2004), 627--654, \href{http://arxiv.org/abs/math.QA/0301027}{math.QA/0301027}.

\bibitem{Fuchs-Schellekens-Schweigert} Fuchs J.,   Schellekens B.,  Schweigert C.,
Quasi-Galois symmetries of the modular $S$-matrix,  {\it  Comm. Math. Phys.} {\bf 176} (1996), 447--465, \href{http://arxiv.org/abs/hep-th/9412009}{hep-th/9412009}.

\bibitem{FuchsRunkelSchweigert-I}  Fuchs J.,   Runkel I.,  Schweigert C.,
TFT construction of RCFT correlators. I.~Partition functions,  {\it Nuclear Phys. B} {\bf 646} (2002), 353--497, \href{http://arxiv.org/abs/hep-th/0204148}{hep-th/0204148}.

\bibitem{FuchsRunkelSchweigert-II}  Fuchs J.   Runkel I.,  Schweigert C.,
TFT construction of RCFT correlators. II.~Unoriented world sheets,  {\it Nuclear Phys. B} {\bf 678} (2003), 511--637, \href{http://arxiv.org/abs/hep-th/0306164}{hep-th/0306164}.


\bibitem{GilDahmaneHassan}   Hammaoui D.,  Schieber  G.,  Tahri  E.H.,
Higher Coxeter graphs associated to af\/f\/ine $su(3)$ modular invariants,  {\it J.  Phys. A: Math. Gen.} {\bf 38} (2005), 8259--8286, \href{http://arxiv.org/abs/hep-th/0412102}{hep-th/0412102}.

\bibitem{EstebanGil}   Isasi E.,   Schieber  G.,
From modular invariants to graphs: the modular splitting method, {\it J. Phys. A: Math. Theor.}  {\bf 40} (2007), 6513--6537, \href{http://arxiv.org/abs/math-ph/06090642}{math-ph/0609064}.

\bibitem{Kac-Peter-81}   Kac V.G.,   Peterson D.H.,
Spin and wedge representations of inf\/inite-dimensional Lie algebras and groupss, {\it Proc. Nat. Acad. Sci. U.S.A.}  {\bf 78} (1981), 3308--3312.

\bibitem{Kac-Peter}   Kac V.G.,   Peterson D.H.,
Inf\/inite-dimensional Lie algebras,  theta functions and modular forms,   {\it Adv. in Math.}  {\bf 53} (1984), 125--264.

\bibitem{KacWakimoto} Kac V.G.,  Wakimoto M.,
Modular and conformal invariance constraints in representation theory of af\/f\/ine algebras,   {\it Adv. in Math.} {\bf  70} (1988), 156--236.

\bibitem{KazhdanLusztig} Kazhdan D.,  Lusztig G.,
Tensor structures arising from af\/f\/ine Lie algebras.~III,    {\it J. Amer. Math. Soc.} {\bf 7} (1994),  335--381.

\bibitem{Kuperberg:spiders} Kuperberg G., Spiders for rank 2 Lie algebras,  {\it Comm. Math. Phys.} {\bf 180} (1996), 109--151, \href{http://arxiv.org/abs/q-alg/9712003}{q-alg/9712003}.

\bibitem{LevsteinLiberati} Levstein F., Liberati J.I.,
Branching rules for conformal embeddings, {\it  Comm. Math. Phys.} {\bf 173} (1995), 1--16.

\bibitem{Longo-Rehren} Longo R., Rehren K.-H.,
Nets of subfactors,  {\it Rev. Math. Phys.} {\bf 7} (1995), 567--598, \href{http://arxiv.org/abs/hep-th/9411077}{hep-th/9411077}.

\bibitem{Ocneanu:Unpublished}   Ocneanu  A., Seminars, 1996,  unpublished.

\bibitem{Ocneanu:paths} Ocneanu  A.,
Paths on Coxeter diagrams: from Platonic solids and singularities to minimal models and subfactors,  Notes taken by  S.~Goto,
{\it Fields Institute Monographs}, Editors Rajarama Bhat et al., AMS, 1999.


\bibitem{Ocneanu:Bariloche}  Ocneanu  A.,
The classif\/ication of  subgroups of quantum
${\rm SU}(N)$, in Proceedings of Bariloche Summer School 2000 ``Quantum symmetries
in theoretical physics and mathematics'' (January 10--21, 2000, S.C. de Bari\-loche),
Editors R.~Coquereaux, A.~Garc\'{\i}a and R.~Trinchero R., {\it  Contemp. Math.} {\bf  294} (2000), 133--160.

\bibitem{Ocneanu:MSRI} Ocneanu  A., Higher Coxeter systems,
Talk at the Workshop
``Subfactors and Algebraic Aspects of Quantum Field Theory''
(December 4--8, 2000, MSRI),\\ available at \url{http://msri.org/publications/ln/msri/2000/subfactors/ocneanu/1/index.html}.

\bibitem{Ostrik}   Ostrik  V.,
Module categories, weak Hopf algebras and modular invariants, {\it Transform.  Groups}  {\bf  8}  (2003), 177--206, \href{http://arxiv.org/abs/math.QA/0111139}{math.QA/0111139}.

\bibitem{PetkovaZuber:Oc}   Petkova V.B.,   Zuber J.-B.,
The many faces of Ocneanu cells,  {\it Nuclear  Phys. B}  {\bf 603} (2001), 449--496, \mbox{\href{http://arxiv.org/abs/hep-th/0101151}{hep-th/0101151}}.

\bibitem{PetkovaZuberNP1995}  Petkova V.B.,  Zuber J.-B.,
From CFT to graphs,  {\it Nuclear Phys. B} {\bf  463} (1996), 161--193, \href{http://arxiv.org/abs/hep-th/9510175}{hep-th/9510175}.

\bibitem{PetkovaZuberNP1997} Petkova V.B., Zuber J.-B.,
Conformal f\/ield theory and graphs,  in Proceedings 21st 
International Colloquium on Group Theoretical Methods in Physics  (July 15--20, 1996,  Goslar),
Physical Applications and Mathematical Aspects of Geometry, Groups, and Algebras,
Editors H.~D. Doebner, P.~Nattermann and W.~Scherer,
World Scientif\/ic,  Singapore, 1997, Vol.~2, 627--632,
\href{http://arxiv.org/abs/hep-th/97011032}{hep-th/9701103}.

\bibitem{Schellekens:c24} Schellekens A. N.,
Meromorphic $c = 24$ conformal f\/ield theories, {\it Comm. Math. Phys.} {\bf 153} (1993), 159--185, \href{http://arxiv.org/abs/hep-th/9205072}{hep-th/9205072}.

\bibitem{SchellekensWarner}   Schellekens A.N.,   Warner  N.P.,
Conformal subalgebras of Kac--Moody algebras, {\it  Phys. Rev. D}  {\bf 34} (1986), 3092--3096.

\bibitem{SchellekensYankielowicz} Schellekens A.N., Yankielowicz S.,
Simple currents, modular invariants and f\/ixed points, {\it Internat. J. Modern Phys. A} {\bf 5} (1990),  2903--2952.

\bibitem{Verstegen} Verstegen D.,
Conformal embeddings, rank-level duality and exceptional modular invariants, {\it Comm. Math. Phys.} {\bf 137} (1991), 567--586.

\bibitem{Wolf} Wolf J. A.,
The geometry and structure of isotropy irreducible homogeneous spaces, {\it Acta Math.} {\bf 120} (1968), 59--148, Erratum,  {\it Acta Math.} {\bf 152} (1984), 141--142.

\bibitem{Xu:1998}  Xu F.,
New braided endomorphisms from conformal inclusions, {\it Comm. Math. Phys.} {\bf 192} (1998), 349--403.

\bibitem{Xu:mirror} Xu F.,
An application of mirror extensions, \href{http://arxiv.org/abs/0710.4116}{arXiv:0710.4116}.

\end{thebibliography}
\end{document}